\documentstyle[amsfonts]{book}

\newfont{\lv}{msbm10}

\newcommand{\N}{{\lv N}}	
\newcommand{\Z}{{\lv Z}}	
\newcommand{\Q}{{\lv Q}}	
\newcommand{\R}{{\lv R}}	
\newcommand{\C}{{\lv C}}	
\newcommand{\F}{{\lv F}}	
\font\header=cmbx10 scaled\magstep1

\newlength{\praefixwidth}
\newlength{\statementwidth}
\newcommand{\blackbox}{{\hfill \vbox{\hrule height6pt width6pt}}}
\newcommand{\kasten}[1]{\hspace*{-3pt}\rule[-3pt]{6pt}{6pt}}
\newcommand{\exkasten}[1]{\hspace*{-3pt}\framebox[8pt]{\rule[-3pt]{6pt}{6pt}}}

\newcommand{\theoremEnde}{\blackbox\end{list}\vspace{3mm}}
\newcommand{\Vor}{\underline{Vor.:}\hfill\begin{minipage}[t]{\statementwidth}} 
\newcommand{\Bew}{\end{minipage}\vspace{5mm}\\ %
\begin{list}{\underline{Bew.:}}{ %
\settowidth{\labelsep}{\ \ } \settowidth{\labelwidth}{Bew.:} %
\setlength{\leftmargin}{\labelwidth} \addtolength{\leftmargin}{\labelsep} %
\setlength{\rightmargin}{0mm} \setlength{\listparindent}{0mm} %
\setlength{\partopsep}{0mm} \setlength{\topsep}{0mm}}
\item}
\newcommand{\Def}{ %
{\mbox{}\\[2mm]\header Definition\\}
\begin{list}{}{ %
\settowidth{\labelsep}{\ \ } \settowidth{\labelwidth}{Beweis :} %
\setlength{\leftmargin}{\labelwidth} \addtolength{\leftmargin}{\labelsep} %
\setlength{\rightmargin}{0mm} \setlength{\listparindent}{0mm} %
\setlength{\partopsep}{0mm} \setlength{\topsep}{0mm}}
\item}
\newcommand{\DefEnde}{\end{list}\vspace{5mm}}
\newcommand{\Bsp}{ %
{\mbox{}\\[2mm]\header Beispiel\\}
\begin{list}{}{ %
\settowidth{\labelsep}{\ \ } \settowidth{\labelwidth}{Beweis :} %
\setlength{\leftmargin}{\labelwidth} \addtolength{\leftmargin}{\labelsep} %
\setlength{\rightmargin}{0mm} \setlength{\listparindent}{0mm} %
\setlength{\partopsep}{0mm} \setlength{\topsep}{0mm}}
\item}
\newcommand{\BspEnde}{\end{list}\vspace{5mm}}
\newcommand{\Folg}{ %
{\mbox{}\\[2mm]\header Folgerung\\[3mm]}
\begin{list}{}{ %
\settowidth{\labelsep}{\ \ } \settowidth{\labelwidth}{Bew.:} %
\setlength{\leftmargin}{\labelwidth} \addtolength{\leftmargin}{\labelsep} %
\setlength{\rightmargin}{0mm} \setlength{\listparindent}{0mm} %
\setlength{\partopsep}{0mm} \setlength{\topsep}{0mm}}
\item}
\newcommand{\FolgEnde}{\end{list}\vspace{5mm}}
\newcounter{nummer}
\newcommand{\Aufzaehlung} {%
\begin{list}{\arabic{nummer}.}{\usecounter{nummer} %
\settowidth{\labelsep}{\ } \settowidth{\labelwidth}{0.} %
\setlength{\leftmargin}{\labelwidth} \addtolength{\leftmargin}{\labelsep} %
\setlength{\itemsep}{0mm} }}
\newcommand{\Aufzaehlungsende}{\end{list}}

\newcommand{\Dd}{{\bf D}}

\newcommand{\Md}{{\bf M}}
\newcommand{\Nd}{{\bf N}}

\newcommand{\rk}{\mbox{\rm rank}}

\newcommand{\mod}{\mbox{\rm mod}}
\newcommand{\Irr}{\mbox{\rm Irr}}

\newcommand{\Det}{\mbox{\rm det}}
\newcommand{\PIM}{\mbox{\rm PIM}}
\newcommand{\IBr}{\mbox{\rm IBr}}
\newcommand{\IPr}{\mbox{\rm IPr}}

\newcommand{\Ind}{\mbox{\rm Ind}}

\newcommand{\Hom}{\mbox{\rm Hom}}

\newcommand{\Gs}{{G_0}}
\newcommand{\Ks}{{K_0}}

\newenvironment{bew}{{\bf Proof.}}{\blackbox \medskip\noindent}

\font\sixdb = msbm6
\font\eightdb = msbm8

\font\twelvedb = msbm10 scaled 1200

\newfam\dbfam
\textfont\dbfam=\twelvedb
         \scriptfont\dbfam=\eightdb
                    \scriptscriptfont\dbfam=\sixdb
\def\lv{\fam\dbfam\twelvedb}
\newtheorem{definit}{Definition}[section]
\newtheorem{lem}[definit]{Lemma}
\newtheorem{prop}[definit]{Proposition}
\newtheorem{satz}[definit]{Theorem}
\newtheorem{rem}[definit]{Remark}
\newtheorem{beisp}[definit]{Example}

\newtheorem{algorithm}[definit]{Algorithm}

\newcommand{\GCD}{\mbox{\rm gcd}}
\def\address{\vskip2\baselineskip\bgroup
\parskip=0pt
\parindent=0pt
\obeylines\addrreess}
\def\addrreess#1{\noindent\ignorespaces#1\bigskip\egroup}
\begin{document}

\thispagestyle{empty}
\title{Computational Modular Character Theory}
\author{G.~HISS \and C.~JANSEN \and K.~LUX \and R.~PARKER}
\date{\vfill \ \ November 11, 1993}
\maketitle

\setcounter{page}{4}
\thispagestyle{empty}
\mbox{\ \ \rule[- 7pt]{0pt}{ 20pt}}
 
\clearpage
\thispagestyle{empty}
\cleardoublepage

\thispagestyle{empty}
\vspace*{2cm}
\begin{center}
\parbox[t]{18em}{\sl \Large To Joachim Neub\"user\\
and Herbert Pahlings}
\end{center}
 
\clearpage
\thispagestyle{empty}
 
\cleardoublepage
 
\thispagestyle{empty}

\section*{Acknowledgements}

\addcontentsline{toc}{chapter}{\textbf{Acknowledgements}}

During the development of the {\sf MOC}-system and the 
preparation of this book we were
supported by various mathematicians and institutions.

Above all, we wish to thank Joachim Neub\"user 
and Herbert Pahlings who brought us together at Aachen. 
Pahlings was the supervisor of the PhD-thesis of Lux and of the
Diploma thesis of Jansen, which were concerned with
extensions and applications of the {\sf MOC}-system.
Neub\"user and Pahlings provided a
great working atmosphere at their institute, the Lehrstuhl~D
f\"ur Mathematik, and supported our project during the whole time.
We therefore dedicate this book to them.

It is a great pleasure for us to thank H.~W.~Lenstra
for his immediate interest in a particular number theoretical
problem related to the {\sf MOC}-system, for solving this problem and
for allowing us to include his proof of this previously
unpublished result in our book.
Meinolf Geck read a preliminary version of the manuscript
with great care. His various remarks and observations concerning 
the theory of modular characters have been incorporated. Together with
his detailed comments and suggestions on the style of the exposition,
this lead to many improvements and added very much to the 
legibility of our book. We wish to thank him for all his work.
Finally we thank Thomas Breuer for 
suggestions concerning abelian number fieds. 

We also gratefully acknowledge the financial support we received from
the Deutsche Forschungsgemeinschaft and from the
Bundesland Nordrhein-Westfalen in form of the
Bennigsen-Foerder-Preis.


\clearpage
\thispagestyle{empty}

\cleardoublepage
 
\thispagestyle{empty}

\section*{Preface}

\addcontentsline{toc}{chapter}{\textbf{Preface}}

This book introduces some new ideas for the
computer calculation of modular characters of finite
groups. These ideas underly the computer algebra system
{\sf MOC}, which has been developed by the authors.
Since there have been quite a few results obtained with
the aid of {\sf MOC}, it seems appropriate to publish 
the concepts and algorithms on which it is built.
The sources of the programs of this system are freely available.
Should anyone wish to use them for any purpose, apply to any of the authors.

This book consists of six chapters and an appendix. In the first
chapter we recall and collect known results about
modular character tables of finite groups and introduce the
{\sf MOC}-system. In the 
second chapter we offer a quick introduction, mostly without
proofs, to the modular character theory of finite groups.
It is our intention to lead the beginner and non-expert quickly
to a state of being able to produce new results
with {\sf MOC}\@. In Chapter~$3$ we describe the computational
concepts underlying our system. They allow us to transform the
problems we are concerned with into problems of solving systems of
integral linear equations and inequalities. These concepts are of independent
interest as they can also be
applied in a generic way to the modular representation theory of finite groups
of Lie type. In the fourth chapter we introduce some 
principal data structures. In 
particular we suggest a new character table format in which
all entries are integers. We also present a previously 
unpublished result due to H.~W.~Lenstra on certain integral bases
for abelian number fields. 
The heart of our book is the fifth chapter, where the algorithms
are introduced. Some of them should be of independent interest,
for example our algorithms for dealing with rational linear
equations. Finally, in Chapter~$6$ we have worked out a non-trivial
example in detail, namely the calculation of the $7$-modular decomposition
numbers of the double covering group of Conway's largest group.

\pagebreak

In the appendix we have collected some results obtained in the
course of examples throughout the book. These results include
the $5$-modular decomposition numbers of the Conway group
$\mbox{\em Co}_2$, including the Brauer character table
and the $7$-modular decomposition numbers of $2\mbox{\em Co}_1$\@.
Both results are published for the first time.

Finally we give a bibliography, which is not intended to be exhaustive,
of papers and monographs cited in our book and also of some related work. 

We hope that even the expert in the field of modular representation
theory will find our book interesting, because of 
our different point of view---we must compute with the objects.

\bigskip

\noindent \parbox{5.5cm}{\em Aachen, Heidelberg, Shepton Mallet \\
November 1993} \hfill
\parbox{4.1cm}{G.~H.,  C.~J.,  K.~L.,  R.~P.}
 
 
\clearpage
\thispagestyle{empty}

\cleardoublepage
 
\markboth{LIST OF NOTATIONS}
         {LIST OF NOTATIONS}

\thispagestyle{empty}
\section*{List of notations}

\addcontentsline{toc}{chapter}{\textbf{List of notations}}

\newlength{\abstand}
\setlength{\abstand}{\parsep}
\begin{center} {\bf  1.\ Numbers} \end{center}

\begin{list}{}{\setlength{\leftmargin}{30mm}
               \setlength{\labelwidth}{27mm}
               \setlength{\labelsep}{3mm}
               \setlength{\itemsep}{0pt}
               \setlength{\parsep}{2pt}}
\item [$|M|$ \hfill{\ }] Number of elements of the set~$M$
\item  [$p$ \hfill{\ }] A prime number
\item [$\GCD(n,m)$ \hfill{\ }] Greatest common divisor of
      the integers~$n$ und~$m$
\item [$n_p$,~$n_{p'}$ \hfill{\ }] $p$-part respectively $p$-free part of the integer~$n$
\item [$\Z$ \hfill{\ }] Integers
\item [$\N$ \hfill{\ }] Positive integers
\item [$\N_0$ \hfill{\ }] $\N \cup \{ 0 \}$
\item [$\Q$ \hfill{\ }] Rational numbers
\item [$\Q_n$ \hfill{\ }] $\Q(\zeta_n)$, where~$\zeta_n$ denotes a
      primitive $n$-th roots of unity
\item [$\R$ \hfill{\ }] Real numbers
\item [$\C$ \hfill{\ }] Complex numbers
\item [$\bar{z}$ \hfill{\ }] Complex conjugate of $z \in \C$
\item [$\F_q$ \hfill{\ }] Finite field with $q$ elements
\item [${\cal G}(L/K)$ \hfill{\ }] Galois group of Galois field
      extension $K \subseteq L$

\end{list}

\begin{center} {\bf  2.\ Matrices} \end{center}
\begin{list}{}{\setlength{\leftmargin}{30mm}
               \setlength{\labelwidth}{27mm}
               \setlength{\labelsep}{3mm}
               \setlength{\itemsep}{0pt}
               \setlength{\parsep}{2pt}}
\item [$X^t$ \hfill{\ } \hfill{\ }] Transpose of matrix~$X$
\item [$\Det(X)$ \hfill{\ } \hfill{\ }] Determinant of~$X$
\item [$E_n$ \hfill{\ } \hfill{\ }] $n \times n$-Identity matrix
\end{list}

\begin{center} {\bf  3.\ Groups} \end{center}
\begin{list}{}{\setlength{\leftmargin}{30mm}
               \setlength{\labelwidth}{27mm}
               \setlength{\labelsep}{3mm}
               \setlength{\itemsep}{0pt}
               \setlength{\parsep}{2pt}}
\item [$G$ \hfill{\ }] Finite group
\item [$H \leq G$ \hfill{\ }] $H$ is a subgroup of~$G$
\item [$C_G(X)$ \hfill{\ }] Centralizer of $X \subseteq G$
\item [$N_G(H)$ \hfill{\ }] Normalizer of subgroup~$H$ in~$G$
\item [$G_p$, $G_{p'}$ \hfill{\ }] Set of $p$- respectively~$p'$-elements 
      of~$G$
\item [$|x|$ \hfill{\ }] Order of~$x \in G$
\end{list}

\pagebreak

\begin{center} {\bf  6.\ Rings and modules} \end{center}
\begin{list}{}{\setlength{\leftmargin}{30mm}
               \setlength{\labelwidth}{27mm}
               \setlength{\labelsep}{3mm}
               \setlength{\itemsep}{0pt}
               \setlength{\parsep}{2pt}}
\item [$R$ \hfill{\ }] Commutative ring with~$1$
\item [$R^*$ \hfill{\ }] Group of units of~$R$
\item [$\rk_R(X)$ \hfill{\ }] $R$-rank of the free $R$-module~$X$
\item [$RG$ \hfill{\ }] Group algebra of~$G$ over~$R$
\item [$\Hom_{RG}(X,Y)$ \hfill{\ }] Set of ${RG}$-homomorphisms between the
$RG$-modules $X$ and~$Y$
\item [${[X,Y]}$ \hfill{\ }] $\rk_R(\Hom_{RG}(X,Y))$ for $RG$-modules~$X$ and~$Y$
\item [$X^G$ \hfill{\ }] Induced module
\item [$X_H$ \hfill{\ }] Restricted module
\item [$X^*$ \hfill{\ }] $RG$-module dual to $X$
\item [$(K,R,F)$ \hfill{\ }] $p$-modular system:
                  \begin{description}
                  \item $R$: Ring of algebraic integers
                  \item $K$: Field of fractions of~$R$
                  \item $F$: Residue class field of~$R$ of characteristic $p \neq 0$
                  \end{description}
\end{list}

\begin{center} {\bf  8.\ Characters} \end{center}
\begin{list}{}{\setlength{\leftmargin}{30mm}
               \setlength{\labelwidth}{27mm}
               \setlength{\labelsep}{3mm}
               \setlength{\itemsep}{0pt}
               \setlength{\parsep}{2pt}}
\item [$\vartheta_{\cal X}$ \hfill{\ }] Brauer character afforded by
      the representation~${\cal X}$
\item [$\Phi_{\varphi}$ \hfill{\ }] Projective indecomposable character
      corresponding to the irreducible Brauer character~$\varphi$
\item [$\PIM$ \hfill{\ }] Projective indecomposable character
\item [$\vartheta^{2+}$ \hfill{\ }] Symmetric square of character~$\vartheta$
\item [$\vartheta^{2-}$ \hfill{\ }] Skew square of character~$\vartheta$
\item [$\omega_\chi$ \hfill{\ }] Central character corresponding to
      $\chi \in \Irr(G)$
\item [$\Irr(G)$ \hfill{\ }] Set of irreducible $K$-characters of~$G$
\item [$\IBr(G)$ \hfill{\ }] Set of irreducible Brauer characters of~$G$
      with respect to $(K,R,F)$
\item [$\IPr(G)$ \hfill{\ }] Set of projective indecomposable characters of~$G$
\item [$B$ \hfill{\ }] A union of blocks of~$G$
\item [$\Irr(B)$ \hfill{\ }] Set of irreducible $K$-characters of~$G$ 
      belonging to~$B$
\item [$\IBr(B)$ \hfill{\ }] Set of irreducible Brauer characters of~$G$      with respect to $(K,R,F)$ 
      belonging to~$B$
\item [$\IPr(B)$ \hfill{\ }] Set of projective indecomposable characters of~$G$ 
      belonging to~$B$
\item [$\Gs(KG)$ \hfill{\ }] Set of all $\Z$-linear combinations of $\Irr(G)$
\item [$\Gs(FG)$ \hfill{\ }] Set of all $\Z$-linear combinations of $\IBr(G)$
\item [$\Ks(FG)$ \hfill{\ }] Set of all $\Z$-linear combinations of $\IPr(G)$
\item [$\Gs(B)$ \hfill{\ }] Set of all $\Z$-linear combinations of $\Irr(B)$
\item [$\Gs(B)^+$ \hfill{\ }] Set of all $\N$-linear combinations of $\Irr(B)$
\item [$\Gs(\hat{B})$ \hfill{\ }] Set of all $\Z$-linear combinations of $\IBr(B)$ 
\item [$\Gs(\hat{B})^+$ \hfill{\ }] Set of all $\N$-linear combinations of $\IBr(B)$ 
\item [$\Ks(B)$ \hfill{\ }] Set of all $\Z$-linear combinations of $\IPr(B)$
\item [${\cal C}(G,B)$ \hfill{\ }] Set of class functions of~$G$ with values in~$K$, spanned by $\Irr(B)$
\item [${\cal C}_{p'}(G)$ \hfill{\ }] Set of class functions of~$G$ 
      which vanish on $p$-singular elements
\item [${\cal C}_{p'}(G,B)$ \hfill{\ }] Set of class functions of~$G$ 
with values in~$K$, spanned by $\IBr(B)$
\item [$\psi_B$ \hfill{\ }] Projection of class function
      $\psi$ to ${\cal C}(G,B)$ respectively ${\cal C}_{p'}(G,B)$
\item [${\cal C}(G)$ \hfill{\ }] Set of class functions of~$G$ with values in~$K$
\item [${[M]}$ \hfill{\ }] Matrix of values of $M \subset {\cal C}(G)$
\item [$\chi_H$ \hfill{\ }] Class function restricted to  subgroup~$H$
\item [$\chi^G$ \hfill{\ }] Class function induced to~$G$
\item [$\bar{\chi}$ \hfill{\ }] Class function complex conjugate to~$\chi$
\item [$\bar{M}$ \hfill{\ }] Set of complex conjugates of $M \subset {\cal C}(G)$
\item [${\langle \mbox{\ }, \mbox{\ } \rangle}$ \hfill{\ }] 
      Inner product of ${\cal C}(G)$
\item [$\langle M , N \rangle$ \hfill{\ }] Matrix of inner products 
      between the elements of the sets $M, N \subset {\cal C}(G)$ 
\item [$\hat{\psi}$ \hfill{\ }] Restriction of class function~$\psi$ to the $p'$-elements\item [$\widehat{\Irr}(B)$ \hfill{\ }] $\{\hat{\chi} \mid \chi \in \Irr(B) \}$   
\item [${\bf D}_B$ \hfill{\ }] Decomposition matrix of~$B$
\item [$d_{\chi\varphi}$ \hfill{\ }] Decomposition number
\item [${\bf BS}$ \hfill{\ }] Basic set of Brauer characters
\item [${\bf PS}$ \hfill{\ }] Basic set of projective characters
\item [${\bf BS}_0$ \hfill{\ }] Special basic set of Brauer characters
\item [${\bf BA}$ \hfill{\ }] Basis of Brauer atoms dual (with respect to
      $\langle \mbox{\ }, \mbox{\ } \rangle$) to ${\bf PS}$
\item [${\bf PA}$ \hfill{\ }] Basis of projective atoms dual to ${\bf BS}$
\item [${\bf PA}_0$ \hfill{\ }] Basis of projective atoms dual to ${\bf BS}_0$
\end{list}

\setlength{\parsep}{\abstand}

\clearpage
\thispagestyle{empty}

\cleardoublepage  

\tableofcontents
 
\clearpage
\thispagestyle{empty}
\cleardoublepage

\setcounter{page}{1}
\pagenumbering{arabic}

\setcounter{page}{1}
\chapter{Introduction}

\section{Modular characters}

Representation theory is an effective and indispensable
tool in the investigation of finite groups. The 
table of ordinary irreducible characters of a finite group~$G$
encodes an enormous amount of information about~$G$\@.
For example, it can be used in many cases to show that~$G$
is a Galois group over some abelian extension of the rational
numbers \cite{matzatK}\@. 
The beautiful book of Isaacs~\cite{isa} gives an impression
of the innumerable applications of character theory.

Ordinary characters for arbitrary finite groups where introduced by
Frobenius in generalization of characters of abelian groups.
Frobenius himself already calculated the character tables
of the symmetric and alternating groups 
and of the simple Mathieu groups. 
To accomplish this work he was led to invent the
method of inducing characters \cite{frobenius}\@. 
Schur not much later gave another 
proof for the character tables of the symmetric and alternating groups,
and continued his investigations to complete the character tables
of their covering groups \cite{schur}\@.
The ordinary character tables of the finite simple groups of Lie type are
now almost known by the work of Deligne and Lusztig \cite{luszbuch}\@. 
The famous ``Atlas of Finite Groups'' \cite{atlas} collects,
among a wealth of other valuable information,
the character tables of the most interesting finite simple groups
up to a certain order, in particular those for the sporadic simple groups.

Brauer started the investigation of representations
of a finite group over a field whose characteristic divides 
the group order. In \cite{brauernesbittI}, he introduced the idea of
modular characters, now called Brauer characters,
and the idea of decomposition numbers. The knowledge of
decomposition numbers is equivalent to the knowledge of the
irreducible Brauer characters. One of the main
tasks of modular representation theory is to find 
methods  for the calculation of the decomposition
numbers for a given finite group.

Let us mention a few of the various applications 
of modular representation theory. 
To begin with, it is useful for the determination of the 
maximal subgroups of a group. For example, building upon
Aschbacher's subgroup theorem \cite{as1} and some results
on modular representations of finite simple groups, Kleidman and
Liebeck were able to determine the maximal subgroups of most
of the finite classical groups in \cite{kleidlieb}\@.
There are also applications outside the theory of finite groups,
for example in geometry and combinatorics, where finite geometries
are studied via their automorphism groups (see for example \cite{ronsmi})
or in topology, where topological spaces are investigated
via their fundamental groups (see for example \cite{hissszczep})\@.
In many potential applications the problem can be
reduced to the case of the simple groups, their
automorphism groups and covering groups.
These groups are known by the classification of the finite simple groups.
It would therefore be very desirable to have a description
of all irreducible Brauer characters for this class of groups.
Such a large collection of examples would doubtlessly lead to new
conjectures and theorems in representation theory and then in turn 
perhaps to new and simpler ways of obtaining these Brauer characters.
As Michler has indicated in his survey article~\cite{Michler}, 
such a collection of modular character tables is also extremly 
useful in checking various famous conjectures of Brauer and Alperin.

\section{Known results}

Let us now describe some of the history of the Brauer character
tables which are known up to the present. 
Brauer and Nesbitt in \cite{brauernesbittI}
already calculated the modular irreducible characters
for the groups {\em PSL}$(2,q)$ in the defining characteristic case, i.e.,
for fields of characteristic $p$ dividing~$q$\@. 
Janko in \cite{j1} found the Brauer characters for his 
first sporadic group $J_1$ in every odd characteristic
and used one of the $11$-modular irreducible characters
to construct his group as a matrix group
of degree~$7$ over a field of characteristic~$11$\@.
Fong in~\cite{Fon} determined the $2$-decomposition numbers 
for this group with some ambiguities.
In~\cite{Jam} James calculated all but one of the irreducible Brauer characters
for the Mathieu groups, the only exception being one 
character in characteristic~$2$ for $M_{24}$\@. This last character
was found only much later with the help of the {\sf Meat-Axe} \cite{MeatAxe}\@.
Humphreys in \cite{Huma} and Benson in \cite{benm22}
investigated the covering groups of the Mathieu groups.
The Higman--Sims group has been studied by Humphreys in \cite{Hum},
and Thackray in his thesis~\cite{thackray} determined the
$2$-modular characters of the McLaughlin group.
In \cite{fblock} Feit introduced some new methods by 
utilizing Green correspondence
in the determination of decomposition numbers, thereby 
calculating many substantial examples, including some for the
triple covering group of the McLaughlin group and the third
Janko group~$J_3$\@.
Feit's methods have been extended and led to the results
of~\cite{BBB}\@. 
By the combined efforts of various authors 
the decomposition numbers are now completely
known for all groups appearing in the Atlas up to page~$100$,
the largest of which is the McLaughlin group.
These results are collected in the first part of a modular Atlas
of finite groups  \cite{modATL}\@.

\section{\sf MOC}

The computer algebra system {\sf MOC} (for MOdular
Characters)  which we are going to 
describe in this book has been developed by the authors in
order to automate some of the more elementary methods known for
calculating decomposition numbers. 
Parts of the material presented here has already been sketched in
\cite{lupa}\@.
An account of the methods we have in mind
has been given by James and Kerber in~\cite[\S 6.3]{jake}\@.
In a sense, modular character theory has been a part of
``computational algebra'' right from the beginning, and most
of the elementary methods of [loc.~cit.] have already
been described by Brauer and Nesbitt in \cite{brauernesbittI}\@.
Many of the algorithms of {\sf MOC} which are mostly concerned with systems
of integral equations or inequalities are not new,
but do not seem to be widely known or used. A subsidiary purpose
of our book is therefore to show how powerful these methods can be, 
particularly the methods of solving systems
of integral linear equations by $p$-adic expansion and
the Gomory-algorithm for solving integral linear inequalities. 

Some of the concepts developed for
{\sf MOC} have turned out to be of theoretical interest
in the investigation of other problems, too. 
The work of Fong--Srinivasan and Dipper--James has
led Geck and one of the authors
to the formulation of a theorem on basic sets of groups
of Lie type \cite{geckhiss}\@. The concept of basic set as 
used in this paper was introduced in {\sf MOC} in order to 
mechanize the calculation of decomposition numbers for individual
finite groups. It is now clear from our experience that
many concepts and algorithms can be used also for
series of groups of Lie type, since the methods work
for generic character tables not just for individual groups.

Some other aspects of our computational work are worth mentioning.
First of all, our need for simple data structures
for calculation in cyclotomic fields led
to a conjecture that abelian number fields always have a nice
integral basis consisting of sums (over an orbit) of roots of unity. 
This conjecture has now been proved by H.~W.~Lenstra.
We are grateful to him for allowing us to 
present his proof in our book. These special bases for abelian
number fields are now also used in the Aachen 
{\sf GAP} system \cite{GAP}\@.

The principal design of {\sf MOC} is at present different 
from other computer algebra systems available. {\sf MOC} consists
of a collection of small one-purpose FORTRAN-programs. We have,
for example one program for matrix multiplication and one for
solving integral linear equations. The various
programs communicate via files. 
More complex tasks are achieved by using Bourne-shell scripts
calling a sequence of {\sf MOC} programs.
The shell scripts are an essential part
of the {\sf MOC} system. One advantage of this design
is the simple extendability of the system. 
It allows the easy addition of new commands to {\sf MOC}\@.
A further advantage is the portability to other UNIX machines. 
Another important feature of {\sf MOC} is supported by this design. 
Every run of a program is recorded in
a file which is attached to a particular group and prime.
The {\sf MOC} programs also write their output
onto this log-file automatically. This makes it possible to
repeat particular calculations and, what is most important,
to write down proofs for the results obtained by the machine.
Even the proof writing can be done automatically to a certain extent.
This feature of {\sf MOC} has been used to obtain most of the
proofs of~\cite{BBB}\@.

We emphasize the fact that {\sf MOC} only works with Brauer characters
and not with representations (c.f.~the {\sf Meat-Axe} \cite{MeatAxe})
and uses only elementary methods
to the extent described in \cite[\S~6.3]{jake}\@.
Nevertheless, in connection with the more advanced methods
described in Section~\ref{AdvancedMethods},
one can get results which are out of reach
of systems working with representations because of
the large degrees of the representations involved. On the other hand, 
for some groups, the {\sf Meat-Axe} and
{\sf Condensation} (see \cite{rybacon,peak})
are used successfully to solve problems
left open by {\sf MOC}\@. This applies in particular for
calculations in characteristics~$2$ and~$3$,
where elementary methods are much less efficient than 
for larger primes. The three systems working
together provide a very powerful computational tool for the
calculation of irreducible Brauer characters for individual
finite groups.

\section{Series of groups}

Apart from individual groups, there is also considerable knowledge
on some series of simple groups. For the reason of completeness
we briefly sketch the known results. 
Various authors, among them Robinson, Kerber and James have established
a theory of modular representations for the symmetric and alternating groups
(see for example \cite{jamlec,Kerber})\@.
Burkhardt did some work for groups of Lie type; in particular he found
all decomposition numbers for the simple groups $PSL(2,q)$ \cite{bursl},
all decomposition numbers for the Suzuki groups in odd 
characteristics \cite{bursu}, and some decomposition numbers
for the Suzuki groups and the unitary groups $U(3,2^f)$ in even
characteristic \cite{bursu,burpsu}.
Fong in \cite{Fon} started the investigation for the small
Ree groups ${^2G}_2(3^{2m+1})$, but he had to leave some problems
open. These were later solved by Landrock and Michler \cite{lami}\@. 
The work on these groups was completed
in \cite{reetree}\@.
 
The investigation for the groups of Lie type naturally divides into
the case where the modular characteristic is the defining characteristic
of the group, i.e., the characteristic of the underlying field,
and the case where the two
characteristics are different.
In the defining characteristic case the representation theory of
algebraic groups is relevant, in particular the problem of finding
the multiplicities of the simple modules in the Weyl modules.
An affirmative answer
to Lusztig's conjecture (see \cite{luszconj})
would solve most of the problems for
the finite groups of Lie type in this case.
 
In the non-defining characteristic case, a theory has been
evolving starting with a series of papers 
by Fong and Srinivasan.
The highlight of their work is the complete description of all the
Brauer trees for the classical groups in \cite{fongcl}\@. Dipper and James have
been developing a representation theory for the general
linear groups which has led to the completion of the
decomposition numbers for $GL_n(q)$, $n \leq 10$ \cite{glndec}\@.
Furthermore, some of the decomposition numbers for
$G_2(q)$ and $Sp_4(q)$ are known by the work of Hiss,
Shamash and White. 
Geck has determined most of the irreducible
Brauer characters for $U(3,q)$ and the
Steinberg triality groups ${^3D(4,q^3)}$ in \cite{geck,gecktrial}\@.
Finally, the Brauer trees for the Ree groups are known 
completely~\cite{reetree}\@.
 
\section{Special defect groups}

In general the problem of finding the Brauer character table of a
finite group is considerable harder than that of finding the ordinary
irreducible characters. This is partly due to the fact that
in contrast to the ordinary character theory there is no inner
product on the set of Brauer characters which tells the
irreducibility of a Brauer character. However, in particular
cases, there are deep theoretical results which can be used
to good advantage. If a group is $p$-soluble, the
Fong--Swan theorem states that every irreducible
Brauer character is liftable to an ordinary character.
This leads immediately to an algorithm
for the determination of the irreducible Brauer characters.
 
Brauer has introduced the concepts of blocks, which leads to
a subdivision of the problems. The blocks correspond to the
primitive idempotents of the centre of the modular group algebra.
The decomposition matrix is a matrix direct sum of the decomposition
matrices for the block algebras, and so it suffices to consider
each block in turn. Brauer has associated to each block a conjugacy
class of $p$-groups, the defect groups of the block. Their isomorphism
type measures
the difficulty or complexity of finding the Brauer characters for
the block. If the defect groups are trivial, the block is called
a block of defect zero, having as decomposition matrix a $(1 \times 1)$-matrix
with single entry~$1$\@. If the defect groups are cyclic, the
Brauer--Dade theorem imposes
strong restrictions on the decomposition numbers.
In this case, the decomposition matrix can be encoded in a certain
tree, the Brauer tree.
There are other defect groups of a particular type, namely
dihedral $2$-groups and quaternion groups, where all possible
decomposition matrices occuring for blocks with these defect groups
are known by work of Brauer, Olsson and Erdmann \cite{erdmann}\@.

Brou\'e and Puig have introduced the notion of nilpotent
blocks in~\cite{broupui}, 
generalizing the concept of blocks with cyclic defect groups.
Puig \cite{nilpuig} has given a remarkable extension of the Brauer--Dade 
theorem.
His results, which are neatly explained in \cite{kuels},
imply the knowledge of the decomposition numbers for nilpotent blocks.

It has been conjectured by Donovan
that, given a finite $p$-group~$D$,
there are only finitely many matrices which are decomposition
matrices for blocks of finite groups with defect groups
isomorphic to~$D$\@.
Scopes \cite{scopes} has shown the conjecture to be true
if we restrict the class of groups to symmetric groups.
For cyclic defect groups Donovan's conjecture is true by 
the results of Brauer and Dade. Nevertheless, to write down
all decomposition matrices belonging to blocks with a given
cyclic defect group, one still uses the classification of finite
simple groups (c.f.~\cite{ftree})\@. It is our opinion that
finite simple groups and in particular the sporadic groups will
play an important role in the description of decomposition 
matrices for general defect groups, even if Donovan's conjecture
is shown to be true some day. This, we believe, justifies the
application of our methods and in particular the usage of {\sf MOC}\@.

\clearpage
\thispagestyle{empty}

\chapter{Background from modular character theory}
\markboth{MODULAR CHARACTER THEORY}{BRAUER CHARACTERS}
\label{modchartheory}

In this chapter we introduce our notation and recall some basic facts
of character theory. We follow the exposition of Goldschmidt
\cite{gold} and also Chapter~$15$ of Isaacs' book \cite{isa}\@.
This chapter is intended to give the beginner a quick introduction
to those parts of the theory which are essential for the understanding
of the following. We give no proofs but references to the
appropriate sources. We assume the reader to be familiar with
ordinary character theory of finite groups.

Throughout this chapter let~$G$ be a finite group, and~$p$ a fixed
prime number. 

\section{Brauer characters}
\markboth{MODULAR CHARACTER THEORY}{BRAUER CHARACTERS}

We use almost the same set-up as Isaacs in Chapter~$15$ of~\cite{isa}\@.
Let~$R$ denote the ring of algebraic integers over the rationals,
and let~$K$ be the quotient field of~$R$\@.
Choose a maximal ideal~$I$ of~$R$ containing~$p$\@.
Finally let $F = R/I$ denote the residue class field of~$R$,
a field of characteristic~$p$ which is algebraically closed\@. 
The canonical epimorphism $R \rightarrow F$
induces an isomorphism~$\rho$
from the set of $|G|_{p'}$-th roots of unity to~$F^*$
(see \cite[Lemma (15.1)]{isa}), 
which is used to define Brauer characters as follows.

An element of~$G$ of order not divisible by~$p$ is called
{\sl $p$-regular} or a {\sl $p'$-element}. An element of order divisible
\index{$p$-regular element}
\index{$p'$-element}
by~$p$ is called {\sl $p$-singular}. The set of all $p$-regular
elements is denoted by $G_{p'}$\@. 
\index{$p$-singular element}
\begin{definit}
{\rm
Let ${\cal X}: G \rightarrow 
\mbox{\em GL}_n(F)$ be an $F$-representation of~$G$\@.
The {\sl Brauer character}~$\vartheta_{\cal X}$ of~${\cal X}$
\index{Brauer character}
is a map $\vartheta_{\cal X} : G_{p'} \rightarrow K$\@. Let $g \in G_{p'}$,
and let $\zeta_1, \ldots , \zeta_n$ denote the eigenvalues of
${\cal X}(g)$\@. We then define
$\vartheta_{\cal X}(g) = \sum_{i=1}^n \rho^{-1}(\zeta_i)$\@.
}
\end{definit}

\noindent If we apply the canonical epimorphism to $\vartheta_{\cal X}(g)$,
we obtain the trace of the matrix ${\cal X}(g)$\@.
We say that the Brauer character $\vartheta_{\cal X}$ is 
{\sl afforded} by the representation ${\cal X}$\@. The function
complex conjugate to $\vartheta_{\cal X}$ is also a Brauer
character. It is afforded by the representation dual to ${\cal X}$\@.
We sometimes think of a Brauer character as being defined
on all of~$G$ by letting its values be~$0$ on $p$-singular
elements. This point of view is taken primarily in this
and the next chapter in order to simplify notation. In the 
implementation of the {\sf MOC}-system a Brauer character is
of course only stored as a function on $G_{p'}$\@.

A Brauer character is called irreducible if it corresponds to
an irreducible $F$-re\-pre\-sent\-ation of~$G$\@. It is a basic
fact that the set of irreducible Brauer characters of~$G$
is finite and linearly independent as functions from~$G$ to~$K$
(\cite[Theorem 15.5]{isa})\@. 

Note that there are several maximal ideals~$I$ inside~$R$
containing~$p$\@.  A different choice of~$I$
would in general give another  Brauer character
for the same $F$-representation. However, all {\sf MOC}-calculations are
independent of the choice of~$I$\@. The reader is invited to
check this assertion, as we go along.

\section{The decomposition matrix}
\markboth{MODULAR CHARACTER THEORY}{THE DECOMPOSITION MATRIX}

The rings $(K,R,F)$ are fixed in the following. They form
our {\sl $p$-modular system}\@.
\index{$p$-modular system}
 
The set of ordinary irreducible $K$-characters of~$G$ is denoted by
$\Irr(G)$, the set of irreducible Brauer characters 
by $\IBr(G)$\@.
Since~$K$ is a splitting field for the group ring~$KG$, 
the number of irreducible 
$K$-characters equals the number of conjugacy classes of~$G$\@.

The set of all $\Z$-linear combinations of ordinary characters
is a free abelian group with basis $\Irr(G)$\@. It is isomorphic to
the Grothendieck group of the category of finitely generated $KG$-modules 
(see \cite[Proposition (16.10)]{CRI}) and is therefore
denoted by $\Gs(KG)$\@. The latter fact is not needed in our book;
it is only stated to motivate the notation.
Similarly, $\Gs(FG)$ denotes the set of all $\Z$-linear 
combinations of irreducible Brauer characters.
Elements in these groups are called
{\sl generalized} or {\sl virtual characters}.
\index{generalized characters}
\index{virtual characters}
A function on~$G$ which is constant on conjugacy classes is called 
a {\sl class function}.
\index{class function}
Characters are class functions and thus $\Gs(KG)$ and
$\Gs(FG)$ are subgroups of ${\cal C}(G)$, the set of $K$-valued
class functions of~$G$\@.
We let ${\cal C}_{p'}(G)$ denote the set of class functions of~$G$
with values in~$K$, vanishing on $p$-singular elements.
Then, by our definition, $\Gs(FG)$ is contained in
${\cal C}_{p'}(G)$\@. 

If $\chi \in \Gs(KG)$, we define~$\hat{\chi}$ by
\begin{equation}
\label{pprimerestriction}
\hat{\chi}(g) = \left\{ 
\begin{array}{cc}
\chi(g), & \mbox{if $g$ is $p$-regular,} \\
0, & \mbox{otherwise}.
\end{array}\right. 
\end{equation}
Now comes the important observation.
If~$\chi$ is an ordinary character, $\hat{\chi}$ is a Brauer character.
This can be seen as follows. Let ${\cal X}$ be a $K$-representation
of~$G$ affording the character~$\chi$\@. Let~$\tilde{R}$ denote 
the localization of~$R$ at the maximal ideal~$I$\@. Then ${\cal X}$
can be realized over~$\tilde{R}$, i.e., there is some
representation ${\cal Y}$ equivalent to ${\cal X}$, such that
the matrices representing the group elements have all entries
in~$\tilde{R}$ (\cite[Theorem 15.8]{isa})\@. We remark that  
${\cal X}$ can even be realized over~$R$ (\cite[Satz V.12.5]{HupI}), 
but this fact is less elementary.
The canonical epimorphism $R \rightarrow F$ can be extended uniquely
to an epimorphism $\tilde{R} \rightarrow F$\@. We thus obtain
a representation $\hat{\cal Y}: G \rightarrow \mbox{\em GL}_n(F)$
by applying this epimorphism to every entry of a representing
matrix. Then the Brauer character afforded by $\hat{\cal Y}$ obviously
is~$\hat{\chi}$\@. 

We now consider the decomposition homomorphism:
$$\begin{array}{rcl} d: \Gs(KG) & \longrightarrow & \Gs(FG) \\
                     \chi    & \longmapsto     & \hat{\chi} \end{array}$$
The matrix of~$d$ with respect
to the bases $\Irr(G)$ and $\IBr(G)$ is
called the {\sl decomposition matrix} of~$G$ and denoted by~$\Dd$\@.
\index{decomposition matrix}
The entries of~$\Dd$ are the non-negative integers~$d_{\chi\varphi}$
defined by 
$$\hat{\chi} = \sum_{\varphi \in {\rm IBr}(G)} d_{\chi\varphi}\,\varphi.$$
The purpose of {\sf MOC} is to support the calculation of
decomposition matrices. 

It is an important fact, which follows from Brauer's characterization
of characters, that the decomposition homomorphism is surjective
(\cite[Theorem 15.14]{isa})\@. This means that every irreducible
Brauer character is a $\Z$-linear combination of 
$\{\hat{\chi} \mid \chi \in \Irr(G)\}$\@. The fact that $\Irr(G)$
is a basis of ${\cal C}(G)$ now implies that the set of
irreducible Brauer characters spans ${\cal C}_{p'}(G)$\@. As
$\IBr(G)$ is linearly independent, we obtain that $|\IBr(G)|$
equals the number of $p$-regular conjugacy classes of~$G$\@.

\section{Projective characters and orthogonality relations}
\markboth{MODULAR CHARACTER THEORY}{PROJECTIVE CHARACTERS AND ORTHOGONALITY}

We are now going to introduce the projective characters.
\begin{definit}
{\rm
Let $\varphi \in \IBr(G)$\@. The {\sl projective indecomposable
character} corresponding to~$\varphi$
\index{projective indecomposable character}
is the ordinary character $\Phi_\varphi$ defined by
$$\Phi_\varphi = \sum_{\chi \in {\rm Irr}(G)} d_{\chi\varphi}\,\chi.$$
A projective indecomposable character will be called a PIM in the 
following.  We denote by $\IPr(G)$ the set of all PIMs and by
$\Ks(FG)$ the group of all $\Z$-linear combinations of $\IPr(G)$\@.
A $\Z$-linear combination with non-negative coefficients of PIMs
is called a {\sl projective character}, or a {\sl genuine projective character}.
\index{projective character}
}
\end{definit}

\noindent It is now most important to observe that a PIM vanishes on
$p$-singular classes (\cite[(6.9)(a)]{gold})\@. In particular,
$\Ks(FG) \subseteq {\cal C}_{p'}(G)$\@. 

This definition of PIMs is suggested by Brauer's reciprocity law. 
It can be given without introducing projective modules
and some important properties of projective characters can be derived
on this more elementary level. Of course, to get deeper results on
projective characters one has to give the usual definition
via projective modules.

For $z \in \C$ let $\bar{z}$ denote the number complex conjugate
to~$z$\@. The usual inner product on ${\cal C}(G)$ 
is denoted by $\langle \mbox{\ } , \mbox{\ } \rangle$\@.
It is defined by
\begin{equation}
\label{iK} 
   \langle \lambda, \mu \rangle = \frac{1}{|G|} \sum_{g \in G} 
   \lambda(g) \overline{\mu(g)}.
\end{equation}
We recall the {\sl orthogonality relations}. 
\index{orthogonality relations}

\begin{satz}
\label{orthorel}
{\rm (Orthogonality relations for Brauer characters):} \\
The two sets $\IPr(G)$ and $\IBr(G)$ are dual bases of 
${\cal C}_{p'}(G)$ with respect to $\langle \ , \ \rangle$\@.
\end{satz}
\begin{bew}
This follows from the orthogonality relations for ordinary characters
and the fact that a PIM vanishes on $p$-singular classes
(see \cite[(6.10)]{gold})\@.
\end{bew}

\noindent It follows from Theorem~\ref{orthorel} that
$\langle \chi, \psi \rangle \in \Z$, if 
$\chi \in \Ks(FG)$ and $\psi \in \Gs(FG)$\@. 
Furthermore, $\langle \chi, \psi \rangle \in \N$, if
$\chi$ is a projective character and $\psi$ a Brauer character.
The latter observation has the following converse,
which will be used in the construction of projective
characters.
\begin{prop}
\label{projconstruct}
Let~$\Psi$ be an ordinary character of~$G$ which vanishes on
$p$-singular classes. Then $\Psi \in \Ks(FG)$\@.
If, furthermore, 
$$\langle \varphi, \Psi \rangle \geq 0 \quad\quad \mbox{\rm for all }
\quad\quad \varphi \in \IBr(G),$$
then $\Psi$ is a genuine projective character.
\end{prop}
\begin{bew}
By the orthogonality relations, $\IPr(G)$ is a $K$-basis of
${\cal C}_{p'}(G)$\@. Write
$$\Psi = \sum_{\varphi \in {\rm IBr}(G)} a_\varphi\,\Phi_\varphi.$$
Then $a_\varphi = \langle \varphi, \Psi \rangle$, which is an integer
since $\varphi$ is a $\Z$-linear combination of restricted
ordinary characters.
\end{bew}

\noindent This proposition is used later on
to show that certain constructions
yield genuine projective characters. The fundamental principle
of the {\sf MOC}-system is to produce projective characters
to approximate the columns of the decomposition matrix.

In order to simplify calculations based on the orthogonality
relations, we introduce some more notation allowing us to
write the orthogonality relations as matrix equations.
We first choose representatives $g_1 = 1, g_2, \ldots , g_s$
for the conjugacy classes of~$G$\@. We assume that the first~$s'$
of these elements are $p$-regular, and that $g_{s'+1}, \ldots , g_s$
are $p$-singular. Let~$C$ denote the $(s \times s)$-diagonal
matrix with entries $|C_G(g_1)|^{-1}, |C_G(g_2)|^{-1}, 
\ldots , |C_G(g_s)|^{-1}$\@. If $\Md = \{\lambda_1, \ldots ,
\lambda_m\}$ is a set of class functions of~$G$, we write
$$[\Md] = \left( \lambda_i(g_j)\right)_{1 \leq i \leq m, 1 \leq j \leq s}$$
for the $(m \times s)$-matrix giving the values of the functions in~$\Md$\@.
For example, $[\Irr(G)]$ is the table of ordinary irreducible
characters of~$G$\@. We always assume that any 
set of class functions is given in a particular order.

If $\Nd = \{\mu_1, \ldots , \mu_n\}$ is another set of
class functions, then $\langle \Md, \Nd \rangle$ is the matrix
$\left( \langle \lambda_i, \mu_j \rangle \right)_{1 \leq i \leq m, 
1 \leq j \leq n}$ of mutual inner products of the elements
of $\Md$ with the elements of~$\Nd$\@. Finally,
$\overline{\Md} := \{\bar{\lambda}_1, \ldots , \bar{\lambda}_m\}$,
where $\bar{\lambda}_i$ is the class function complex conjugate
to~$\lambda_i$\@.

With these notations we have:
\begin{equation}
\label{scalprod}
\langle \Md, \Nd \rangle = [\Md]\,C\,[\overline{\Nd}]^t.
\end{equation}
Writing~$E_\ell$ for the identity matrix of size~$\ell$, the
two orthogonality relations can be written as follows:
\begin{eqnarray}
\label{orthorelmat}
\left[\Irr(G)\right]\,C\,[\overline{\Irr(G)}]^t & = & E_s \\[1.5mm]
\left[\IPr(G)\right]\,C\,[\overline{\IBr(G)}]^t & = & E_{s'}.
\end{eqnarray}

\section{Blocks}
\markboth{MODULAR CHARACTER THEORY}{BLOCKS}

Let us recall the fundamental notions of Brauer's block theory.
Following Isaacs \cite[p.~272]{isa},
we first introduce the {\sl Brauer graph} of~$G$\@.
\index{Brauer graph}
Its vertex set is $\Irr(G)$\@. Two vertices $\chi$ and $\psi$ are
linked by an edge, if there exists $\varphi \in \IBr(G)$ such that
$d_{\chi\varphi} \neq 0 \neq d_{\psi\varphi}$\@.
The set of characters corresponding to a connected component of the 
Brauer graph is called a {\sl block} of~$G$\@.
\index{block}

If~$B$ is a union of blocks, 
we write $\Irr(B)$ and $\IBr(B)$ for the sets of irreducible
ordinary characters respectively Brauer characters in~$B$\@.
Let $B_1, \ldots , B_b$ denote the blocks of~$G$\@. Then
$$\Irr(G) = \bigcup_{i=1}^b \Irr(B_i),$$
a disjoint union. If ${\cal C}(G,B_i)$ denotes the space of
class functions spanned by $\Irr(B_i)$, we obtain
$${\cal C}(G) = \bigoplus_{i=1}^b {\cal C}(G,B_i).$$

For the set of irreducible Brauer characters we have 
$$\IBr(G) = \bigcup_{i=1}^b \IBr(B_i),$$
again a disjoint union.
By our definition of projective characters the ordinary 
constituents of a PIM are all
contained in a unique block. In this sense every PIM belongs to
a block. We write $\IPr(B_i)$ for the set of PIMs lying
in block~$B_i$\@.
Let ${\cal C}_{p'}(G,B_i)$ denote the $K$-space 
of class functions spanned by $\IBr(B_i)$\@. We thus have
$${\cal C}_{p'}(G) = \bigoplus_{i=1}^b {\cal C}_{p'}(G,B_i).$$
Therefore, the decomposition matrix splits into a direct
sum of block decomposition matrices~${\bf D}_B$,
and we can restrict our attention to the particular 
blocks, thus greatly
simplifying the problem.

The preceding remark of course only makes sense if one can determine
the blocks of~$G$ from the ordinary character table before
knowing the decomposition numbers. This is indeed the case.
For $\chi \in \Irr(G)$ let~$\omega_\chi$ denote the central
character corresponding to~$\chi$ defined by
\begin{equation}
\label{centralcharacter}
\omega_\chi(g) = \frac{|G|\chi(g)}{|C_G(g)|\chi(1)}, \quad g \in G.
\end{equation}
Then $\omega_\chi(g)$ is an algebraic integer and hence is in~$R$
for all $g \in G$ (see \cite[Theorem~3.7]{isa})\@.
Then (see \cite[Theorem 15.27]{isa}) two characters 
$\chi, \psi \in \Irr(G)$ are in the same
block, if and only if
\begin{equation}
\label{blockdistribution}
\frac{|G|\chi(g)}{|C_G(g)|\chi(1)} \equiv \frac{|G|\psi(g)}{|C_G(g)|\psi(1)}
\left(\mbox {mod } I \right)
\mbox{ for all } g \in G_{p'}.
\end{equation}
This criterion is independent of the choice of the maximal ideal~$I$
of~$R$ containing~$p$ (\cite[Theorem 15.18]{isa})\@.

The easiest blocks to deal with are the so-called blocks of
defect~$0$\@. They contain exactly one ordinary irreducible
character. This is at the same time the unique irreducible
Brauer character and the unique PIM in the block.
Thus the decomposition matrix of a block of defect~$0$
is the $(1 \times 1)$-matrix with single entry~$1$\@. 
An irreducible character~$\chi$ lies in a block of defect~$0$,
if and only if~$p$ does not divide $|G|/\chi(1)$
(\cite[Theorem~(15.29)]{isa})\@. In particular, if~$p$
does not divide the order of~$G$, the decomposition matrix
is the identity matrix.

It is often convenient to know the number of irreducible
Brauer characters in a block. Since the decomposition map is
surjective, the block decomposition matrices have maximal
$\Z$-rank and we have:
\begin{rem}
{\rm 
Let $\widehat{\Irr}(B) = \{ \hat{\chi} \mid \chi \in \Irr(B) \}$\@.
Then the $\Z$-rank of $[\widehat{\Irr}(B)]$ equals $|\IBr(B)|$\@.
}
\end{rem}

More generally, suppose we have a $K$-linear map
$$\alpha: {\cal C}(G) \longrightarrow {\cal C}(G),$$
which leaves~$\Irr(B)$ and $\IBr(B)$ invariant
and commutes with the decomposition map
$d: {\cal C}(G) \rightarrow {\cal C}(G)$\@.
Such homomorphisms arise, for example, from a group
automorphism of~$G$ or from the duality operation
on the modules of~$G$\@. We are then interested in
the number of characters in~$B$ fixed by~$\alpha$\@. 
We write $^\alpha\!\vartheta$ for the $\alpha$-image of a
character~$\vartheta$ and $\Irr(B)^\alpha$, respectively $\IBr(B)^\alpha$
for the set of $\alpha$-invariant elements.
The following result generalizes an observation of
Brauer and Feit (see \cite[(7.14)(b)]{gold}), which was communicated
to the authors by Kn\"orr.
\begin{prop}
Let~$\rho$ be the character of~$B$ defined by:
$$\rho(g) = 
\sum_{\chi \in {\rm Irr}(B)} \chi(g) {^\alpha\!\chi}(g^{-1}), \quad g \in G.$$
Then
$$|\Irr(B)^\alpha| = \langle {\bf 1}_G, \rho \rangle
\quad\quad \mbox{\rm and} \quad\quad
|\IBr(B)^\alpha| = \langle {\bf 1}_G, \hat{\rho} \rangle.$$
\end{prop} 
\begin{bew}
We have $\langle {\bf 1}_G, \rho \rangle =
\sum_{\chi \in {\rm Irr}(B)} \langle \chi, {^\alpha\!\chi} \rangle$,
from which the first equation follows.
To prove the second, let ${\bf D}_B = (d_{\chi\varphi})$
denote the decomposition matrix of~$B$\@. For $\varphi \in \IBr(G)$
let $\Phi_\varphi$ denote the PIM corresponding to~$\varphi$\@.
We then have:
$$\langle {\bf 1}_G, \hat{\rho} \rangle =
\sum_{\chi \in {\rm Irr}(B)} \langle {\chi}, {^\alpha\!\hat{\chi}} \rangle =$$
$$\sum_{\chi \in {\rm Irr}(B)} \sum_{\varphi \in {\rm IBr}(B)} 
d_{\chi\varphi} \langle \chi, {^\alpha}\!\varphi \rangle = 
\sum_{\varphi \in {\rm IBr}(B)} 
\langle \Phi_\varphi, {^\alpha}\!\varphi \rangle.$$
The result follows.
\end{bew}

\section{Generation of characters}
\markboth{MODULAR CHARACTER THEORY}{GENERATION OF CHARACTERS}
\label{generation}

We briefly describe the principal methods to produce
Brauer characters and projective characters.

\subsection{Restriction and induction}
\label{restandind}

Suppose~$H$ is a subgroup of $G$\@.
If $\vartheta \in {\cal C}(G)$ is a class function 
of~${G}$, the restricted class function $\vartheta_H$
is defined by:
\begin{equation}
\label{restrictiontosubgroup}
\vartheta_H(h) = \vartheta(h), \quad\quad h \in H.
\end{equation}
If $\vartheta$ is an ordinary or a Brauer character, 
then so is $\vartheta_H$\@. 
If $\vartheta \in {\cal C}({H})$ is a class function 
of~${H}$, the induced class function $\vartheta^G$ 
is defined as follows.
For $g \in G$ choose representatives $h_1, \ldots , h_s$
for the $H$-conjugacy classes contained in the $G$-conjugacy
class of~$g$\@.
Then
\begin{equation}
\label{inductionfromsubgroup}
\vartheta^G(g) = \sum_{i = 1}^{s}
\frac{|C_G(g)|}{|C_H(h_i)|}\,\vartheta(h_i).
\end{equation} 
The Frobenius reciprocity law asserts that the two operations
introduced above are adjoint operations with respect to the 
inner product $\langle \mbox{\ }, \mbox{\ } \rangle$\@.
This means that, given class functions $\vartheta$, $\rho$
of~$H$ respectively~$G$, we have
$$\langle \vartheta^G, \rho \rangle = \langle \vartheta, \rho_H \rangle.$$

Geck has observed that Frobenius reciprocity, together with 
Proposition~\ref{projconstruct}, can be used to
show that induction takes projective characters to projectives.
Let~$\Phi$ be a projective character of~$H$\@. Then
$\Phi^H$ is an ordinary character of~$G$; looking at the
definition~(\ref{inductionfromsubgroup}), we see that it
vanishes on the $p$-singular conjugacy classes.
If~$\varphi$ is an irreducible Brauer character of~$G$,
then $\langle \varphi, \Phi^G \rangle = \langle \varphi_H, \Phi \rangle \geq 0$,
since~$\varphi_H$ is a genuine Brauer character of~$H$, and
$\Phi$ is projective. Thus $\Phi$ is a genuine projective
character of~$G$\@.

Let $\vartheta$ denote a Brauer character of~$H$ afforded by
some re\-pre\-sentation ${\cal X}$ over~$F$\@.
Then $\vartheta^G$ is a Brauer character of~$G$; indeed, it is
the Brauer character afforded by the representation of~$G$
induced by~${\cal X}$ (see \cite[\S 10A, Exercise~10 in \S 18]{CRI})\@.
This observation in turn shows that the restriction of a projective
character to a subgroup is a projective character.

The following theorem tells us how to find enough projective
characters.
\begin{satz}
{\rm (Fong \cite[Lemma~1]{fongsplit}:)} Let ${\bf P}$ be the set of
projective characters of~$G$ obtained by
inducing the PIMs of all maximal subgroups. 
Then $\Ks(FG)$ is the $\Z$-span of ${\bf P}$\@. 
\end{satz}
\begin{bew}
This follows from Brauer's induction theorem. 
See, for example \cite[II 4.2,Corollary~2]{serre}\@.
\end{bew}
 
\subsection{Restriction to $p$-regular classes}

If $\vartheta \in {\cal C}({G})$ is a class function  
of~${G}$, let $\hat{\vartheta}$ be defined by~(\ref{pprimerestriction})\@.
Since the decomposition map is surjective, the set
$\widehat{\Irr}(G) = \{\hat{\chi} \mid \chi \in \Irr(G)\}$ spans $\Gs(FG)$
over the integers.

\subsection{Restriction to blocks}

Let~$B$ denote a union of $p$-blocks of~$G$\@.
Let $\psi = \sum_{\chi \in {\rm Irr}(G)} m_{\psi,\chi} \chi \in \Gs(KG)$\@.
Then the restriction of~$\psi$ to~$B$ is the generalized
character
\begin{equation}
\label{blockresordinary}
\psi_B = \sum_{\chi \in {\rm Irr}(B)} m_{\psi,\chi} \chi.
\end{equation}
Thus $\psi_B$ is just the projection of $\psi$ from ${\cal C}(G)$
onto ${\cal C}(G,B)$, the space of class functions spanned
by $\Irr(B)$\@. 

Similarly, if~$\vartheta$ is a generalized Brauer character,
$\vartheta_B$ is the projection of $\vartheta$ onto
${\cal C}_{p'}(G,B)$\@. Thus if
$\vartheta = \sum_{\varphi \in {\rm Irr}(G)} 
d_{\vartheta,\varphi} \varphi,$
then 
\begin{equation}
\label{blockresmodular}
\vartheta_B = \sum_{\varphi \in {\rm Irr}(B)} 
d_{\vartheta,\varphi} \varphi.
\end{equation}
We remark, that in each case we have to know a basis of 
the spaces ${\cal C}(G,B)$ respectively ${\cal C}_{p'}(G,B)$
in order to calculate the projections, but of course any
basis will do.

\subsection{Tensor products}

The product of any two characters is again a character. 
This is true for ordinary and Brauer characters.
For projective characters the following more general
statement is true. Namely, if $\vartheta$ is a Brauer
character and $\Psi$ is a projective character, then
the tensor product $\vartheta \Psi$ is a projective 
character. As Geck has observed, a proof of this fact can 
be given similar to those in \S~\ref{restandind}
by just using Proposition~\ref{projconstruct}. 

This statement for projective characters is extremely important for
the purposes of {\sf MOC}\@. It is used mainly in the case
when $\Psi$ is an ordinary irreducible character of
defect~$0$\@. These are of course the projective characters
which are readily available. The above statement means
in particular, that the product of two ordinary irreducible
characters is projective, if one of the factors is a
defect~$0$ character. The projectives obtained in this
way only depend on the character table of~$G$ and not on
the set of subgroups of~$G$\@. If we are able to determine
the irreducible Brauer characters of~$G$
by using only projectives of this sort, the decomposition
numbers are determined by the character table.

\subsection{Symmetrizations}

This paragraph is somewhat less elementary than the preceding
ones, since the facts we are going to describe
do not seem to follow by just considering characters.

Let~$\chi$ be an ordinary or modular character of~$G$\@.
The symmetric and skew square of~$\chi$ are the class
functions
$$\chi^{2+}(g) = \frac{1}{2}(\chi(g)^2 + \chi(g^2))$$ 
\mbox{\rm respectively}
$$\chi^{2-}(g) = \frac{1}{2}(\chi(g)^2 - \chi(g^2)), \quad\quad g \in G.$$
The symmetric or skew square of a character is again a
character, afforded by the symmetric respectively skew power
of the module affording~$\chi$ (if $p \neq 2$)\@. This observation
can be generalized to yield, for every positive integer~$r$
and every partition~$\lambda$ of~$r$ a {\sl symmetrized character}
\index{symmetrized character}
$\chi^\lambda$\@.

We now outline this method following James~\cite{jamsym}\@.
Let $V$ denote the $SG$-module affording the character~$\chi$\@.
Here, $S$ is one of the fields~$K$ or~$F$\@. Let $r$ denote a 
positive integer satisfying $r \leq n$, where $n = \mbox{\rm Dim}_S(V)$\@.
Let $V^{\otimes r}$ denote the $r$-fold tensor product of~$V$
on which~$G$ acts diagonally. Now~$V$ and hence $V^{\otimes r}$ 
is a module for the general linear group $GL(V)$ in its natural
action on~$V$ and thus
every $GL(V)$-section (a quotient module of a submodule)
of $V^{\otimes r}$ gives rise to an
$SG$-module by way of restriction. If the characteristic~$p$ of~$S$
is~$0$ or larger than~$r$, then $V^{\otimes r}$ is a semisimple
$GL(V)$-module. For every partition $\lambda$ of~$r$
there is an irreducible $GL(V)$-module $V_\lambda$, the 
{\sl Weyl module} corresponding to~$\lambda$, and $V^{\otimes r}$
\index{Weyl module}
is a direct sum of the Weyl modules.
If~$p$ is positive and less than or 
equal to~$r$, then $V^{\otimes r}$ still has a filtration
by Weyl modules~$V_\lambda$, but they are no longer direct summands.
Furthermore, in general the Weyl modules are not irreducible anymore.
Every Weyl module $V_\lambda$ has a unique top 
composition factor~$F_\lambda$\@. The other composition
factors of~$V_\lambda$ are of the form $F_\mu$ for partitions~$\mu$
which are strictly smaller than~$\lambda$ in the dominance 
order of partitions.
Thus, if the partitions of~$r$ are arranged in reverse lexicographic 
order, the matrix $M_{r,p}$ giving the multiplicities of the $F_\lambda$
in the $V_\mu$ is lower unitriangular.

Let~$\psi^\lambda$ denote the character of~$G$ afforded by
the Weyl module $V_\lambda$\@. It can be calculated with the help of
the character table of the symmetric group $S_r$\@. 
Let $\rho = (\rho_1,\rho_2, \ldots ,\rho_r)$, $\rho_1 \geq \rho_2 \geq 
\dots \geq \rho_r \geq 0$, be a partition of~$r$\@. 
Let $C(\rho)$ denote the centralizer in $S_r$ of an element
with cycle type $\rho$, and let $[\lambda](\rho)$ denote the value
of the irreducible character of~$S_r$ corresponding to~$\lambda$
on an element with cycle type~$\rho$\@. 
Then, for $p$-regular $g \in G$,
\begin{equation}
\label{weylsymmetrization}
\psi^\lambda(g) = \sum_{\rho \vdash r}\frac{1}{|C(\rho)|}[\lambda](\rho)
\prod_{i=1}^r \psi(g^{\rho_i}).
\end{equation}

If $p$ is non-zero and does not exceed~$r$, than we can improve
on the above symmetrization formula, since the Weyl modules
are not irreducible in general. 
Let $\Sigma{r,p} = M_{r,p}^{-1}\,[\Irr(S_r)]$\@. This matrix is called by
James an ``extended $p$-modular character table of~$S_r$''\@.
The character $\varphi^\lambda$ of~$G$ afforded 
by the simple module $F_\lambda$ can then be calculated by the
formula~(\ref{weylsymmetrization}) with $[\lambda]$ replaced
by the corresponding row of $\Sigma{r,p}$\@. 

In \cite{jamsym} James has calculated the matrices $M_{r,p}$
and $\Sigma_{r,p}$ for $r \leq 6$ and $p =2, 3$\@. For example,
labelling the rows and columns of the following matrices with
the partitions of~$3$ in reverse lexicographic order,
we have
$$
M_{3,3} = \left(\begin{array}{ccc} 1 & 0 & 0 \\ 1 & 1 & 0 \\ 0 & 1 & 1
\end{array}\right)
$$
and
$$
\Sigma_{3,3} = \left(\begin{array}{rrr} 1 & -1 & 1 \\ 1 & 1 & -2 \\ 0 & 0 & 3 
\end{array}\right).
$$

If $\psi$ is a $3$-modular character of a group~$G$ of degree
at least~$3$, the three $3$-modular symmetrizations are the following:
$$\varphi^{(1,1,1)}(g) = \frac{1}{6}\psi(g)^3 - \frac{1}{2}\psi(g^2)\psi(g)
+ \frac{1}{3}\psi(g^3)$$ \\[1.0mm]
$$\varphi^{(2,1)}(g) = \frac{1}{6}\psi(g)^3 + \frac{1}{2}\psi(g^2)\psi(g)
- \frac{2}{3}\psi(g^3)$$ \\[1.0mm]
$$\varphi^{(3)}(g) =  \psi(g^3).$$
The formulae of James have been implemented in the CAS system
\cite{CAS}\@.
In his thesis \cite{grabm}, Grabmeier has extended the calculations
of James to the case of $r \leq 9$ and $p = 2,3$\@.

\clearpage
\thispagestyle{empty}

\chapter{Computational Concepts}
\label{concepts}

We are now going to introduce the notions of basic sets and atoms
and derive some of their most important properties.
These notions are fundamental for the algorithms of {\sf MOC},
in particular those which try to prove the irreducibility of a 
Brauer character or the indecomposability of a projective character.
They allow us to translate some of the problems occuring 
during the calculation of decomposition matrices into the
language of Linear Algebra.
Furthermore, the {\sf MOC} methods generalize all elementary methods 
for calculating decomposition numbers described in the literature
(see \cite[Section~VI.3]{jake} for a collection af these methods)\@.
They can be systematically applied and are mechanized so as 
to be suitable for the treatment with computers.

We also formulate two fundamental computational problems which
we believe to be of independent interest.

It turned out that these concepts, originally designed for
application to sporadic groups, could also be applied successfully
to the study of modular characters of finite groups of Lie
type (see \cite{habil})\@. They allow to deal with all 
groups of a fixed rank and Dynkin type simultaneously. The reason for this
is that in a sufficient general situation we can find basic
sets which can be described independently of the underlying field
(see \cite{geckhiss})\@.

\section{Basic sets}

We continue with the notation introduced in the preceding section.
Here, we let~$B$ denote a fixed union of $p$-blocks of~$G$\@.
Let $\Irr(B) = \{\chi_1, \ldots ,\chi_n\}$\@.
The subgroups of $\Gs(KG)$, $\Gs(FG)$ and $\Ks(FG)$ generated
by the characters in~$B$ are denoted by $\Gs(B)$, $\Gs(\hat{B})$ and
$\Ks(B)$ respectively.
We write $\Gs(B)^+$ for the set of genuine characters
of $\Gs(B)$, and use an analogous notation for 
Brauer and projective characters. 

\begin{definit}
\label{basissys}
{\rm 
A {\sl basic set of Brauer characters} ${\bf BS}$ 
\index{basic set of Brauer characters}
({\sl basic set of projective characters} ${\bf PS}$)
\index{basic set of projective characters}
is a $\Z$-basis of $\Gs(\hat{B})$ (of $\Ks(B)$)
consisting of Brauer (projective) characters.
}
\end{definit}
We emphasize the fact that a basic set consists of proper
characters, in contrast to the usage in the literature.
One of the reasons is the experience that one should never give
away the information of a character being proper. We 
invest a great deal of time to reach the aim of having basic sets
at every stage of the calculations. 
In particular, we introduce an algorithm which for a given
set of proper characters determines a $\Z$-basis of their span
consisting again of proper characters. 

The first question we investigate, is the following.
Given a set of Brauer characters, how can we decide whether
it is a basic set?

\begin{lem}
\label{basicsettest}
Let $S \subset \Gs(\hat{B})^+$\@. Then~$S$ is a basic set
if and only if~$S$ is linearly independent over~$\Z$, and
every ordinary character of~$B$, restricted to the $p$-regular classes,
is a $\Z$-linear combination of elements of~$S$\@.
\end{lem}
\begin{bew}
The necessity of the condition is clear. The converse follows
from the surjectivity of the decomposition homomorphism.
\end{bew}

\begin{beisp}
\label{beisp1}
{\rm 
Here, and in the following examples, we 
denote irreducible characters by their degrees, which in our
examples identifies them uniquely. 
Consider the principal $7$-block~$B$ of Conway's simple
group~$Co_1$\@. 
It contains the following~$27$ ordinary irreducible characters:
 $${\footnotesize \left\{\begin{array}{rrrrr}
\underline{1}, &  \underline{276}, &  \underline{299}, &  
\underline{17250}, &  \underline{80730}, \\
\underline{94875}, &  \underline{822250}, &  \underline{871884}, &  
\underline{1821600}, &  \underline{2055625}, \\
\underline{9221850}, &  \underline{16347825}, &  \underline{21528000}, &  
\underline{21579129}, &  \underline{24667500}, \\
          {31574400}, &  \underline{57544344}, &  \underline{66602250}, &  
\underline{85250880}, &  \underline{150732800}, \\
\underline{163478250}, &            {191102976}, &  \underline{207491625}, &  
          {215547904}, &            {219648000}, \\
          {299710125}, &            {326956500} &  &  &
\end{array}\right\}.}$$
If we restrict these to the $7$-regular classes, 
the subset ${\bf BS}$ of the~$21$ underlined characters
is linearly independent, and the six
remaining characters not in ${\bf BS}$ 
can be written as $\Z$-linear combinations
of these according to the following table
(we omit the hat on top of the degrees):
{\small
$$
\begin{array}{r@{\hspace{ 6pt}}*{22}{r@{\hspace{ 2.9pt}}}}
 31574400 &  . &  . &  . &  . &  . &  -1 &  . &  . &  1 &  . &  . &  -1 &  1 &  . &  1 &  . &  . &  . &  . &  . &  . & \\
 191102976 &  1 &  . &  -1 &  . &  -1 &  . &  1 &  . &  -1 &  . &  -1 &  1 &  . &  1 &  . &  . &  . &  . &  . &  1 &  . & \\
 215547904 &  . &  . &  . &  -1 &  -1 &  1 &  1 &  1 &  -1 &  . &  -1 &  1 &  -1 &  . &  . &  . &  1 &  . &  . &  1 &  . & \\
 219648000 &  1 &  . &  . &  1 &  . &  . &  . &  . &  . &  . &  . &  1 &  . &  1 &  -1 &  . &  -1 &  -1 &  1 &  . &  1 & \\
 299710125 &  . &  . &  . &  1 &  . &  -1 &  -1 &  . &  1 &  -1 &  1 &  . &  . &  . &  . &  . &  -1 &  . &  1 &  . &  1 & \\
 326956500 &  . &  -1 &  1 &  -1 &  . &  1 &  . &  1 &  -1 &  1 &  -1 &  1 &  -1 &  1 &  -1 &  1 &  . &  -1 &  . &  1 &  1 & 
\end{array}
$$
}

Here, the columns correspond to the~$21$ characters in increasing order of
their degrees.
Thus the block contains exactly~$21$ irreducible Brauer characters
and our~$21$ characters are a basic set.
}
\end{beisp}
 
We have seen how to prove for a given set of
Brauer characters that it is a basic set.
Now the question arises, of how to do the same for projective
characters. By the theorem of Fong cited in \S~\ref{restandind} 
we know that in principle enough projective characters can be obtained
to span the space $\Ks({B})$\@.

We assume that we have constructed a set of projective characters
and selected from these a maximal linearly independent subset ${\bf PS}$\@.
How can we decide whether or not ${\bf PS}$ is a 
basic set? 

\begin{lem}
\label{projbastest}
Let ${\bf BS}$ and ${\bf PS}$ be two sets
of Brauer characters respectively projective
characters. We assume that each contains exactly $|\IBr(B)|$
elements. 
Let
\begin{equation}
\label{mutualscal}
U =  \langle {\bf BS}, {\bf PS} \rangle
\end{equation}
be the matrix of their mutual scalar products.
Then ${\bf BS}$ and ${\bf PS}$ are basic sets if and only if~$U$ 
is invertible over~$\Z$\@.
\end{lem}
\begin{bew}
We have
\begin{equation}
\label{bsinbr}
[{\bf BS}] = U_1\,[\IBr(B)] 
\end{equation}
and
\begin{equation}
\label{psinpr}
[{\bf PS}] = U_2^t\,[\IPr(B)]. 
\end{equation}
with integral matrices $U_1$ and~$U_2$\@.
From~(\ref{scalprod}) and the orthogonality relations we obtain:
$$U = [{\bf BS}]\,C\,[\overline{\bf PS}]^t = 
      U_1 [\IBr(B)]\,C\,[\overline{\IPr(B)}]^t U_2 =
      U_1 U_2.$$
Now ${\bf BS}$ respectively ${\bf PS}$ are bases if and only if~$U_1$
and~$U_2$ are invertible over~$\Z$ and the assertion follows.
\end{bew}

\begin{beisp}
\label{beisp2}
{\rm 
Consider again the principal $7$-block of~$Co_1$\@.
Let ${\bf BS}$ be as in Example~\ref{beisp1}\@.
After some calculations, which we do not want to comment on right here,
(but see Chapter~\ref{beweise}), we have found
a set ${\bf PS}$ of projective characters, such that
$U = \langle {\bf BS}, {\bf PS} \rangle$ 
is the matrix given below, where the columns correspond
to the projective characters.
Now $\Det\,U =  1$, 
so that ${\bf PS}$ is indeed a basic set of projectives by
Lemma~\ref{projbastest}\@.
$$
\begin{array}{r@{\hspace{10pt}}*{10}{r@{\hspace{ 4.8pt}}}*{11}{r@{\hspace{ 3.4pt}}}} \hline
  \Psi:  &  1 &  2 &  3 &  4 &  5 &  6 &  7 &  8 &  9 & 10 &
 11 & 12 & 13 & 14 & 15 & 16 & 17 & 18 & 19 & 20 & 21 \rule[- 7pt]{0pt}{ 20pt} \\ \hline
 1 &  1 &  . &  . &  . &  . &  . &  . &  . &  . &  . &  . &  . &  . &  . &  . &  . &  . &  . &  . &  . & . \rule[  0pt]{0pt}{ 13pt} \\
 276 &  1 &  1 &  . &  . &  . &  . &  . &  . &  . &  . &  . &  . &  . &  . &  . &  . &  . &  . &  . &  . & . \\
 299 &  1 &  . &  1 &  . &  . &  . &  . &  . &  . &  . &  . &  . &  . &  . &  . &  . &  . &  . &  . &  . & . \\
 17250 &  1 &  . &  1 &  1 &  . &  . &  . &  . &  . &  . &  . &  . &  . &  . &  . &  . &  . &  . &  . &  . & . \\
 80730 &  1 &  . &  . &  . &  1 &  . &  . &  . &  . &  . &  . &  . &  . &  . &  . &  . &  . &  . &  . &  . & . \\
 94875 &  1 &  1 &  . &  . &  . &  1 &  . &  . &  . &  . &  . &  . &  . &  . &  . &  . &  . &  . &  . &  . & . \\
 822250 &  1 &  . &  1 &  . &  . &  . &  1 &  . &  . &  . &  . &  . &  . &  . &  . &  . &  . &  . &  . &  . & . \\
 871884 &  1 &  . &  . &  . &  . &  . &  . &  1 &  . &  . &  . &  . &  . &  . &  . &  . &  . &  . &  . &  . & . \\
 1821600 &  1 &  1 &  . &  1 &  1 &  . &  . &  . &  1 &  . &  . &  . &  . &  . &  . &  . &  . &  . &  . &  . & . \\
 2055625 &  . &  . &  . &  1 &  . &  . &  . &  . &  . &  1 &  . &  . &  . &  . &  . &  . &  . &  . &  . &  . & . \\
 9221850 &  1 &  . &  . &  . &  . &  1 &  1 &  . &  . &  . &  1 &  . &  . &  . &  . &  . &  . &  . &  . &  . & . \\
 16347825 &  . &  . &  . &  1 &  2 &  . &  . &  . &  1 &  . &  . &  1 &  . &  . &  . &  . &  . &  . &  . &  . & . \\
 21528000 &  1 &  . &  . &  . &  . &  1 &  . &  1 &  . &  . &  . &  . &  1 &  . &  . &  . &  . &  . &  . &  . & . \\
 21579129 &  . &  . &  . &  . &  . &  . &  . &  . &  . &  . &  . &  . &  . &  1 &  1 &  . &  . &  . &  . &  . & . \\
 24667500 &  1 &  . &  . &  . &  1 &  . &  . &  . &  . &  . &  . &  1 &  . &  1 &  . &  . &  . &  . &  . &  . & . \\
 57544344 &  2 &  . &  1 &  1 &  . &  . &  1 &  . &  . &  . &  1 &  . &  . &  1 &  . &  1 &  . &  . &  . &  . & . \\
 66602250 &  . &  . &  . &  2 &  1 &  . &  . &  . &  1 &  1 &  . &  . &  . &  1 &  . &  . &  2 &  1 &  1 &  . & . \\
 85250880 &  2 &  . &  1 &  2 &  . &  . &  . &  . &  . &  1 &  . &  . &  . &  2 &  . &  1 &  1 &  . &  . &  . & . \\
 150732800 &  . &  . &  . &  1 &  . &  . &  . &  . &  . &  1 &  . &  . &  1 &  1 &  . &  1 &  1 &  1 &  . &  1 & 1 \\
 163478250 &  3 &  1 &  . &  1 &  . &  1 &  . &  . &  . &  . &  1 &  . &  1 &  1 &  . &  1 &  1 &  . &  . &  1 & 1 \\
 207491625 &  1 &  . &  . &  1 &  . &  . &  . &  . &  . &  1 &  . &  . &  1 &  3 &  . &  . &  3 &  1 &  1 &  1 & . 
\rule[- 7pt]{0pt}{  5pt} \\ \hline
\end{array}
$$
}
\end{beisp}

\pagebreak

\section{Atoms}

\begin{definit}
\label{atomsys}
{\rm 
A $\Z$-basis ${\bf BA}$ (${\bf PA}$) of~$\Gs(\hat{B})$ ($\Ks(B)$) is called 
\index{system of atoms}
a {\sl system of atoms of Brauer (projective) characters},
if every Brauer (projective) character is
a linear combination with non-neg\-ative coefficients
of elements of~${\bf BA}$ (${\bf PA}$)\@.
A {\sl Brauer atom (projective atom)} is a generalized 
\index{atom}
\index{Brauer atom}
\index{projective atom}
Brauer (projective) character,
which is a member of a system of atoms.
}
\end{definit}
In contrast to the definition of basic sets, 
we do not insist that our system of atoms
consists of proper characters.
The next result shows how to obtain a system of 
projective atoms from a basic set of Brauer characters.
Of course, Brauer characters and projective characters can be 
interchanged in the formulation of this and many of the 
following results. 

\begin{lem}
\label{PAausBS}
Let~${\bf BS}$ be a basic set of Brauer characters. If~${\bf PA}$
denotes the basis of $\Ks(B)$ dual to~${\bf BS}$ with 
respect to the bilinear form
$\langle \mbox{\ },\mbox{\ } \rangle$,
then~${\bf PA}$ is a system of projective atoms.
\end{lem}
\begin{bew}
We have $$[{\bf BS}] = U\,[\IBr(B)],$$ where the entries 
of~$U$ are non-negative integers, and $\Det\,U = \pm 1$\@.  
From $$[{\bf PA}]\,C\,[\overline{\bf BS}]^t = E_m$$ 
and the orthogonality relations we obtain
$$[\IPr(B)] = U^t\,[{\bf PA}],$$ 
and thus our assertion.
\end{bew}

\noindent Given ${\bf BS}$, it is not necessary for our purposes, 
though in principle 
possible, to know the elements of the dual basis~${\bf PA}$ 
explicitly as class functions.
Suppose we are given two basic sets ${\bf BS}$ and ${\bf PS}$\@.
The corresponding dual bases of atoms are denoted by
${\bf PA}$ and ${\bf BA}$\@.
Suppose furthermore, that ${\bf P}$ and ${\bf B}$ are
two sets of projective characters respectively Brauer 
characters. Then to express the elements of ${\bf P}$
in terms of the basis ${\bf PA}$ is equivalent to 
calculating the inner products of the elements of ${\bf P}$
with the characters in ${\bf BS}$\@.
Namely, let $V = \langle {\bf P}, {\bf BS} \rangle$\@. 
Then
$$[{\bf P}] = V\,[{\bf PA}].$$
Similarly, if $U = \langle {\bf PS}, {\bf B} \rangle$,
then 
\begin{equation} 
\label{BinBA} 
[{\bf B}] = U^t\,[{\bf BA}]. 
\end{equation}
In particular, if we take ${\bf BS}$ for the set
${\bf B}$ of Brauer characters in~(\ref{BinBA}), we obtain
the elements of ${\bf BA}$ as class functions by inverting
the matrix~$U^t$\@.

Before we continue with our exposition,
we give an easy but nevertheless very 
useful application of the notions introduced so far.
\begin{lem}
\label{indatom}
Suppose ${\bf PA}$ is a system of projective atoms.
If~$\Phi$ is a projective character contained 
in ${\bf PA}$, then~$\Phi$ is indecomposable.
\end{lem}
\begin{bew}
Let ${\bf PA} = \{\Phi_1, \ldots , \Phi_m\}$ with $\Phi = \Phi_1$\@.
Furthermore, let $\Psi_1, \ldots , \Psi_m$ denote the PIMs.
We write     
$$\Psi_i = \sum_{j=1}^m z_{ij}\Phi_j.$$
Then all $z_{ij} \geq 0$, and their matrix is invertible.
Since~$\Phi$ is projective, we can write
$$\Phi = \sum_{i=1}^m x_i \Psi_i$$
with $x_i \geq 0$\@. Therefore
$$\Phi_1 = \Phi = \sum_{j=1}^m \left( \sum_{i=1}^m x_i z_{ij}\right) \Phi_j.$$
From this it follows that $x_i = 0$ for all 
$i \neq i_0$ and $x_{i_0} = 1$, i.e.\ 
$\Phi = \Psi_{i_0}$\@.
\end{bew}

\begin{beisp}
\label{beisp3}
{\rm 
The projective character $\Psi_{15}$ of Example~\ref{beisp2}
is an atom, since it has inner product~$1$ with $21579129$, 
but~$0$ with all the other characters of ${\bf BS}$\@.
Hence it is a PIM\@. 
}
\end{beisp}

\noindent Of course there is an analogous statement for Brauer characters.
The above observation allows to prove indecomposability of projective
characters and irreducibility of Brauer characters. It is, however,
only a very special case of a much more general procedure,
which will be described in Chapter~\ref{improve}\@.

\section{The first fundamental problem}
\label{fundamentalproblemI}

We are now ready to formulate the first fundamental problem of our theory.
We assume that we are given two basic sets ${\bf BS}$ and ${\bf PS}$\@.
The corresponding dual bases of atoms are denoted by
${\bf PA}$ and ${\bf BA}$\@. We are interested in all 
matrices which can possibly be decomposition matrices with
respect to the given information.
First we calculate the matrix~$U$ of the scalar products of
the characters of ${\bf BS}$ with the characters of 
${\bf PS}$\@. Writing 
\begin{equation}
\label{bsinbrandpsinpr}
[{\bf BS}] = U_1\,[\IBr(B)] 
\quad\quad \mbox{\rm  and } \quad\quad
[{\bf PS}] = U_2^t\,[\IPr(B)],
\end{equation}
we have, as in the proof of Lemma~\ref{mutualscal}:
\begin{equation}
\label{faczer}
U = U_1\,U_2.
\end{equation}
Notice that by our definition of basic sets,~$U_1$ and~$U_2$
have non-negative entries. We write $X \geq 0$ to indicate
that every entry of the matrix~$X$ over the integers is non-negative.
Having ${\bf BS}$ at our hands, we know how to write the
restriction of the ordinary characters to the $p$-regular classes
in terms of the characters in ${\bf BS}$, say 
$$[\widehat{\Irr}(G)] = V\,[{\bf BS}].$$
Knowing~$U_1$ is then equivalent to knowing the decomposition 
matrix:
$${\bf D}_B = V\,U_1.$$
Thus the decomposition matrix is obtained from a solution of
\begin{equation}
\label{loeszer}
U = U_1\,U_2, \quad\quad U_1, U_2 \mbox{ unimodular }, \quad\quad U_1, U_2 \geq 0.
\end{equation}
Two solutions $(U_1,U_2)$ and $(U_1',U_2')$ of~(\ref{loeszer}) are
called equivalent, if there is a permutation matrix~$X$ such that
$U_1' = U_1X$ and $U_2' = X^tU_2$\@. Equivalent solutions lead to
the same set of Brauer characters arranged in a different order.
If~(\ref{loeszer}) has more than one equivalence class of solutions, 
then we cannot decide which of them determines the decomposition
matrix without further information.

In many cases, however, we have additional conditions which drastically
reduce the number of solutions of~(\ref{loeszer})\@.
Namely, if we are given a set of Brauer characters~${\bf B}$ and 
a set of projective characters~${\bf P}$, we can write:
$$[{\bf B}] = V\,[{\bf BS}] = (V\,U_1)\,[\IBr(B)]$$
and
$$[{\bf P}] = W\,[{\bf PS}] = (W\,U_2^t)\,[\IPr(B)]$$
with integer matrices~$V$ and~$W$\@. The solutions~$U_1$
and~$U_2$ of~(\ref{loeszer}) must therefore satisfy the following
conditions:
\begin{equation}
\label{nebenbed}
V\,U_1 \geq 0,\quad \mbox{ and } \quad W\,U_2^t \geq 0.
\end{equation}
So we have the first fundamental problem:

\medskip

\noindent {\bf Fundamental Problem I:} Given 
$V \in \Z^{s \times m}$ and $W \in \Z^{t \times m}$,
and a unimodular matrix $U \in \N^{m \times m}$,
find all solutions  $U_1, U_2 \in \Z^{m \times m}$ of 
$$U = U_1\,U_2, \quad \quad U_1, U_2 \geq 0.$$
satisfying
$$V\,U_1 \geq 0,\quad \mbox{ and } \quad W\,U_2^t \geq 0.$$

\medskip

\noindent Every solution of the fundamental problem leads to a
possible decomposition matrix by equation~(\ref{bsinbrandpsinpr}):
$[\IBr(B)] = U_1^{-1}\,[{\bf BS}]$\@. 
This shows that at the present state of knowledge the fact that
the determinant of the Cartan matrix ${\bf D}_B^t\,{\bf D}_B$
is a power of~$p$ does not
impose any restrictions on the solutions of the fundamental
problem (cf.\ the remarks in \cite[p.~$273$]{jake})\@.

We do not intend to solve this fundamental problem in its 
full generality. Usually the number of solutions is much too large
to be of any practical use. Instead we gradually improve on
the two basic sets using the conditions~(\ref{loeszer}) and~(\ref{nebenbed})
leading to a matrix~$U$ with smaller entries. How this
is done in practise with the help of methods of
integral linear programming will be described in Chapter~\ref{improve}\@.

\begin{beisp}
\label{beisp4}
{\rm 
We continue with our example of the principal $7$-block of~$Co_1$\@. 
Let ${\bf B}$ be the set of the six ordinary
characters which are not contained in ${\bf BS}$, restricted to
the $7$-regular classes of~$Co_1$\@. Then, of course,
the matrix~$V$ is just the one given in Example~\ref{beisp1}\@. 

We also have some projective characters which cannot be expressed
in terms of ${\bf PS}$ of Example~(\ref{beisp2})
entirely with non-negative coefficients. 
The corresponding matrix~$W$ is as follows.
$$
W = \left[ \begin{array}{*{21}{r@{\hspace{ 3.8pt}}}} 
  1 & -1 & -1 &  1 &  1 &  . &  . &  3 & -2 &  2 & -1
  &  3 &  1 &  1 &  . &  3 &  4 &  7 & -9 &  2 &-8 \\
  1 & -1 &  1 &  . &  . &  2 &  6 &  5 &  1 &  . &  2
  &  . &  5 &  . &  . &  8 & -1 &  . &  2 & -4 & 4 \\
  . &  . &  . &  . &  . &  . &  . &  . &  . &  . &  .
  &  1 &  . &  1 &  . & -1 &  . &  . &  1 &  3 & . \\
  . &  . &  . &  . &  . &  . &  . &  . &  . &  . &  .
  &  1 &  . &  1 &  . &  . &  1 &  1 & -1 & -1 & . \\
  1 & -1 &  . &  . &  . &  . &  . &  . &  . &  . & -1
  &  1 & -1 &  . &  . &  1 &  1 &  1 & -1 &  1 &-2 \\
  . &  . &  . &  . &  . &  . &  . &  . &  . &  . &  .
  &  . &  . &  1 &  . &  . & -1 &  1 &  2 &  . & 1 \\
  . &  . &  . &  1 &  2 &  . &  . &  . & -2 &  . &  .
  &  . &  . &  . &  . & -1 &  . &  . &  . &  . & . \\
  . &  . &  . &  . &  . &  . &  . &  . &  . &  . &  .
  &  . &  . &  1 &  . &  . &  1 & -1 & -2 &  1 &-1 \\
  . &  . &  . &  . &  . &  . &  . &  . &  . &  . &  .
  &  . &  . &  . &  1 &  . &  . &  . &  . &  2 &-2 \\
  . &  . &  . &  . &  . &  . &  . &  . &  . &  . &  .
  &  1 &  . &  . &  . &  . &  1 &  . & -1 &  . &-1 \\
  . &  . &  . &  . &  . &  . &  . &  . &  . &  . &  1
  &  . &  1 &  . &  . &  2 &  1 & -1 & -1 & -1 & . \\
  . &  . &  . &  . &  . &  . &  . &  . &  . &  . &  2
  &  . &  2 &  . &  . &  1 &  . & -1 &  1 & -2 & 1 \\
  . &  . &  . &  . &  . &  . &  . &  . &  1 &  . &  .
  &  . &  . &  . &  . &  . &  . & -1 &  1 &  1 & . \\
  . &  . &  . &  . &  . &  . &  . &  . &  . &  . &  1
  &  . &  1 &  . &  1 &  1 &  1 &  . & -2 &  . &-1 \\
  . &  . &  . &  . &  . &  . &  . &  . &  . &  . &  2
  &  . &  2 &  . &  . &  . &  1 & -2 &  . & -2 & 1 \\
  . &  . &  . &  . &  . &  . &  1 &  . &  . &  . &  .
  &  . &  1 &  . &  . &  1 &  . &  . &  . & -1 & .
\end{array} \right]
$$
With these data, there are only~$4$ solutions of the fundamental
problem~I, namely, up to a permutation of the columns,~$U_1$ 
is one of the matrices given in Table~\ref{possol}\@.
\begin{table}[tb]
\caption{\label{possol} Possible decomposition matrices for the
principal block of $\mbox{\em Co}_1$}
$$
\begin{array}{r@{\hspace{10pt}}*{21}{r@{\hspace{ 7pt}}}} \hline
 1 &  1 &  . &  . &  . &  . &  . &  . &  . &  . &  . &  . &  . &  . &  . &  . &  . &  . &  . &  . &  . & .  \rule[  0pt]{0pt}{ 13pt} \\
 276 &  . &  1 &  . &  . &  . &  . &  . &  . &  . &  . &  . &  . &  . &  . &  . &  . &  . &  . &  . &  . & . \\
 299 &  . &  . &  1 &  . &  . &  . &  . &  . &  . &  . &  . &  . &  . &  . &  . &  . &  . &  . &  . &  . & . \\
 17250 &  . &  . &  1 &  1 &  . &  . &  . &  . &  . &  . &  . &  . &  . &  . &  . &  . &  . &  . &  . &  . & . \\
 80730 &  1 &  . &  . &  . &  1 &  . &  . &  . &  . &  . &  . &  . &  . &  . &  . &  . &  . &  . &  . &  . & . \\
 94875 &  . &  1 &  . &  . &  . &  1 &  . &  . &  . &  . &  . &  . &  . &  . &  . &  . &  . &  . &  . &  . & . \\
 822250 &  . &  . &  1 &  . &  . &  . &  1 &  . &  . &  . &  . &  . &  . &  . &  . &  . &  . &  . &  . &  . & . \\
 871884 &  1 &  . &  . &  . &  . &  . &  . &  1 &  . &  . &  . &  . &  . &  . &  . &  . &  . &  . &  . &  . & . \\
 1821600 &  . &  1 &  . &  . &  . &  . &  . &  . &  1 &  . &  . &  . &  . &  . &  . &  . &  . &  . &  . &  . & . \\
 2055625 &  . &  . &  . &  1 &  . &  . &  . &  . &  . &  1 &  . &  . &  . &  . &  . &  . &  . &  . &  . &  . & . \\
 9221850 &  . &  . &  . &  . &  . &  1 &  1 &  . &  . &  . &  1 &  . &  . &  . &  . &  . &  . &  . &  . &  . & . \\
 16347825 &  . &  . &  . &  . &  1 &  . &  . &  . &  1 &  . &  . &  1 &  . &  . &  . &  . &  . &  . &  . &  . & . \\
 21528000 &  . &  . &  . &  . &  . &  1 &  . &  1 &  . &  . &  . &  . &  1 &  . &  . &  . &  . &  . &  . &  . & . \\
 21579129 &  . &  . &  . &  . &  . &  . &  . &  . &  . &  . &  . &  . &  . &  1 &  1 &  . &  . &  . &  . &  . & . \\
 24667500 &  1 &  . &  . &  . &  1 &  . &  . &  . &  . &  . &  . &  1 &  . &  1 &  . &  . &  . &  . &  . &  . & . \\
 57544344 &  . &  . &  1 &  . &  . &  . &  1 &  . &  . &  . &  1 &  . &  . &  . &  . &  1 &  . &  . &  . &  . & . \\
 66602250 &  . &  . &  . &  1 &  . &  . &  . &  . &  1 &  a &  . &  . &  . &  . &  . &  . &  1 &  . &  . &  . & . \\
 85250880 &  . &  . &  1 &  1 &  . &  . &  . &  . &  . &  1 &  . &  . &  . &  . &  . &  1 &  . &  1 &  . &  . & . \\
 150732800 &  . &  . &  . &  . &  . &  . &  . &  . &  . &  1 &  . &  . &  . &  . &  . &  1 &  . &  . &  1 &  . & . \\
 163478250 &  . &  1 &  . &  . &  . &  1 &  . &  . &  . &  . &  1 &  . &  1 &  b &  . &  . &  . &  . &  . &  1 & . \\
 207491625 &  . &  . &  . &  1 &  . &  . &  . &  . &  . &  a &  . &  . &  . &  . &  . &  . &  1 &  1 &  . &  . & 1 
\rule[- 7pt]{0pt}{  5pt} \\ \hline
\end{array}
$$
\end{table}
There, $a, b \in \{0, 1\}$\@. 
We do not intend to give a proof of this 
fact here. An indication of how this is done is given in Chapter~\ref{improve}.
}
\end{beisp}

If we are left with a small number of non-equivalent
solutions of the fundamental
problem, we can distinguish cases. We take each solution in turn
and pretend that it leads to the true decomposition matrix.
Then we tensor all irreducible Brauer characters with themselves,
and write the resulting characters in terms of the irreducibles.
If this leads to expressions with negative coefficients, our assumption
was wrong and we consider the next solution.

In principle one could extend the fundamental problem by
adding the conditions arising from tensor products.
Suppose ${\bf BS} = \{\vartheta_1, \ldots , \vartheta_m\}$
and ${\bf PS} = \{\Psi_1, \ldots , \Psi_m\}$\@.
Let $a_{ijl} = \langle \vartheta_i\vartheta_j, \Psi_l \rangle$\@.
For $i = 1, 2$, let $U_i^{-1} = (\underline{u}_{rs}^i)_{r,s}$\@.
Then the solution $(U_1,U_2)$ which determines the decomposition
matrix must satisfy:
$$\sum_{r=1}^m \sum_{s=1}^m \sum_{t=1}^m 
a_{rst} \underline{u}_{ir}^1 \underline{u}_{js}^1 \underline{u}_{tl}^2 \geq 0,
\quad \mbox{ for all } i,j,l,$$
since the expression on the left hand side equals
$\langle \varphi_i\varphi_j, \Phi_l \rangle$, where
$[\varphi_i]$ respectively $[\Phi_i]$ is the $i$-th row
of $U_1^{-1}\,[{\bf BS}]$ respectively $U_2^{-t}[{\bf PS}]$\@.

In the ideal case, the matrix~$U$ in~(\ref{mutualscal})
is the identity matrix. Then we have solved our problem.

\begin{prop}
\label{idealfall}
Let ${\bf BS}$ and ${\bf PS}$ be two basic sets such that the
matrix~$U$ of their mutual scalar products
is the identity matrix. Then ${\bf BS}$ consists of
the irreducible Brauer characters and ${\bf PS}$ of the PIMs\@.
\end{prop}
\begin{bew}
This immediately follows from~(\ref{faczer}),~(\ref{bsinbr}) 
and~(\ref{psinpr})\@.
\end{bew}

\noindent As already mentioned above, 
{\sf MOC} tries to reach at this 
final stage by improving the two basic sets 
${\bf BS}$ and ${\bf PS}$ step by step by making use of the 
conditions~(\ref{nebenbed})\@.

We finally indicate 
a proof of a fact mentioned in \cite[Remark~6.3.34]{jake},
namely that the information
gained by inducing Brauer characters from a subgroup is
the same as that obtained by restricting projective
characters to this same subgroup.
It is not really necessary to use the language of basic sets
to do so, but it can be used to quantify the information.

Let~$H$ be a subgroup of~$G$\@. We assume that the decomposition
matrix ${\bf D}_H$ of~$H$ is known. Let~$F$ denote 
the $G$-$H$-induction-restriction matrix, i.e.,
$$[\Irr(G)_H] = F\,[\Irr(H)],$$
where $\Irr(G)_H$ denotes the set of restrictions to~$H$ of
the irreducible characters of~$G$\@. 
Then, if $\Irr(H)^G$ is the set of 
characters of~$G$ induced from the irreducible characters of~$H$,
we have 
$$[\Irr(H)^G] = F^t\,[\Irr(G)],$$
by Frobenius reciprocity.
Let ${\bf BS}$ be a basic set for the Brauer characters of~$G$\@.
Since the decomposition map is surjective, there is an 
integral matrix~$Y$ such that
$$[{\bf BS}] = Y\,[\widehat{\Irr}(G)].$$
We choose such a~$Y$ and fix it in the following.
Then 
$$[\IBr(G)] = X\,[\widehat{\Irr}(G)],$$
with $X = U_1^{-1}\,Y$\@. By the surjectivity of
the decomposition homomorphism and the definition of the PIMs, 
this implies that
$$[\Irr(G)] = X^t\,[\IPr(G)].$$
Now $[\IPr(H)] = {\bf D}_H^t\,[\Irr(H)]$, and so
\begin{eqnarray*}
[\IPr(H)^G] & = & {\bf D}_H^t\,[\Irr(H)^G] \\
            & = & {\bf D}_H^t\,F^t\,X^t\,[\IPr(G)] \\
            & = & {\bf D}_H^t\,F^t\,X^t\,U_2^{-t}\,[{\bf PS}] \\
            & = & {\bf D}_H^t\,F^t\,Y^t\,U^{-t}\,[{\bf PS}].
\end{eqnarray*}
Thus the information gained by inducing the projective characters
from $H$ is the following condition for~$U_2$ coming from~(\ref{nebenbed})
$$({\bf D}_H^t\,F^t\,Y^t\,U^{-t})\,U_2^t \geq 0.$$
On the other hand, restricting the irreducible Brauer characters
of~$G$ down to~$H$ we obtain
\begin{eqnarray*}
[\IBr(G)_H] & = & X\,[\widehat{\Irr}(G)_H] \\
            & = & X\,F\,[\Irr(H)] \\
            & = & U_1^{-1}\,Y\,F\,{\bf D}_H\,[\IBr(H)] \\
            & = & U_2\,U^{-1}\,Y\,F\,{\bf D}_H\,[\IBr(H)],
\end{eqnarray*}
giving exactly the same conditions as above.

\section{The second fundamental problem}
\label{fundamentalproblemII}

In general, the two matrices of~(\ref{nebenbed}) have many rows.
But if one of these rows is a sum of others, then the
row obviously contains no condition at all and can be deleted.
As our experience shows, a quick 
procedure to throw away these obsolete rows is of enormous help.
For example, in the course of the calculation of the $7$-modular
characters of the Conway group $2.Co_1$ we had to produce more
than $10000$ projective characters in order to find the few
which are of any use.
This leads us to the formulation of our second fundamental 
computational problem. Before we do this, we give an
ad hoc definition, which is only needed to formulate
the problem.

\begin{definit}{\rm Let~$W$ be a finite subset of $\Z^{1 \times s}$,
i.e.,\ a set of rows of integers of length~$s$\@. 
A subset $W_0 \subseteq W$ is called {\sl essential}, if for all $w \in W$
\index{essential subset}
there are non-negative integers $n_v$, $v \in W_0$,
such that $w = \sum_{v \in W_0} n_v v$.
}
\end{definit}
If the zero vector is not in the positive $\Z$-span of~$W$, then 
there is a unique minimal essential subset of~$W$\@.

\medskip

\noindent {\bf Fundamental Problem~II:} Let~$W$ be a finite 
subset of $\Z^{1 \times s}$\@. 
Determine an essential set $W_0 \subseteq W$, with
$|W_0|$ as small as possible.

In general the above problem is slightly more complicated since
we do not have a set of vectors but a sequence with some vectors
repeated. 

\section{Relations}
\label{relations}

If we are given a basic set of Brauer characters ${\bf BS}$ and we know
how to write the basic set characters in terms of the irreducible
Brauer characters, we can derive the decomposition matrix as
remarked at the beginning of Section~\ref{fundamentalproblemI}\@.
The following notation is convenient in our situation.

\begin{definit}
{\rm 
A {\sl relation} is a linear 
\index{relation}
dependence over~$\Z$ between generalized characters.
More specificly: If~${\bf B}$ is any basis of $\Gs(\hat{B})$,
and if~$\vartheta \in \Gs(\hat{B})$,
then the expression of~$\vartheta$ in terms of~${\bf B}$
is called a {\sl relation with respect to}~${\bf B}$\@.
\index{relation with respect to a basis}
A similar notion is used for projective characters.
}
\end{definit}

\begin{beisp}
\label{beisp5}
{\rm
The matrix in Example~\ref{beisp1}
gives the relations of the six
ordinary characters not contained in ${\bf BS}$ in terms of this
basic set.
The matrix~$W$ in Example~\ref{beisp4} gives the relations of some projectives
in terms of the basic set of projectives ${\bf PS}$ of 
Example~\ref{beisp1}\@.
}
\end{beisp}

If we are given a set~${\bf BS}$ of ordinary irreducible characters,
restricted to the $p$-regular classes,
then in order to show that they form a basic
set it suffices, as we have seen above, to express every other
irreducible ordinary character, restricted to the $p$-regular classes,
as a $\Z$-linear combination of the elements of~${\bf BS}$\@.
If this can be done, we have at the same time found the relations
which are needed to reconstruct the whole decomposition matrix.

Sometimes the relations with respect to a basis are stored
rather than the characters themselves. This is for reasons of space
but also to simplify and speed up the calculation of inner
products. Let ${\bf BS}$ be a basic set of Brauer characters
with dual basis of projective atoms ${\bf PA}$\@.
Suppose furthermore, that ${\bf P}$ and ${\bf B}$ are
two sets of projective characters respectively Brauer
characters. Then in order to find the matrix 
$\langle {\bf P}, {\bf B} \rangle$ of mutual inner products,
it suffices to write 
\begin{equation}
\label{BinBS}
[{\bf B}] = V\,[{\bf BS}]
\end{equation}
and
\begin{equation}
\label{PinPA}
[{\bf P}] = W\,[{\bf PA}],
\end{equation}
in order to find 
\begin{equation}
\label{lemma19}
\langle {\bf P}, {\bf B} \rangle = W\,V^t.
\end{equation}

\begin{beisp}
\label{beisp6}
{\rm 
The scalar product of a Brauer character in the set  
$$\{31574400, 191102976, 215547904, 219648000, 299710125, 326956500\}$$
of Example \ref{beisp1} 
with a projective of ${\bf PS}$ of Example~\ref{beisp2}
is given by matrix multiplication 
of the rows of matrix in Example~\ref{beisp1}
with the columns of 
the matrix of Example~\ref{beisp2}\@.
Thus, for example, the Brauer character $191102976$ has inner
product $1 - 1 - 1 + 1 - 1 - 1 + 3 = 1$ with the projective
$\Psi_1$ and inner product $1 - 1 = 0$ with projective~$\Psi_2$\@.
}
\end{beisp}

\section{Special basic sets}
\label{special}

We keep the notation of the previous section. 
The situation is considerably simpler if we have basic sets of
Brauer characters of a special nature.

\begin{definit}
{\rm 
A basic set ${\bf BS}$ of Brauer characters is called {\sl special},
\index{special basic set}
if it is of the form ${\bf BS} = \{\varphi_1, \ldots ,\varphi_m\}$, where the~$\varphi_i$
are restrictions of ordinary characters of~$G$ to the
$p$-regular conjugacy classes, i.e., $\varphi_i = \hat{\chi}_{j_i}$
for some $\chi_{j_i} \in \Irr(B)$\@.
}
\end{definit}

Let ${\bf BS}_0$ denote a special basic set and let ${\bf PA}_0$
be the corresponding set of projective atoms. Then write
\begin{equation}
\label{specialrel}
[\widehat{\Irr}(B)] = S\,[{\bf BS}_0].
\end{equation}
The matrix~$S$ thus gives the relations arising from the
restrictions of the ordinary characters to the $p$-regular
conjugacy classes with respect to the special basic set
${\bf BS}_0$\@. In the present version of {\sf MOC} we 
always fix a special basic set and store the matrix~$S$\@.
In particular, {\sf MOC} will only work if there exists
such a special basic set. 

It is not clear at all whether there is always a special basic set.
In $p$-soluble groups there is always one by the theorem of Fong and
Swan. If the ordinary irreducible characters in a $p$-soluble
group are sorted in increasing order of their degrees, 
then the first maximal linear independent subsequence of characters 
consists of the irreducible Brauer characters. 
In blocks with cyclic defect group there is a special basic set
by the results of Brauer and Dade. The authors have checked all blocks for all
sporadic simple groups and all possible~$p$, and found a special
basic set in all cases. In~\cite{geckhiss} it is shown that there
is a special basic set for groups of Lie type in non-defining characteristics
under some additional hypotheses\@. In searching for a counterexample
we would first try a group of Lie type in defining characteristic.

A special basic set has various advantages. First of all,
each basic set character contains only modular constituents
of a single block, and so the subset ${\bf BS}_0 \cap B$ is a 
special basic set of~$B$\@. This is used to find the contribution
of a generalized Brauer character to~$B$ by~(\ref{blockresmodular})\@.

The second advantage of a special basic set is the fact that
it allows to express a projective character in terms of
${\bf PA}_0$ without much work.
Namely, let ${\bf P}$ denote a set of projective characters
expressed in terms of the ordinary irreducible characters.
That is to say we know the matrix
\begin{equation}
\label{PinIrr}
Y = \langle {\bf P}, \Irr(B) \rangle.
\end{equation}
We then obtain the matrix 
$$Y_0 = \langle {\bf P}, {\bf BS}_0 \rangle$$
by just deleting the columns of~$Y$ corresponding 
to characters in $\Irr(B) \setminus {\bf BS}_0$\@.
Thus 
\begin{equation}
\label{PinPA0}
[{\bf P}] = Y_0\,[{\bf PA}_0].
\end{equation}
Provided we have a special  basic set
of Brauer characters, this observation allows to restrict 
(in an obvious sense) the decomposition matrix to it,
and then give the columns a new interpretation: These are
just the PIMs expressed in the dual basis of atoms.

If one tries to calculate the decomposition matrix by constructing
projective characters and expressing them in terms of ordinary
irreducibles, it suffices to calculate their scalar products
with the elements belonging to a special basic set.
This is often an enormous saving of space and time.

\begin{beisp}
\label{beisp7}
{\rm
The basic set ${\bf BS}$ of Example~\ref{beisp1} is special,
and hence the columns of the matrix of Example~\ref{beisp2}
give the coefficients in the expression of the corresponding
projectives in terms of the basis ${\bf PA}$ dual to
${\bf BS}$\@.
}
\end{beisp}

By comparing 
(\ref{specialrel}), (\ref{PinIrr}) and~(\ref{PinPA0})
with (\ref{BinBS}), (\ref{PinPA}) and~(\ref{lemma19}),
we have
$$Y = Y_0\,S^t.$$
This can be interpreted as finding the matrix~$Y$
if a set of projectives is given in terms of ${\bf PA}_0$\@.
This will be important later on, since in the process
of improving the projectives ${\sf MOC}$ uses their
expressions via the matrix~$Y_0$\@. 

We finally specialize to the case that ${\bf P} = {\bf PS}$
is a basic set of projectives.
Then, writing 
$$X = \langle {\bf PS}, \Irr(B) \rangle, 
\quad\quad \mbox{\rm and} \quad\quad 
X_0 = \langle {\bf PS}, {\bf BS}_0 \rangle,$$
we have 
$$X = X_0\,S^t.$$
By the orthogonality relations we also have
\begin{equation}
\label{BS0inBA}
[{\bf BS}_0] = X_0^t\,[{\bf BA}].
\end{equation}
If ${\bf B}$ is a set of Brauer characters, $[{\bf B}] = V\,[{\bf BS}_0]$,
we obtain
$$[{\bf B}] = VX_0^t\,[{\bf BA}]$$
which can also be expressed as
$$\langle {\bf B}, {\bf PS} \rangle = VX_0^t.$$
This means that if ${\bf B}$ is given
in terms of ${\bf BA}$, i.e., by its inner products
with ${\bf PS}$, say ${\bf B} = W\,[{\bf BA}]$,
we can express ${\bf B}$ in terms of
${\bf BS}_0$ as follows:
$$[{\bf B}] = WX_0^{-t}\,[{\bf BS}_0].$$

\section{Triangular decomposition matrices}

The solution of the Fundamental Problem~I is drastically simplified
if the matrix~$U$ is a lower triangular matrix.
Since~$U$ is unimodular, the entries on the main diagonal
are all equal to~$1$\@.
The two matrices~$U_1$ and~$U_2$ are then
unimodular triangular matrices, too. One obtains~$U_1$, say,
by iteratively subtracting columns of~$U$ from the right
(compare also \cite[6.3.22]{jake} and the subsequent remarks)\@.
Closely related to this special case is the notion of wedge shape
of the decomposition matrix.

\begin{definit}
{\rm
The decomposition matrix of~$B$ has {\sl wedge shape}, if the ordinary
\index{wedge shape of decomposition matrix}
irreducible characters can be ordered such that the decomposition
matrix has the following form:
$$\left(\begin{array}{cccccc}
1& & & & &   \\
 &1& & & &   \\
 & &.& & &   \\
 &*& &.& &   \\
 & & & &.&   \\
 & & & & &1  \\
 & & & & &   \\
 \multicolumn{6}{c}{\Large *} \\
 & & & & &  
\end{array} \right),$$
where the entries above the main diagonal are all~$0$\@.
}
\end{definit}

\noindent In a symmetric group every decomposition matrix has wedge shape
(see e.g.~\cite[Theorem 6.3.60]{jake})\@. The same is true for the
general linear and unitary groups in non-defining characteristics, 
as Dipper~\cite{dippI} and Geck~\cite{geckuni} have shown.

We obviously have the following connection.
\begin{lem}
The decomposition matrix has wedge shape, if and only if
there is a special basic set of Brauer characters
${\bf BS}_0$ and a basic set ${\bf PS}$ of projective characters
such that the matrix~$U$ of~{\rm (\ref{mutualscal})} is a lower
triangular matrix. \blackbox
\end{lem} 


\clearpage
\thispagestyle{empty}

\chapter{Data structures}

In this chapter we describe the principal data 
structures underlying the {\sf MOC}-system. We also present a
new and remarkable result due to H.~W.~Lenstra on certain
integral bases of abelian number fields. This result had been
suggested by experimental evidence gained in the 
preparation of the {\sf MOC}-system.

Based on the theory described in Chapter \ref{concepts},
{\sf MOC} consists of a collection of FORTRAN-programs and of unformatted
FORTRAN-files, which contain characters and  information about the characters.
The programs call
the  subroutine RRR to perform reading and writing of files and
to do the basic operations on long integers.
The files are used for exchanging information between 
the different programs and for recording the calculations, which 
have been done. The representation of the various data is the same for
the external and the internal format. 
The programs solve very specific tasks, like
tensoring two lists of Brauer characters or multiplying two matrices.
{\sf MOC} lacks a memory management and a language
on top of it. Instead more complicated tasks are achieved by concatenating 
suitable programs via UNIX Bourne-shell scripts.
We give an example for such a shell script in the next chapter.
We begin by describing the basic data types, which the programs use,
and the operations, which are performed on them.

\section{Long integers}

The fundamental data type in {\sf MOC} is the long integer. A long integer $n$
is represented as a list of coefficients with respect to the basis $10^4$\@.
If $n \ge 0$ then
\[
n = \sum_{i=0}^{k} a_{i} 10^{4i}, \quad 0 \le a_i \le 9999, \quad a_k \neq 0.
\]
The list of coefficients $a_{k-i+1}$ is stored in a 1-dimensional
FORTRAN-array $r[\mbox{\ }]$ and $10000$ is added to $a_0$\@. That is to say,
$r[i] = a_{k - i + 1}$, $1 \leq i \leq k$, and $r[k+1] = a_0 + 10000$\@.
For example, if
$n = 123456789$ then 
\[
r[1]=1, \quad r[2]=2345, \quad r[3]=16789.
\]
In case $n$ is negative, we store the coefficients $-a_{k-i+1}$ and
add $20000$ to $-a_0$\@. If $n=-123456789$, then
\[
r[1]=1, \quad r[2]=2345, \quad r[3]=26789.
\]
In order to do the basic operations on long integers, the
FORTRAN programs call the arithmetic subroutine RRR, which performs
addition, subtraction, multiplication
and division of long integers and returns the result to the calling program.
The algorithms used by RRR can be found for example in Knuth's book
\cite[Chapter~4.3]{knuthII}\@.

\section{Files}

The data stored in memory or on file consists of numbers  coded
as described above and of  
separators which are numbers in the range between 30000 and 31000.
The separators are used as labels, which indicate the sort of data,
following the separator.  The separator 30900 for example is
followed by class functions evaluated on certain conjugacy classes.

{\sf MOC} keeps track of the calculation performed for a specific group
$G$ and a specific prime $p$ on several files. We shortly mention the
four most important ones.

\begin{itemize}
\item[(a)]{The file {\tt G.p} contains the following information:
A special basic set ${\bf BS}_0$
 stored under the label 30900, the matrix of relations $S$
stored under 30550 expressing the ordinary characters in terms
of ${\bf BS}_0$, a basic set of projectives ${\bf PS}$
in terms of the projective atoms dual to the special basic set,
an actual basic set  of Brauer characters ${\bf BS}$
in terms of the Brauer atoms dual to  ${\bf PS}$. Furthermore
it contains information related to the distribution of the characters
into blocks. Changing from a representation of characters in a certain
basis to a representation in a different basis or writing the
characters as class functions is  achieved by using the formulae
given in Section \ref{special}.}
\item[(b)]{The file {\tt G.p.bras} contains a list of Brauer characters
stored under the label 30500 expressed in terms of ${\bf BS}_0$.
Whenever new Brauer characters are generated they are put onto this file.}
\item[(c)]{The file {\tt G.p.proj} contains a list of projective characters
stored under the label 30700 expressed in terms of the 
ordinary irreducible characters.
Whenever new projective characters are generated they are put onto this file.}
\item[(d)]{The file {\tt G.p.info} contains  information concerning the
calculations already performed by {\sf MOC}. This can be used 
to give a conventional proof of the results obtained by {\sf MOC}.
We describe in a later section this important feature 
in more detail.}
\end{itemize}

\section{Integers in abelian number fields}

The second basic data type deals with character values,
which are algebraic integers from abelian number fields.
Our representation of character values differs from those used in other
systems, such as for example {\sf CAS} \cite{CAS}\@.
This is motivated by the fact that Brauer character tables tend to contain
more irrationalities than ordinary character tables. One of the important
consequences of the new {\sf MOC}-format for character values
is that the basic algorithms involving
Brauer characters can be reformulated as algorithms for ${\Z}$-lattices,
avoiding
explicit calculations in cyclotomic fields. The {\sf MOC}-format also seems
to give  a more compact representation of a character table.

Let $K$ be an algebraic number field  of degree $d$ over $\Q$
with ring of integers $R$.  For performing
arithmetic in $R$ we choose an integral basis $b_1, \ldots, b_d$.
The elements of $R$ are represented by the coefficients in this basis.
The addition of two elements is now easily performed by adding their coefficient
vectors. 
Next we produce the multiplication matrix of the basis elements,
i.e., we store the coefficients of all products $b_i \cdot b_j$.
The multiplication of two arbitrary elements of $R$ is now achieved by
matrix multiplication of the two coefficient vectors with
the multiplication matrices of the basis elements. 
The fields $K$ we mainly deal with have a small degree $d$ over
$\Q$.
Therefore, we decided to store the multiplication tables
explicitly, since the cost of storing is relatively small compared 
to the speedup we
gain in performing the multiplication by the method described above.

\section{Character tables}

The character values of finite groups are sums of roots of unity and are
therefore algebraic integers of abelian number fields. As described above
the character values are stored as coefficient vectors corresponding
to certain integral bases, which are to be described below.
We shall see that this can be done in such a way that
the resulting  character table is still a square matrix.
The new format of a character table is  based upon the following considerations.

We begin by introducing some notation. If $K \subseteq L$ is a Galois
extension of fields, we let ${\cal G}(L/K)$ denote its Galois group.
Let $n$ be a positive integer. We write $\Q_n$ for
the $n$-th cyclotomic field, and $\varphi(n)$ for the degree of
$\Q_n$ over~$\Q$, i.e, $\varphi$ is the Euler function. 
If $f$ is an integer coprime to $n$ let $*f$
denote the Galois automorphism of $\Q_n$ mapping an $n$-th 
root of unity to its $f$-th power.
 
Let $G$ be a finite group and $x \in G$\@.
Define  ${\Q}(x) = {\Q} (\chi(x) \mid \chi \in \Irr(G))$ to be
the column field of x.
If $x$ has order $n$ then $\Q(x)$ is contained in $\Q_n$\@. 

\begin{definit}
{\rm
Two elements of $G$ are called algebraically conjugate if
the cyclic subgroups they generate are conjugate in $G$. The
equivalence class of $x \in G$  with respect to this relation is called 
the algebraic conjugacy class of $x$ and denoted by $Cl_a(x)$\@. 
Note that it is a union of conjugacy
classes of $G$.
}
\end{definit}
Let $y$ be an element of $Cl_a(x)$. Then $y$ is conjugate to $x^f$
with $f$ relatively prime to $n$.
If $\chi$ is a character of $G$ then $\chi(y) = \chi(x)^{*f}$. 

\begin{prop}
Let $x$ be an element in $G$ of order $n$.
Suppose the column field $\Q(x)$ has degree $d$ 
over $\Q$. Then $d$ equals the number of (ordinary)
conjugacy classes
in the algebraic conjugacy class $Cl_a(x)$\@.
\end{prop}
\begin{bew}
We put $Y = \{ x^i \mid 1 \leq i \leq n, \gcd(i,n) = 1 \}$ 
and write~$S_Y$ for the group of permutations of~$Y$\@.
The normalizer $N = N_G(\langle x \rangle)$ acts on~$Y$ by conjugation. 
The kernel of the corresponding homomorphism
\[\rho: N \rightarrow S_Y \]
is $C_G(x)$ and so $|\rho(N)| = |N/C_G(x)|$\@.

The orbits of $N$ on $Y$ are in one to one correspondence
with the $G$-conjugacy classes intersecting~$Y$ non-trivially.
The size of the orbits is $|N/C_G(x)|$ and the number
of orbits is $\varphi(n)/|N/C_G(x)|$\@.

Define the monomorphism 
\[ \psi: {\cal G}(\Q_n /\Q) \rightarrow S_Y, \]
\[ *k \mapsto (y \mapsto y^k)\] 
and observe that
$\psi({\cal G}(\Q_n/\Q(x))) = \rho(N)$\@. This follows from the remark
preceding the proposition and the fact that two elements
$a, b \in G$ are conjugate, if and only if $\chi(a) = \chi(b)$
for all $\chi \in \Irr(G)$\@. 
Hence $d = |\Q(x) : \Q| = \varphi(n)/|\rho(N)|$ and the assertion
follows.
\end{bew}

\noindent We are now ready to define the {\sf MOC}-character table.
\begin{definit}
{\rm 
\label{MOCtable}
Let $x_1, \ldots,x_r$ be representatives for the algebraic conjugacy classes
of $G$. Let $d_i$ denote the degree of $\Q(x_i)$ over $\Q$\@.
For each~$i$, we choose an integral basis $b_1(x_i), \ldots,b_{d_i}(x_i)$
of $\Q(x_i)$\@.

The columns of the {\sf MOC}-table are indexed by 
\[b_1(x_1),\ldots,b_{d_1}(x_1),b_1(x_2), \ldots,b_{d_2}(x_2),\ldots,
b_1(x_r),\ldots,b_{d_r}(x_r)
\]
 and the entries of $\chi \in \Irr(G)$
at the columns $b_1(x_i), \ldots,b_{d_i}(x_i)$ are 
given by the coefficients of
the decomposition of $\chi(x_i)$ in the basis $b_1(x_i), \ldots, b_{d_i}(x_i)$.
By the preceding proposition the {\sf MOC}-character table
is a square integral matrix.
}
\end{definit}
\begin{beisp}
{\rm
Let $A_5$ be the alternating group on five letters. The ordinary character
table is given by :

\begin{verbatim}
      CAS                               MOC            
                 
           2   2  2  .  .  .     
           3   1  .  1  .  .     
           5   1  .  .  1  1     
                
              1a 2a 3a 5a 5b        
          2P  1a 1a 3a 5b 5a     
          3P  1a 2a 1a 5b 5a      
          5P  1a 2a 3a 1a 1a            1  2  3  5  5
           2          
      X.1  +   1  1  1  1  1            1  1  1  1  0
      X.2  +   3 -1  .  A *A            3 -1  0  1  1
      X.3  +   3 -1  . *A  A            3 -1  0  0 -1
      X.4  +   4  .  1 -1 -1            4  0  1 -1  0
      X.5  +   5  1 -1  .  .            5  1 -1  0  0
   
\end{verbatim}
Here, $A=(1+\sqrt{5})/2$ and $*A=(1-\sqrt{5})/2$\@.
The integral basis of $\Q(5A)$ in the {\sf MOC}-table
is $\{1, \epsilon + \epsilon^4\}$, where
$\epsilon= e^{2\pi i/5}$\@. 
}
\end{beisp}
We now indicate how a usual character table can be transformed into
a {\sf MOC}-character table and vice versa.
\begin{definit}
{\rm Let ${\cal G}(\Q(x_i)/\Q) = 
\{ \sigma_{i,1},\ldots,\sigma_{i,d_i} \}$,
where the notation is as in Definition \ref{MOCtable}\@. 
Let $A_i$ be the $(d_i \times d_i)$-matrix whose $(j,k)$-th entry is
\[
\sigma_{i,j}(b_k(x_i)).
\]
Let $d = \sum_{i=1}^r d_i$ and let
$A$ be the $(d \times d)$-block diagonal matrix whose $i$-th block
is $A_i$.} 
\end{definit}

\noindent Note that $(\Det(A_i))^2$ equals the discriminant 
of $\Q(x_i)$ over $\Q$; in particular, $A$ is invertible.
With a suitable ordering of the classes, we obtain the usual
character table from the {\sf MOC}-character table by a matrix multiplication:
\[
\mbox{\rm  Usual ordinary table} = \mbox{{\sf MOC}-{\rm Table}} \cdot A.
\]
Note that this observation can be used to calculate the block distribution
of the characters by using the {\sf MOC}-table as follows.  
Let $\chi \in \Irr(G)$ and let $\omega_\chi$ be the corresponding 
central character.
Let $c^i_\chi$ denote the row of coefficients of $\omega_\chi(x_i)$
expressed in the integral basis $b_1(x_i), \ldots,b_{d_i}(x_i)$\@.
Then 
\[
[ \omega_\chi] = [c^1_\chi, \ldots , c^r_\chi] \cdot A.
\]
Since the block matrices $A_i$ corresponding to $p$-regular classes have
determinant not divisible by $p$ 
\cite[Lemma I.(10.1), Satz I.(2.11), Korollar III.(2.10)]{neukirch}, 
two characters $\chi$ and $\psi$ are in the same $p$-block of $G$, if and
only if 
\[
c_\chi^i \equiv c_\psi^i (\mod\, p) 
\] 
for all $i$ corresponding to $p$-regular classes. We thus have proved
the stronger result on block distribution (see \cite[(7.10)]{gold}), namely
that two characters lie in the same $p$-block if and only if the
values of their central characters on $p$-regular classes differ
by elements of $pR$\@. 

The {\sf MOC}-system contains a data base for the
subfields $L$ of cyclotomic fields up to a certain degree.
In this data base a field $L$ is  stored
by the following information:
The minimal $n$ such that $L$ is contained in 
the cyclotomic field ${\Q}_n$,
the degree $d$
of $L$ over ${\Q}$, generators for the  Galois group
${\cal L} = {\cal G}({\Q}_n / L)$ and an integral
basis of $L$.
Each basis element is a sum over the elements in an orbit of
${\cal L}$ acting on the $n$-th roots of unity. Such elements are 
called orbit sums in the following. 
It is only necessary to store one representative for each orbit of
${\cal L}$ whose corresponding sum is a basis element.

Only recently it has been shown by
H.W. Lenstra that there always exists such a basis
for a given abelian extension over ${\Q}$ \cite{lenstra}\@.
We are grateful to him for allowing us to reproduce this previously
unpublished result. It is given in the next section.

\section{Lenstra's theorem on integral bases}

We keep the notation of the preceding section.
Furthermore, we put $\zeta_m = \exp(2\pi i/m)$\@.

\begin{satz}
{\rm (H.\ W.\ Lenstra~\cite{lenstra}:)} Let~$L/\Q$ be a finite abelian extension
of the rational numbers. Let~$f$ be the smallest positive integer
such that~$L$ is contained in the $f$-th cyclotomic field $\Q_f$
(the conductor of~$L$)\@. Finally let~$H$ be the subgroup of 
the Galois group ${\cal G}(\Q_f/\Q)$ which fixes~$L$\@.
Then there is an integral basis of~$L$,
which consists of orbit sums of $H$ on $f$-th roots of unity.
\end{satz}
\begin{bew}
Let~$p$ be a prime and $k \in \N$\@. We arrange the set of $p^k$-th 
roots of unity into the following matrix:
$$\begin{array}{ccccc}
1 & \zeta_p & \zeta_p^2 & \ldots & \zeta_p^{p-1} \\
\zeta_{p^k} & \zeta_{p^k}\zeta_p & \zeta_{p^k}\zeta_p^2 & \ldots 
& \zeta_{p^k}\zeta_p^{p-1} \\
\zeta_{p^k}^2 & \zeta_{p^k}^2\zeta_p & \zeta_{p^k}^2\zeta_p^2 & \ldots 
& \zeta_{p^k}^2\zeta_p^{p-1} \\
\multicolumn{1}{c}{\vdots} & \multicolumn{1}{c}{\vdots} & 
\multicolumn{1}{c}{\vdots} & \dots & \multicolumn{1}{c}{\vdots} \\
\zeta_{p^k}^{p^{k-1} - 1} & \zeta_{p^k}^{p^{k-1} - 1}\zeta_p 
& \zeta_{p^k}^{p^{k-1} - 1}\zeta_p^2 & \ldots 
& \zeta_{p^k}^{p^{k-1} - 1}\zeta_p^{p-1} \\
\end{array}$$
We introduce an equivalence relation $\sim$ on 
$T = \langle \zeta_{p^k} \rangle \setminus \{1\}$ as follows:
For $\eta, \eta' \in T$
we put $\eta \sim \eta'$ if and only if either $\eta = \eta'$ or else
$\eta^p = {\eta'}^p \neq 1$\@. 
The equivalence classes are the one element sets $\{\zeta_p^i\}$, 
$1 \leq i \leq p - 1$ and the $p$-element sets
constituted by the elements of a row
different from the first of the above matrix.
If~$C$ is one of the equivalence classes containing more than
one element and if $\eta \in C$, then 
$$\sum_{\mu \in C} \Z\mu = \eta\,\Z[\zeta_p],$$
and thus is a free $\Z$-module of rank $p - 1$\@.
By counting dimensions, we find
\begin{equation}
\label{zerlegung1}
\Z[\zeta_{p^k}] = 
\bigoplus_{C \in T/\sim} 
\left( \sum_{\eta \in C} \Z\eta \right).
\end{equation}
We now generalize these observations.
Write $$f = \prod_{p\,\,{\rm prime}} p^{k(p)}.$$ Then
\begin{equation}
\label{tensor1}
\Z[\zeta_f] \cong \bigotimes_{p\,\,{\rm prime}} \Z[\zeta_{p^{k(p)}}].
\end{equation}
From this and the decomposition~(\ref{zerlegung1})
we obtain a similar direct decomposition of $\Z[\zeta_f]$
which we are now going to describe.
Let $f_0$ be the product of the different
prime divisors of~$f$\@. Put $S = \{ \eta \in \langle \zeta_f \rangle \mid f_0
\mbox{\rm\ divides } |\eta|\}$\@. 
For $\eta \in S$ let $d(\eta)$ denote the largest squarefree number such
that $d(\eta)^2 \mid |\eta|$\@.
We introduce an equivalence relation~$\sim$ on~$S$ as follows:
For $\eta, \eta' \in S$ we put $\eta \sim \eta'$ if and only if
$\eta^{d(\eta)} = {\eta'}^{d(\eta)}$\@.
If $\eta \sim \eta'$, then $|\eta| = |\eta'|$ and so our relation is
indeed symmetric\@.
From~(\ref{zerlegung1}) and~(\ref{tensor1}), we obtain
\begin{equation}
\label{zerlegung2}
\Z[\zeta_f] = \bigoplus_{S/\sim} \left( \sum_{\eta \in C} \Z\eta \right).
\end{equation} 
If $C$ is an equivalence class, $\eta \in C$ and $d = d(\eta)$,
then $\sum_{\mu \in C} \Z\mu = \eta \Z[\zeta_d]$,
and thus  the $\Z$-rank of $\sum_{\mu \in C} \Z\mu$ equals $\varphi(d)$\@.
Thus there is a subset $B_C \subset C$ which is a $\Z$-basis of
$\sum_{\mu \in C} \Z\mu$, e.g.,
$B_C = \{ \eta \zeta_d^i \mid 0 \leq i < \varphi(d)\}$, or
$B_C = \{ \eta \zeta_d^i \mid \GCD(i,d) = 1\}$\@.

It is clear that~$H$ fixes~$S$ and permutes the equivalence classes.
Suppose first that~$H$ acts fixed point freely on the set of 
equivalence classes of~$S$, i.e., no non-identity element of~$H$
fixes an equivalence class. Then, by~(\ref{zerlegung2}) and the subsequent
remarks, it is clear that we may choose an integral basis 
of $\Z[\zeta_f]^H$, the set of integers in~$L$,
consisting of orbit sums of~$H$ on~$S$\@.
Thus we are done in this case.

Suppose there is some $1 \neq \sigma \in H$ and $\eta \in S$
such that $\sigma(\eta) \sim \eta$\@. Let~$d = d(\eta)$\@.
Then, by definition,
$\sigma(\eta^d) = \sigma (\eta)^d = \eta^d$\@. 
The order of $\eta^d$ is divisible
by $f_0$, 
and so $\sigma (\zeta_{f_0}) = \zeta_{f_0}$\@. This implies
$\sigma \in H_0$, where $H_0$ is the intersection of~$H$
with the kernel~$M$ of the natural map
$$(\Z/f\Z)^* \rightarrow (\Z/f_0\Z)^*.$$
The order of~$M$ is $\prod_p p^{k(p) - 1}$, and the
Sylow $p$-subgroup is cyclic, except if $p = 2$, and $k(2) \geq 3$\@.
Let $f' = f/2^{k(2)}$ denote the odd part of~$f$\@.
If $k(2) \geq 3$, the Sylow $2$-subgroup of~$M$ is 
$M_0 \times M_1$, with $M_0 = \{*a \mid a \equiv 1(\mod\,f') \mbox{\rm\ and } 
a \equiv \pm 1(\mod\,2^{k(2)})\} $ and 
$M_1 = \{*b \mid b \equiv 1(\mod\,4f')\}$\@.

The condition that~$f$ be the conductor of~$L$ is equivalent to
saying that for all~$p$ dividing~$f$, the group~$H$ does not
contain the kernel of the natural map
\begin{equation}
\label{fuehrer}
(\Z/f\Z)^* \rightarrow (\Z/\frac{f}{p}\Z)^*.
\end{equation}
This kernel has order~$p$ and is generated by $*(1 + f/p)$\@.

So if $p$ is odd or $p = 2$ and $k(2) \leq 2$, then the Sylow
$p$-subgroup of $H_0$ is trivial, since otherwise 
$*(1 + f/p)$ would be contained in $H$,
which is not the case by the remarks above. 
Thus $H_0$ is a $2$-group, $k(2) \geq 3$ and so $H_0 \leq M_0 \times M_1$\@. 
Let $a \in \Z$ such that $*a \in H_0$ has order~$2$\@. 
Then $a \equiv 1 (\mod\,f')$ and one of the following holds:
$$ a \equiv \left\{ \begin{array}{l} -1 \\ -1 + 2^{k(2) - 1} \\
1 + 2^{k(2) - 1} \end{array}\right\} (\mod\,2^{k(2)}).$$
By the remarks preceding and following~(\ref{fuehrer}), 
$H_0 \cap M_1 = 1$, and thus the last possibility cannot occur 
and we must have $H_0 = \langle a \rangle$, i.e., $H_0$ is a cyclic
group of order~$2$\@.

So if some non-trivial element $\sigma_0$ of $H$ fixes the equivalence
class containing~$\eta$, then $\sigma_0 = *a$\@. Furthermore,
if $2$ is the exact power of~$2$ dividing $|\eta|$, then
$\sigma_0 (\eta) = \eta$, and if $4$ is the exact power of~$2$
dividing $|\eta|$, then $\sigma_0 (\eta) = - \eta$\@.
On the other hand, if $\sigma (\eta) = \eta$ or $\sigma (\eta) = - \eta$
for $\sigma \in H$ then $\sigma (\eta) \sim \eta$\@.
Suppose that $8$ divides the order of $\eta$\@. Write
$d(\eta) = 2c$, $c$ odd. Then $\eta \sim \sigma_0 (\eta)
= \eta^a$ implies $\eta^{2c} = \eta^{2ca}$, and so $\eta^{2c(a - 1)} = 1$\@.
However, $2c(a - 1)$ is not divisible by~$8$, a contradiction.

Let $H_1 = \{ *a \in H \mid a \equiv 1 (\mod 4)\}$\@. 
Then $H_0 \cap H_1 = 1$, and so $H = H_0 \times H_1$\@.
We now proceed as follows. Since $H_1$ acts fixed point freely
on the equivalence classes, we can choose an integral basis
of $\Z[\zeta_f]^{H_1}$ consisting of orbit sums of $H_1$ on~$S$\@.
If an orbit sum contains an element from an equivalence class
fixed by~$H_0$, we consider two cases. Either the element has order
divisible by~$4$, then $H_0$ acts as~$-1$ on the corresponding
$\Z H_1$-module. Otherwise, $H_0$ fixes every element in the orbit,
in which case we include the corresponding orbit sum, which is now 
an $H$-orbit sum, into our basis. 

The remaining orbit sums contain only elements of equivalence classes
not fixed by~$H_0$\@. The corresponding $H$-orbit sums are
taken to yield the desired integral basis.
\end{bew}

\clearpage
\thispagestyle{empty}

\chapter{Algorithms}
\label{improve}

In this chapter we describe the principal algorithms
utilized by {\sf MOC}\@.
We first explain two algorithms for dealing with systems of integral linear 
equations using $q$-adic approximation. Then we describe
the methods of integral linear programming 
used to attack the two fundamental problems. Furthermore we
show how all {\sf MOC}-computations can be documented 
and how these documentations are used to give a more
conventional proof for the results obtained.
Finally we
discuss some more advanced methods which can be applied if the
{\sf MOC}-system does not give the complete answer.

\section{Integral linear equations}

Most of the time during a run of {\sf MOC} we are dealing with
the problem of expressing a character as a $\Z$-linear combination
in terms of a basic set. In other words we have to solve a system
of integral linear equations. 
Dixon \cite{Dix82} independently suggested the following way of 
solving this problem. 

\begin{algorithm}
\label{DEC}
{\rm (DEC) 

\noindent {PROBLEM:} Given $m+1$ rows of integers 
\[
w, b_1,b_2, \ldots, b_m \in \Z^{1 \times n},
\]
where $b_1,b_2, \ldots, b_m$ are linearly independent over $\Z$.
If possible write $w$ as a $\Z$-linear combination of these.

\noindent {SOLUTION:} Let $T$ be the matrix 
\[
T = \left( \begin{array}{c}
                   b_1 \\
                   \vdots \\
                   b_m
            \end{array} \right)
\]
Let $q$ be a (large) prime such that $\bar T$, the reduction modulo $q$
of $T$, has rank $m$. Set $v:=w$, $j:=0$ and $w_0:=0$ (zero vector).
Set {MAXJ} to a reasonable positive integer, representing the 
maximum number of loops, we are willing to spend our time with.

\begin{description}
\item{{\bf Step 1.}} Check whether $\bar v$, the reduction modulo $q$
of~$v$, is in the $\F_q$-span of
$\bar b_1,\bar b_2, \ldots ,\bar b_m$.
If not, STOP, noting that $w$ is not in the $\Q$-span of
$b_1,b_2, \ldots, b_m$. Otherwise go to Step~$2$\@.

\item{{\bf Step 2.}} 
Define the integers $z_{ij}$,
$1 \le i \le m$ by
\[
\bar v = \sum_{i=1}^{m} \bar z_{ij} \bar b_i, \quad \mbox{where} \quad
-(q-1)/2 \le z_{ij} \le (q-1)/2,
\]
and set 
\[
w_{j+1} := \sum_{s=0}^{j}(z_{1s}, \ldots,z_{ms})q^sT.
\]
Then every coefficient of $w - w_{j+1}$ is divisible by $q^{j+1}$. 
If $w - w_{j + 1} = 0$, then STOP\@. Else set 
\[
v:= {(w-w_{j+1})/q^{j+1}}
\]
and increase $j$ by $1$. If $j$ is larger than {MAXJ} then STOP,
otherwise continue with Step 1.
\end{description}
}
\end{algorithm}
\begin{rem}{\rm
Note that $w$ is in the ${\Q }$-span  of
$b_1,b_2, \ldots , b_m$ if and only if
\begin{equation}
\label{padic}
w = \sum_{j=0}^{\infty} (z_{1,j}, \ldots, z_{m,j})q^jT
\end{equation}
in the $q$-adic metric, where the integers $z_{ij}$ are obtained by
the above algorithm. Furthermore $w$ is in the $\Z$-span of
$b_1,b_2, \ldots , b_m$ if and only if (\ref{padic}) is finite.
When the algorithm terminates because $j$ gets larger then {MAXJ}, this may
have three reasons.
\begin{itemize}
\item[(a)]{ $w$ is in the $\Z$-span but some coefficient is too large.}
\item[(b)]{$w$ is in the $\Q$-span but not in the $\Z$-span.}
\item[(c)]{$w$ is not in the $\Q$-span.}
\end{itemize}
Our experience tells us that $q = 101$ and $\mbox{\rm MAXJ} = 20$ are 
reasonable choices for the problems we are concerned with.
}
\end{rem}

The first step in the computation of the decomposition numbers
is to find a basic set for $G_0(FG)$. We start with the restrictions
of the ordinary irreducible characters 
to the  $p$-regular elements which form a generating set for
$G_0(FG)$. We then have to derive from these 
a basic set of $G_0(FG)$. There are algorithms for constructing 
bases of $\Z$-lattices from generating sets, but in our situation we
have the additional restriction that our basis has to consist of
proper Brauer characters.  
This is achieved by an algorithm due to Parker, called FBA, 
which is based upon Algorithm \ref{DEC}.
We describe  FBA in the following general setup. 

\begin{algorithm}
\label{FBA}
{\rm (FBA)

\noindent {PROBLEM:} Let $X$ be a
$\Z$-lattice with finite basis $\varphi_1, \ldots, \varphi_m$.
Define $X^+$ to be the subset of $X$ consisting of all non-negative
linear combinations of the basis elements $\varphi_1,\ldots,\varphi_m$.
Suppose that we are given a finite set 
${\bf B}=\{\vartheta_1,\ldots,\vartheta_n\} \subset X^+$
generating a sublattice $Y$ of $X$.
Derive from ${\bf B}$ a basis of $Y$ lying in $X^+$.

\medskip

\noindent {SOLUTION:} Assume inductively that a subset 
$\tilde{B}=\lbrace\beta_1, \ldots, \beta_{k-1}\rbrace \subset X^+$
has been found with span equal to $\langle \vartheta_1,
\ldots, \vartheta_s \rangle_{\Z}$, $s \le n$, and such that
$\beta_1, \ldots, \beta_{k-1}$
are linearly independent modulo some large prime $q$.

\begin{description}
\item{{\bf Step 1.}} If $s=n$ then STOP and output $\tilde{B}$. Otherwise
take the next element 
$\vartheta_{s+1} \in {\bf B}$
and try to express it in terms of $\beta_1, \ldots,\beta_{k-1}$ using DEC with
the prime $q$. The following cases can occur and are checked in 

\pagebreak

\item{{\bf Step 2.}}
\begin{itemize}
\item[(a)]{If $\vartheta_{s+1} \in \langle \beta_1, \ldots ,\beta_{k-1}
\rangle_{\Z}$, increase $s$ by~$1$ and go to Step 1.
}
\item[(b)]{ If $ \beta_1, \ldots, \beta_{k-1},
\vartheta_{s+1}$ are linearly independent modulo $q$,
increase $s$ and $k$ by $1$, replace $\tilde{B}$
by $\lbrace\beta_1, \ldots, \beta_{k-1},\vartheta_{s+1}\rbrace$ and
go to Step~$1$\@.
}
\item[(c)]{ If $\beta_1, \ldots, \beta_{k-1},
\vartheta_{s+1}$ are linearly dependent modulo $q$ but linearly independent
over $\Z$, replace $q$ by a prime $\ell$ such that
$\beta_1, \ldots,\beta_{k-1}, \vartheta_{s+1}$ are linearly independent
modulo $\ell$\@.
Replace $\tilde{B}$
by $\lbrace\beta_1, \ldots, \beta_{k-1},\vartheta_{s+1}\rbrace$,
increase $s$ and $k$ by $1$ and go to Step~$1$\@.
}
\item[(d)]{ If $\beta:=\vartheta_{s+1} \in \langle \beta_1, \ldots, \beta_{k-1} \rangle_{\Q}$
but  $\beta \not\in \langle \beta_1, \ldots , \beta_{k-1} \rangle_{\Z}$,
we use the first terms of the $q$-adic expression for the 
coefficients of $\beta$ determined by DEC to find the
smallest possible natural number $m_1$ with

\[
m_1 \beta = \sum_{i=1}^{k-1} z_i \beta_i,
\]
where $z_i \in \Z$. 
Let $m = \min\lbrace \gcd(m_1,z_i) \mid 1 \le i \le k-1\rbrace$ and 
suppose $m = \gcd(m_1,z_j) > 0$\@. 
Write $ m = a m_1 + b z_j$, $a, b \in \Z$\@.
For $i$ different from $j$, let $c_i$ be the smallest integer
such that $c_i + b z_i m_1^{-1} \ge 0$ and  define
\[
\tilde{\beta_j} = b \beta + a \beta_j + \sum^{k - 1}_{i \ne j} c_i \beta_i.
\]
Then $m_1 \tilde{\beta_j}$ is contained in $X^+$ and hence
$\tilde{\beta_j} \in X^+$, too.
Replace $\beta_j$ by $\tilde{\beta_j}$ in $\tilde{B}$\@.  
Increase $n$ by~$1$ and put $\vartheta_{n+1} = \beta_j$\@.
Do not increase $s$ in this case. Go to Step~$1$\@.
}
\end{itemize}
\end{description}
}
\end{algorithm}

\begin{satz} {\rm Algorithm \ref{FBA}} terminates and finds a 
$\Z$-basis of~$Y$ contained in $X^+$\@.
\end{satz}
\begin{bew}
For any subset $B \subseteq {\bf B}$ we define the weight 
$w(B)$ of~$B$ to be the pair $(\mbox{\rm cr}(B),d(B))$\@.
Here, $\mbox{\rm cr}(B)$ is the co-rank of the $\Z$-span of~$B$
in~$Y$\@. The number $d(B)$ is the discriminant of $\langle B \rangle_\Z$, 
which is defined in the following way.
Choose any $\Z$-basis $b_1, \ldots , b_s$ of $\langle B \rangle_\Z$, 
and let~$M$ denote the $(m \times s)$-matrix
of the coefficients of $b_1, \ldots , b_s$ expressed 
in the basis $\varphi_1, \ldots , \varphi_m$ of~$X$\@.
Then $d(B) = \Det\,MM^t$\@. This number is positive
and independent of the basis chosen for $\langle B \rangle_\Z$\@. 

We order the set of weights $W({\bf B}) = \{w(B) \mid B \subseteq {\bf B}\}$ 
lexicographically. Observe that $B$ generates~$Y$ if and only if
$w(B)$ is minimal in $W({\bf B})$\@. 

In each of the cases occuring in Step~$2$ of Algorithm~\ref{FBA},
where we change~$\tilde{B}$, we decrease its weight.
This is clear in cases~(b) and~(c), since we increase the
rank of the generated lattice. 
In case~(d) let~$S$ denote the matrix transforming 
$\{\beta_1, \ldots , \beta_j, \ldots , \beta_{k-1} \}$ into 
$\{\beta_1, \ldots , \tilde{\beta}_j, \ldots , \beta_{k-1} \}$\@. 
Then~$S$ is a diagonal matrix with diagonal entries equal to~$1$,
except that row~$j$ of~$S$ has the entries
$$c_1 + bz_1/m_1, \ldots, c_{j-1} + bz_{j-1}/m_1, a + bz_j/m_1, $$
$$c_{j+1} + bz_{j+1}/m_1, \ldots , c_{k-1} + bz_{k-1}/m_1.$$ 
Thus $\Det (S) = a + bz_j/m_1 = m/m_1 < 1$\@.
This implies that the new set $\tilde{B}$ has a smaller discriminant.

Finally, by adding $\beta_j$ to ${\bf B}$ the minimum of
$W({\bf B})$ remains unchanged, since the new enlarged ${\bf B}$
still generates the same lattice~$Y$\@. 
This proves that the algorithm terminates. Since in
each step the set $\tilde{B}$ is linearly independent over the
integers, the final such set is a $\Z$-basis of~$Y$\@.
\end{bew}

\section{Improving basic sets}
\label{Improvingbasicsets}

Having established the major part of the theory we are now prepared to give
some algorithms for dealing with the two fundamental problems stated in
Sections \ref{fundamentalproblemI} and \ref{fundamentalproblemII}\@.
These algorithms are implemented in the {\sf MOC}-system
and are parts of the program IMPROVE, which
uses methods of linear programming. 
There are some remarks concerning optimization at the end of this section. 

Let us 
fix some notation. Let ${\bf PS}  = \{ \Phi_1,
\ldots, \Phi_s \} $ be the 
basic set of projectives and ${\bf BS} =
\{ \varphi_1, \ldots, \varphi_s \}$ the 
basic set of Brauer characters.
We denote by ${\bf PA} = \{ \varphi_1^*, \ldots, \varphi_s^*\}$ and ${\bf BA} =
\{\Phi_1^*, \ldots, \Phi_s^*\}$ the dual systems of atoms. Let 
${\bf P}$ and ${\bf B}$ denote additional
sets of projective respectively Brauer characters. 
We denote by $U = (\langle \varphi_i, \Phi_j \rangle_
{ij})$ the matrix of scalar products between ${\bf BS}$ and ${\bf PS}$. 

We are now going to restate the first fundamental problem. 
Recall that we have to find
all square integral matrices $U_1 \geq 0$ and $U_2 \geq 0$
so that $U = U_1 \cdot U_2$ subject to the conditions
(\ref{nebenbed})\@. We consider the
equivalent problem of finding $U_1 \geq 0$ and $U_2 \geq 0$ invertible over
$\Z$ with $U_1^{-1} \cdot U \cdot U_2^{-1} = E_s$. 
Operating by $U_1^{-1}$ means changing ${\bf BS}$ and operating
by $U_2^{-1}$ means changing ${\bf PS}$. 
We are going to change $U$ according
to (\ref{nebenbed}) hoping to get $E_s$ as result.
We have the two following problems to solve:
\begin{enumerate}
\item Prove that $\varphi_i$ is irreducible or that $\Phi_j$ is indecomposable.
This is considered in~\S~\ref{pimtestalg}\@.
\item Improve ${\bf PS}$ or ${\bf BS}$. We shall show that this can be
done in several ways. We present a theorem in~\S~\ref{subtractalg}
which allows to subtract 
irreducible characters from other ones. 
Then in~\S~\ref{triangularalg} there is an algorithm which
computes characters lexicographically smaller. 
We finally present a criterion in~\S~\ref{decomposablealg}
for a character to be decomposable or reducible.
\end{enumerate}

\subsection{Proving indecomposability}
\label{pimtestalg}

In this paragraph we are going to establish a test for proving irreducibility
of a Brauer 
character respectively indecomposability of a projective character.
First of all we need the following

\begin{definit}
\label{partdef}
{\rm 
Let
$$\Phi = \sum_{i=1}^s n_i \varphi^*_i \quad \mbox{ with } \quad n_i \geq 0$$
be a projective
character decomposed into ${\bf PA}$. A generalized projective character 
$$\Phi' = \sum_{i=1}^s n_i' \varphi^*_i$$
is called a {\sl part of $\Phi$} if
\index{part of a character}
$0 \leq n_i' \leq n_i$ for all $i$. A similar notion is used for Brauer
characters.
}
\end{definit}
The projective indecomposable summands of a projective character~$\Phi$
are among the parts of~$\Phi$\@. 
We are able to state the main criterion for indecomposability.

\begin{satz}
\label{pimtest}
{\rm (PIM-test:)} The projective
$\Phi$ is indecomposable if for all parts $\Phi'$ of $\Phi$ with $\Phi' \neq
\Phi$ and $\Phi' \neq 0$ there is a Brauer character $\varphi \in {\bf B}$
with $\langle \varphi, \Phi' \rangle < 0$ or $\langle \varphi, 
\Phi - \Phi' \rangle
< 0$. 
\end{satz}
\begin{bew}
If there is a $\varphi$ having negative scalar product with $\Phi'$ or
with $\Phi - \Phi'$ then one of $\Phi'$,~$\Phi - \Phi'$
is not a genuine projective. If $\Phi'$
is a part of~$\Phi$, then so is $\Phi - \Phi'$\@.
\end{bew}

\noindent Sometimes one can find another test for indecomposability
in the literature.
It is based on the fact that an ordinary character is a generalized projective
character if and only if it is $0$ on all $p$-singular classes. Let $\Phi =
\sum_{i=1}^n a_i \chi_i$ be a projective character expressed as a sum of
ordinary characters. Then we call a character $\Phi' = \sum_{i=1}^n a_i'
\chi_i$ a {\sl subsum of $\Phi$} if $0 \leq a_i' \leq a_i, \Phi' \neq 0$
\index{subsum of a character}
and $\Phi' \neq \Phi$. In the test in question we consider 
all subsums of $\Phi$ and
check if one of them is a generalized projective. If there is no such subsum
then $\Phi$ is indecomposable. Now we are going to compare the two tests.

\begin{prop}
Suppose that ${\bf BS}$ is a special basic set and that ${\bf B}$
consists of the restricted ordinary characters.
Let $\Phi$ be a projective character. Then the PIM-test
presented in {\rm Theorem~\ref{pimtest}}
proves that $\Phi$ is indecomposable, 
if and only if there is no
subsum of $\Phi$ which is a generalized projective.
\end{prop}
\begin{bew}
Let $\Psi$ be an arbitrary generalized projective character and
$$\Psi = \sum_{i=1}^n a_i \chi_i$$
its decomposition into the ordinary irreducible characters. Then we have
$\langle \hat{\chi}_i , \Psi \rangle = a_i$. 
Let $\Phi'$ and $\Phi'' = \Phi - \Phi'$
be parts of $\Phi$. Then $\Phi'$ and $\Phi''$ are generalized projective
characters. Suppose there is no subsum of $\Phi$ which is a generalized
projective character. Then neither $\Phi'$ nor $\Phi''$ is a subsum. Therefore
at least one of them, say~$\Phi'$, has a negative coefficient 
in its decomposition into
ordinary characters. Call this coefficient $a_i$. But $a_i$ is the scalar
product of~$\Phi'$ with $\hat{\chi}_i$, as we have seen above.
Hence the PIM-test of Theorem~\ref{pimtest} proves the 
indecomposability of~$\Phi$\@.

Now suppose that $\Phi'$ is a subsum of $\Phi$ and $\Phi'$ is a generalized
projective character. Then $\langle \hat{\chi}_i, \Phi' \rangle \geq 0$ 
and $\langle \hat{\chi}_i, \Phi - \Phi' \rangle \geq 0$ for all
$1 \leq i \leq n$. Because of the choice of ${\bf BS}$ this particular
$\Phi'$ is checked by our PIM-test as well. 
So it can not decide whether $\Phi$ is indecomposable or not.
\end{bew}

\noindent By adding more Brauer characters to ${\bf B}$ 
(for example Brauer characters
got by induction) our PIM-test may find more PIMs.

\begin{rem}
\label{gomoryhypo1}
\rm
Let $\Phi = \sum_{i=1}^s n_i \varphi_i^*$, $|{\cal B}| = m$
and denote by $v_{k1}, \ldots, v_{ks}$ 
the coefficients of the decomposition
of the $k$-th character in ${\bf B}$, $1 \leq k \leq m,$ when decomposed into
${\bf BS}$. Then $\Phi$ is indecomposable if there are
no $n_1', \ldots, n_s'$ satisfying the following conditions:
\medskip

$$0 \leq n_i' \leq n_i \quad \mbox{ for all } \quad 1 \leq i \leq s,$$
\medskip

$$\quad \sum_{i=1}^s n_i' \leq \sum_{i=1}^s n_i - 1 \quad \mbox{ (that is } \Phi'
\neq \Phi),$$
\medskip

$$\sum_{i=1}^s n_i' \geq 1 \quad \mbox{ (that is } \Phi' \neq 0),$$
\medskip

$$\sum_{j=1}^s - v_{kj} n_j' \leq 0 \quad \mbox{ for all } \quad 1 \leq k \leq m
\quad \mbox{ (that is } \langle \varphi, \Phi' \rangle \geq 0 \mbox{ for } \varphi \in {\bf B}),$$
\medskip

$$\sum_{j=1}^s v_{kj} n_j' \leq \sum_{j=1}^s v_{kj} \langle 
\varphi_j, \Phi \rangle
\mbox{ for all } 1 \leq k \leq m
\mbox{ (that is } \langle \varphi, \Phi - \Phi' \rangle \geq 0).$$
This system of inequations can be written in matrices as follows.
There are no $n_1', \ldots, n_s'$ satisfying $A \cdot x \leq b$ where
$$A = \left( \begin{array}{ccc} 1 & & 0 \\ & \ddots & \\ 0 & & 1 \\
1 & \ldots & 1 \\ -1 & \ldots & -1 \\
-v_{11} & \ldots & -v_{1s} \\ \vdots & & \vdots \\ -v_{m1} & \ldots & -v_{ms} \\
v_{11} & \ldots & v_{1s} \\ \vdots & & \vdots \\ v_{m1} & \ldots & v_{ms}
\end{array} \right) \mbox{ and } b = \left( \begin{array}{c} \langle
\varphi_1 , \Phi \rangle \\ \vdots \\ \langle \varphi_s, \Phi 
\rangle \\ -1 + \sum_{j=1}^s
\langle \varphi_j, \Phi \rangle \\ -1 \\ 0 \\ \vdots \\ 0 \\ \sum_{j=1}^s v_{1j}
\langle \varphi_j, \Phi \rangle \\ \vdots \\ \sum_{j=1}^s v_{mj} \langle
\varphi_j, \Phi \rangle \end{array} \right)$$
and $x = (n_1', \ldots, n_s')^t$. How this system of inequalities can be
solved is discussed in \S~\ref{gomoryint}\@. 
For the reader who is familiar with
linear programming we mention that this system is dual feasible if we use
a dual simplex algorithm and minimize the $0$-function. Then, of course,
the zero vector is a dual solution. Furthermore we observe that in each 
column of $A$ the first non-zero entry is $1$.
\end{rem}

\begin{rem}
\rm
\label{smallsize}
\begin{itemize}
\item[(i)] The number of rows in $A$ depends on $| {\bf B} |$. Therefore it 
is useful to keep $|{\bf B}|$ as small as possible.
\item[(ii)] By interchanging projectives and Brauer characters 
we can test the irreducibility of a Brauer character.
\end{itemize}
\end{rem}

\begin{beisp}
\label{c2mod5}
\rm
In this chapter we are going to calculate the $5$-modular table of the sporadic
simple group {\em Co}$_2$. 
Its ordinary character table can be found in \cite{atlas}
as well as the one of the first maximal subgroup which is isomorphic
to $U_6(2).2$. We make use of the sixth maximal subgroup, too. 
It is isomorphic to $(2_+^{1+6} \times 2^4).A_8$. 
The $5$-modular tables of both maximal 
subgroups can be constructed by making use only 
of tensor products. 

\bigskip

(i) In $U_6(2).2$ we have six blocks of defect~$1$
with the following decomposition matrices:
\[
\begin{array}{rrrrr}
 1^+ &  1 &  . &  . & . \\
 616^- &  . &  1 &  . & . \\
 8064^- &  . &  1 &  1 & . \\
 11264^+ &  1 &  . &  . & 1 \\
 18711_2^- &  . &  . &  1 & 1 \\
& \Psi_4 & \Psi_3 & \Psi_2 & \Psi_1
\end{array}
\quad
\begin{array}{rrrrr}
 1^- &  1 &  . &  . & . \\
 616^+ &  . &  1 &  . & . \\
 8064^+ &  . &  1 &  1 & . \\
 11264^- &  1 &  . &  . & 1 \\
 18711_2^+ &  . &  . &  1 & 1 \\
& \Psi_8 & \Psi_7 & \Psi_6 & \Psi_5 
\end{array}
\]
\[
\begin{array}{rrrrr}
 22^+ &  1 &  . &  . & . \\
 252^- &  . &  1 &  . & . \\
 4928^- &  . &  1 &  1 & . \\
 32768^+ &  1 &  . &  . & 1 \\
 37422^+ &  . &  . &  1 & 1 \\
& \Psi_{12} & \Psi_{11} & \Psi_{10} & \Psi_9
\end{array}
\quad
\begin{array}{rrrrr}
 22^- &  1 &  . &  . & . \\
 252^+ &  . &  1 &  . & . \\
 4928^+ &  . &  1 &  1 & . \\
 32768^- &  1 &  . &  . & 1 \\
 37422^- &  . &  . &  1 & 1 \\
& \Psi_{16} & \Psi_{15} & \Psi_{14} & \Psi_{13}
\end{array}
\]
\[
\begin{array}{rrrrr}
 231^+ &  1 &  . &  . & . \\
 1386^- &  . &  1 &  . & . \\
 5544^+ &  1 &  . &  1 & . \\
 14784^- &  . &  1 &  . & 1 \\
 18711_1^- &  . &  . &  1 & 1 \\
& \Psi_{20} & \Psi_{19} & \Psi_{18} & \Psi_{17}
\end{array}
\quad
\begin{array}{rrrrr}
 231^- &  1 &  . &  . & . \\
 1386^+ &  . &  1 &  . & . \\
 5544^- &  1 &  . &  1 & . \\
 14784^+ &  . &  1 &  . & 1 \\
 18711_1^+ &  . &  . &  1 & 1 \\
& \Psi_{24} & \Psi_{23} & \Psi_{22} & \Psi_{21}
\end{array}
\]
The remaining characters have defect $0$. We write them as
$\Psi_{25}, \ldots, \Psi_{59}$, where their order agrees 
with the one in \cite{atlas}.

\bigskip

(ii) From $(2_+^{1+6} \times 2^4).A_8$ we only need one character of
defect $0$ and the projectives coming from the factor group $A_8$. 
The ordinary character table of this maximal subgroup
was computed by Fischer~\cite{fischer}\@.
We refer the reader to his survey article~\cite{FiCli}, where he explains
his methods. We obtain:
\[
\begin{array}{rrrrr}
 1 &  1 &  . &  . & . \\
 14 &  1 &  1 &  . & . \\
 21 &  . &  . &  1 & . \\
 56 &  . &  1 &  . & 1 \\
 64 &  . &  . &  1 & 1 \\
& \Phi_4 & \Phi_3 & \Phi_2 & \Phi_1
\end{array}
\quad\quad\quad\quad
\begin{array}{rrrrr}
 7 &  1 & . \\
 21_1 &  . & 1 \\
 21_2 &  . & 1 \\
 28 &  1 & 1 \\
& \Phi_6 & \Phi_5
\end{array}
\]
Let $\Phi_{7}$ denote the defect 0 character of degree $15$.

We now use {\sf MOC} to calculate the $5$-modular characters of
the second Conway group $\mbox{\em Co}_2$\@.
The Brauer characters are denoted
by their degrees. The $60$ ordinary characters of {\em Co}$_2$ 
split into $26$ blocks.
The $23$ characters $$9625_1, 9625_2, 23000, 31625_1, 31625_2, 63250, 91125_1,
91125_2, 221375,$$ $$253000, 284625,
442750, 462000, 664125_1, 664125_2, 664125_3,
853875,$$ $$1288000, 1771000_1, 1771000_2, 1992375, 2004750, 2095875$$ 
are of defect $0$, and therefore each forms a block in its own. 
The characters $275$, $44275$,
$113850$, $398475$ and $467775$ belong to the second block. It is a block of
defect $1$ and the decomposition matrix can be found in \cite{BBB}. It is
$$\begin{array}{rrrrr} 275 & 1 & . & . & . \\ 44275 & . & 1 & . & . \\
113850 & 1 & . & 1 & . \\ 398475 & . & 1 & . & 1 \\ 467775 & . & . & 1 & 1 
\end{array}$$
In the third block we find the characters $4025, 12650, 177100, 398475$ and
$558900$. The decomposition matrix can be found in \cite{BBB} as well:
$$\begin{array}{rrrrr} 4025 & 1 & . & . & . \\ 12650 & . & 1 & . & . \\
177100 & 1 & . & 1 & . \\ 398475 & . & 1 & . & 1 \\ 558900 & . & . & 1 & 1 
\end{array}$$
All remaining characters are in the principal block. We have to find
$16$ irreducible Brauer characters and PIMs. We induce the projective 
characters from the two maximal subgroups and call them $\Psi_1, \ldots, 
\Psi_{59}$ and $\Phi_1, \ldots, \Phi_{104}$\@.
{\sf MOC} finds $16$ linearly independent projective
characters lying in the principal block. With Lemma \ref{projbastest} we check that
they form a basic set. They are collected in 
$${\bf PS} =
\{\Psi_{37}, \Psi_{51}, \Psi_{46}, \Psi_{39}, \Psi_{43}, \Psi_{42},
\Psi_{38}, \Psi_{34},\Psi_{49}, $$ $$ \Phi_6, \Psi_{11}, \Psi_{32}, \Psi_{31},
\Psi_{20}, \Psi_8, \Psi_4 \} .$$
Next {\sf MOC} finds a basic set of Brauer characters. It consists of
$${\bf BS} = \{ 1, 23, 253, 1771, 2024, 2277, 7084, 10395_1, 31878, 37422,
129536, $$ $$184437, 212520, 239085_1, 368874, 1291059 \} .$$
Furthermore we consider the following three projectives which do not decompose
into ${\bf PS}$ with non-negative coefficients:
{\small
$$\begin{array}{r|@{\hspace{ 6pt}}*{16}{r@{\hspace{3.5pt}}}} \hline
 & \Psi_{37} & \Psi_{51} & \Psi_{46} & \Psi_{39} & \Psi_{43} & \Psi_{42}
 & \Psi_{38} & \Psi_{34} & \Psi_{49} & \Phi_6    & \Psi_{11} & \Psi_{32} 
 & \Psi_{31} & \Psi_{20} & \Psi_8    & \Psi_4 
   \rule[- 7pt]{0pt}{ 20pt} \\ \hline 
\Phi_4 & -1 & . & . & . & . & . & 1 & . & . & . & . & 2 & . & . & . & 1 
 \rule[  0pt]{0pt}{ 13pt} \\
\Phi_5 & -1 & 1 & 1 & . & . & 1 & . & . & . & . & . & . & . & . & . & . \\ 
\Phi_7 & -1 & . & . & . & . & . & . & . & . & . & 1 & . & . & . & 1 & . 
\rule[- 7pt]{0pt}{  5pt} \\ \hline
\end{array}
$$
}
We have some more Brauer characters:
{\small
$$\begin{array}{r|@{\hspace{ 6pt}}*{16}{r@{\hspace{3.5pt}}}} \hline
 & \varphi_{ 1} & \varphi_{ 2} & \varphi_{ 3} & \varphi_{ 4} 
 & \varphi_{ 5} & \varphi_{ 6} & \varphi_{ 7} & \varphi_{ 8} 
 & \varphi_{ 9} & \varphi_{10} & \varphi_{11} & \varphi_{12} 
 & \varphi_{13} & \varphi_{14} & \varphi_{15} & \varphi_{16}
  \rule[- 7pt]{0pt}{ 20pt} \\ \hline 
   245\,916 & -1 &  1 & -1 &  1 &  . & -1 & -1 &  . &  1 &  1 &
               . &  1 &  . &  . &  . &  .  \rule[  0pt]{0pt}{ 13pt} \\
   312\,984 &  . &  1 & -1 & -1 &  . & -1 &  1 &  . & -1 &  . &
               1 &  . &  1 &  . &  . &  .  \\
   637\,560 &  . &  . &  . &  1 & -1 & -1 &  . &  . &  1 &  . &
              -1 &  2 &  . &  . &  1 &  .  \\
1\,835\,008 &  . & -1 &  . & -1 &  . &  . &  1 &  . & -1 & -1 &
               . &  . &  . &  1 &  1 &  1  \\
2\,072\,576 &  . &  . &  . & -3 &  1 &  . &  3 & -1 & -2 & -1 &
               . & -1 &  1 &  2 &  1 &  1 
\rule[- 7pt]{0pt}{  5pt} \\ \hline
\end{array} 
$$ 
}
These projectives and Brauer characters define relations (see Section
\ref{relations})\@.
With these choices of ${\bf PS}$ and ${\bf BS}$ we get the matrix~$U$
of scalar products displayed in Table~\ref{tabco2}\@.
\begin{table}[tb]
\caption{\label{tabco2} Scalar products in the 
principal block of $\mbox{\em Co}_2$}
{\small
$$\begin{array}{r|@{\hspace{ 3pt}}*{16}{r@{\hspace{2.7pt}}}} \hline
 & \Psi_{37} & \Psi_{51} & \Psi_{46} & \Psi_{39} & \Psi_{43} & \Psi_{42}
 & \Psi_{38} & \Psi_{34} & \Psi_{49} & \Phi_6    & \Psi_{11} & \Psi_{32}
 & \Psi_{31} & \Psi_{20} & \Psi_8    & \Psi_4
   \rule[- 7pt]{0pt}{ 10pt} \\ \hline
      1 & . & . & . & . & . & . & . & . & . & . & . & . & . & . & . & 1 \\
     23 & . & . & . & . & . & . & . & . & . & . & . & . & . & . & 1 & . \\
    253 & . & . & . & . & . & . & . & . & . & . & . & . & . & 1 & . & . \\
   1771 & . & . & . & . & . & . & . & . & . & . & . & . & 1 & . & . & . \\
   2024 & . & . & . & . & . & . & . & . & . & . & . & 1 & . & . & . & 1 \\
   2277 & . & . & . & . & . & . & . & . & . & . & 1 & . & . & . & 1 & . \\
   7084 & . & . & . & . & . & . & . & . & . & 1 & . & . & 1 & . & . & . \\	
 10395_1& . & . & . & . & . & . & . & . & 1 & . & . & . & . & . & . & . \\
  31878 & . & . & . & . & . & . & . & 1 & . & . & 1 & . & . & . & . & . \\
  37422 & . & . & . & . & . & . & 1 & . & . & . & . & . & . & 1 & . & . \\
 129536 & . & . & . & . & . & 1 & . & . & . & 1 & . & . & 1 & 1 & . & . \\
 184437 & . & . & . & . & 1 & . & . & . & . & 1 & . & 1 & . & . & . & 1 \\
 212520 & . & . & . & 1 & . & . & . & 1 & . & . & 2 & . & . & . & . & . \\
239085_1& . & . & 1 & . & . & . & . & . & . & . & . & . & . & . & . & . \\
 368874 & . & 1 & . & . & . & 1 & . & . & 1 & . & . & . & . & 1 & 1 & . \\
1291059 &  1 & 1 & . & . & 1 & . & 1 & 1 & . & . & 1 & . & . & . & 1 & 1
\end{array}$$
}
\end{table}
By Lemma \ref{indatom} we know that $\Psi_{37}, \Psi_{46}$ and
$\Psi_{39}$ are PIMs. We are going to show that $\Psi_{43}$ is indecomposable.
If we decompose $\Psi_{43}$ into the atoms we get 
$$[\Psi_{43}] = (0, 0, 0, 0, 0, 0, 0, 0, 0, 0, 0, 1, 0, 0, 0, 1) \cdot
[ {\bf PA} ].$$
Therefore every part $\Psi'$ of $\Psi_{43}$ must be of the form
$$[\Psi'] = (0, 0, 0, 0, 0, 0, 0, 0, 0, 0, 0, n_{12}', 0, 0, 0,
n_{16}') \cdot [ {\bf PA} ]$$
where $0 \leq n_{12}', n_{16}' \leq 1$. As we should have $\Psi' \neq 0$ and
$\Psi' \neq \Psi_{43}$ we get $n_{12}' = 0$ and $n_{16}' = 1$ or vice versa.
But then we get a negative scalar product with $2072576$. So $\Psi_{43}$ is
indecomposable. We obtain the same result for $\Psi_{42}, \Psi_{38}, \Psi_{49}$
and $\Psi_{32}$ by using the Brauer characters $637560$ and three times 
$2072576$. Now we are going to prove that $\Psi_{34}$ is indecomposable.
We have 
$$[\Psi_{34}] = (0, 0, 0, 0, 0, 0, 0, 0, 1, 0, 0, 0, 1, 0, 0, 1)
\cdot [ {\bf PA} ].$$
Therefore every part of $\Psi_{34}$ must be some
$$[\Psi'] = (0, 0, 0, 0, 0, 0, 0, 0, n_9', 0, 0, 0, n_{13}', 0, 0,
n_{16}') \cdot [ {\bf PA} ]$$
with $0 \leq n_9', n_{13}', n_{16}' \leq 1$. 
Similarly to the preceding proof we have
$1 \leq n_9' + n_{13}' + n_{16}' \leq 2$. But for all those parts there is a
Brauer character having negative scalar product with $\Psi'$ or $\Psi_{34} -
\Psi'$ as the following table shows:
$$\begin{array}{rrr|r} n_9' & n_{13}' & n_{16}' & \mbox{\rm Character} \\ \hline
0 & 0 & 1 & 1835008 \\ 0 & 1 & 0 & 312984 \\ 0 & 1 & 1 & 312984 \\ 1 & 0 & 0 &
312984 \\ 1 & 0 & 1 & 312984 \\ 1 & 1 & 0 & 1835008 \end{array}$$
Therefore, $\Psi_{34}$ is indecomposable. We get similar results for
$\Phi_6, \Psi_{11}, \Psi_{31}$ and $\Psi_{20}$ by applying the five
Brauer characters. Up to now, we know that $13$ of the $16$ projectives 
belonging to ${\bf PS}$ are indecomposable. In the next paragraph we are going
to apply our algorithms to $\Psi_{51}, \Psi_8$ and $\Psi_4$.
\end{beisp}

\subsection{Subtracting indecomposables}
\label{subtractalg}

In this paragraph we are going to present an algorithm for improving basic
sets by subtracting irreducible respectively indecomposable
characters with the help of relations.
We keep the notation introduced at the beginning of this section.
To start with, we need some technical definitions.

If one of our projective characters of ${\bf PS}$ is known to be 
a PIM, we would like to have upper bounds for the number of times
it can possibly be contained in any of the other projectives of the basic set.
Such bounds can easily be determined.
\begin{definit}
\label{maxmult}
{\rm 
Assume that $\Phi_j$ is a PIM for $j \in J$, 
a subset of $\{1, \ldots , s\}$\@.
Define non-negative integers (or infinity) $m_{ij}, 1 \leq i,j \leq s,$ by:
\begin{itemize}
\item[(i)] If $i \not\in J$, let 
$m_{ij} = \infty$.
\item[(ii)] If $i, j \in J$,
let $m_{ij} = \delta_{ij}$.
\item[(iii)] If $i \in J$ and $j \not\in J$,
then $$m_{ij} = \max \{ n \in \N_0 \mid 
\langle \varphi,
\Phi_j - n \cdot \Phi_i \rangle \geq 0 
\mbox{ for all }
\varphi \in {\bf B} \cup {\bf BS} \}.$$
\end{itemize}
Then obviously $m_{ij}$ denotes how often $\Phi_i$ can possibly
be contained in $\Phi_j.$
Therefore $m_{ij}$ is called the {\sl maximal multiplicity of $\Phi_i$ in
\index{maximal multiplicity}
$\Phi_j$}.
}
\end{definit}

Suppose that $\Phi \in {\bf PS}$ is a PIM\@. If we knew the 
corresponding irreducible Brauer character~$\varphi$, we could of course
find the multiplicity of~$\Phi$ in any other projective character~$\Psi$
by calculating $\langle \varphi, \Psi \rangle$\@. 
Of course, such a situation is exceptional. In general, there will be Brauer characters
which have positive scalar product with~$\Phi$, but none
of them is known to be irreducible. However, 
each of them contains~$\varphi$ as a part (see Definition~\ref{partdef})\@.
Let~$\vartheta$ be a Brauer character with positive 
scalar product with~$\Phi$\@.
The idea is to consider a set of parts of~$\vartheta$ containing~$\varphi$\@.
If all these special parts, which are called bits and which are
introduced in Definition~\ref{bitdefinit}
below,  have positive
scalar product with~$\Psi$, then~$\Phi$ is contained in~$\Psi$\@.
We can apply the same idea in a slightly more general situation.
\begin{definit}
{\rm
Let $\Phi = n_1 \varphi^*_1 + \ldots +
n_s \varphi^*_s$, with $n_i \geq 0$ for all $i$, 
be some generalized projective character.
Then $\Phi$ is
called {\sl multiplicity free} if $n_i \in \{ 0,1 \}$ for all $1 \leq
\index{multiplicity free}
i \leq s$.
}
\end{definit}

\begin{definit}
\label{bitdefinit}
{\rm
Let $\varphi \in$ {\bf BS}, $\varphi = n_1 \Phi^*_1 + \ldots + n_s \Phi^*_s$
with $n_1 = \langle \varphi, \Phi_1 \rangle > 0$. 
We are going to define a bit of
$\varphi$ with respect to $\Phi_1$ if either $\Phi_1$ is 
indecomposable or if $\Phi_1$ is multiplicity free.
\begin{itemize}
\item[(i)] Suppose that $\Phi_1, \ldots, \Phi_t$ are PIMs for some
$1 \leq t \leq s$\@. Then
$\varphi' = n_1' \Phi^*_1 + \ldots + n_s' \Phi^*_s$
is called a {\sl bit of $\varphi$ with respect to $\Phi_1$} if the following
\index{bit of a character}
conditions are satisfied:
\begin{itemize}
\item[(a)] $\varphi'$ is a part of $\varphi$, i.e., $0 \leq n_i' \leq n_i$ 
for all~$i$,
\item[(b)] $n_1' = 1, n_2' =\ldots=n_t' = 0, n_i' \leq m_{1i}$ for all $i > t$,
\item[(c)] $\langle \varphi', \Psi \rangle \geq 0$ and 
$\langle \varphi - \varphi', \Psi \rangle \geq 0$
for all $\Psi \in {\bf P}$.
\end{itemize}
\item[(ii)] If $\Phi_1$ is multiplicity free and $\Phi_{t+1}, \ldots, \Phi_s$
are indecomposable,
then $\varphi' = 
n_1' \Phi^*_1 + \ldots + n_s' \Phi^*_s$ is called a {\sl bit of $\varphi$
with respect to $\Phi_1$} if
\begin{itemize}
\item[(a)] $\varphi'$ is a part of $\varphi$, i.e., $0 \leq n_i' \leq n_i$ 
for all~$i$,
\item[(b)] $n_1' = 1$ and  $n_i' \leq m_{i1}$ for all $i > t$,
\item[(c)] $\langle \varphi', \Psi \rangle \geq 0$ and
$\langle \varphi - \varphi', \Psi \rangle \geq 0$
for all $\Psi \in$ {\bf P}.
\end{itemize}
\end{itemize}
}
\end{definit}

\begin{rem}
\label{gomoryhypo2}
\rm
Let $\Phi_1$ be a PIM\@.
If $|{\bf P}| =r$ and $v_{k1}, \ldots, v_{ks}$ are the coefficients of the
$k$-th character in ${\bf P}$ when decomposed into ${\bf PS}, 1 \leq k \leq r$,
then the conditions that $\varphi$ is a bit with respect to $\Phi_1$ can be
expressed by the following system of linear inequalities
$A \cdot x \leq b$ with:
$$A = \left( \begin{array}{ccc} 1 & & 0 \\ & \ddots & \\ 0 & & 1
\\ v_{1,t+1} & \ldots & v_{1,s} \\
\vdots & & \vdots \\ v_{r,t+1} & \ldots & v_{r,s} \\ - v_{1,t+1} & \ldots &
- v_{1,s} \\ \vdots & & \vdots \\ - v_{r,t+1} & \ldots &  - v_{r,s} \end{array}
\right) \mbox{ and } b = \left( \begin{array}{c} \min \{ n_{t+1}, m_{1,t+1} \}
\\ \vdots \\ \min \{ n_s, m_{1,s} \} \\ - v_{1,1} + \sum_{j=1}^s v_{1j}
\langle \varphi, \Phi_j \rangle \\ \vdots \\ - v_{r,1} + \sum_{j=1}^s v_{rj}
\langle \varphi, \Phi_j \rangle
\\ v_{1,1} \\ \vdots \\ v_{r,1} \end{array} \right)$$
and $x = (n_{t+1}', \ldots, n'_s)^{t}$.
We observe that the number of rows in $A$ and $b$ depends on $|{\bf P}|$. A
similar remark applies in case $\Phi_1$ is multiplicity free.
\end{rem}

\noindent The following two lemmas justify our definition of a bit.
\begin{lem}
\label{pimbit}
If $\Phi_1, \ldots, \Phi_t$ are PIMs and $\varphi \in {\bf BS}$ with
$\langle \varphi, \Phi_1 \rangle > 0$ then the
irreducible Brauer character $\beta$ corresponding to $\Phi_1$ is 
a bit of $\varphi$ with respect to~$\Phi_1$\@.
In particular, the set of bits of~$\varphi$ with respect to~$\Phi_1$
is non-empty.
\end{lem}
\begin{bew}
As $\langle \varphi, \Phi_1 \rangle > 0$ the irreducible Brauer character
$\beta = \sum_{i=1}^s n_i'\Phi_i^*$ corresponding 
to $\Phi_1$ is contained in $\varphi$. Thus it is a part of $\varphi$,
hence $0 \leq n_i' \leq n_i$ for all $i$. Clearly, $n_1' = 1$. All other PIMs must
have scalar product $0$ with $\beta$, so $n_i' = 0$ for $2 \leq i \leq t$.
If $i > t$ then $n_i' = \langle \beta, \Phi_i \rangle \leq m_{1i}$.
As $\beta$ is a genuine Brauer character its scalar product with any
projective is non-negative, therefore condition (c) is fulfilled.
Thus $\beta$ is a bit.
\end{bew}

\begin{lem}
\label{mulfreebit}
Suppose that $\Phi_{t+1}, \ldots, \Phi_s$ are PIMs and that
$\Phi_1 \in {\bf PS}$ is multiplicity free. 
Let $\Pi_1, \ldots, \Pi_s$ denote all the PIMs\@.
We assume that
they are numbered so that $\Pi_i = \Phi_i$ for $t+1 \leq i \leq s$.
Furthermore, let $\Phi_1 = v_1 \Pi_1 + \ldots + v_s \Pi_s$\@.
Then the following statements hold:
\begin{itemize}
\item[\rm(i)] $v_i \in \{0, 1\}$ for all $i$.
\item[\rm(ii)] Let $1 \leq k \leq s$ with $v_k = 1$. If $\varphi \in
{\bf BS}$ such that $\langle \varphi, \Pi_k
\rangle > 0$, then $\langle \varphi, \Phi_1 \rangle > 0$ and
the irreducible Brauer character $\beta$ corresponding to $\Pi_k$
is a bit of $\varphi$ with respect to $\Phi_1$.
In particular, the set of bits of~$\varphi$ with respect to~$\Phi_1$
is non-empty.
\end{itemize}
\end{lem}
\begin{bew}
\begin{itemize}
\item[(i)] As $\Phi_1$ is multiplicity free, its decomposition into ${\bf PA}$
leads to coefficients in $\{ 0, 1 \}$. Now every PIM is a non-negative linear
combination of ${\bf PA}$. The result follows immediately.
\item[(ii)] Let $n_1', \ldots, n_s'$ be the coefficients of $\beta$ when
decomposed into ${\bf BA}$. Conditions (a) and (c) of
Definition~\ref{bitdefinit} are fulfilled as
$\beta$ is contained in $\varphi$ and is a genuine Brauer character. We have
$$\langle \beta, \Phi_1 \rangle = \sum_{i=1}^s v_i \langle \beta, \Pi_i \rangle
= v_k = 1,$$
that is $n_1' = 1$.
If $i > t$ then
$$n'_i = \langle \beta, \Phi_i \rangle = \langle \beta, \Pi_i \rangle =
\delta_{ik}.$$
If $t < i \neq k$ then $\delta_{ik} = 0 \leq m_{i1}$. If $i = k$ then $n_k' = 1
=
m_{k1}$ as $\Phi_1$ is multiplicity free and
contains $\Phi_k = \Pi_k$ exactly once. This completes
the proof.
\end{itemize}
\end{bew}

\begin{definit}
{\rm 
Let $\Phi \in {\bf PS}$ be indecomposable or multiplicity free. 
For $\varphi \in {\bf BS}$ with $\langle \varphi, \Phi \rangle > 0$
and a projective character $\Sigma \in {\bf P}$ let
$$m(\Phi, \Sigma, \varphi) = \min \{ \langle \varphi', \Sigma
\rangle \mid \varphi'
\mbox{ is a bit of } \varphi \mbox{ corresponding to } \Phi \}.$$
}
\end{definit}

\noindent Let $\Phi \in {\bf PS}$ be indecomposable or multiplicity free.
There is at least one
$\varphi \in {\bf BS}$ such that $\langle \varphi, \Phi \rangle > 0$\@.
By the preceding lemmas, there is always a bit of~$\varphi$ 
with respect to~$\Phi$ and so $m(\Phi, \Sigma, \varphi)$ is well defined.
By our definition of bits, $m(\Phi, \Sigma, \varphi)$ is a 
non-negative integer.
The set 
$$\{ m(\Phi, \Sigma, \varphi) \mid \varphi \in {\bf BS}, \langle 
\varphi, \Phi \rangle > 0 \}$$
is non-empty. It is, of course, finite, since~$\varphi$ has
only finitely many parts.

\begin{satz}
Let $\Phi\in {\bf PS}$ and $\Sigma \in {\bf P}$\@.
If $\Phi$ is indecomposable, let
$$z = \max \{ m(\Phi, \Sigma, \varphi) \mid
\varphi \in {\bf BS}, \langle \varphi, \Phi \rangle > 0 \}.$$
If $\Phi$ is multiplicity free, let
$$z = \min \{ m(\Phi, \Sigma, \varphi) \mid
\varphi \in {\bf BS}, \langle \varphi, \Phi \rangle > 0 \}.$$
Then $\Sigma - z \cdot \Phi$ is a genuine projective character.
\end{satz}
\begin{bew}
As remarked above, $z$ is well-defined.
Suppose first, that $\Phi$ is indecomposable.
Let $\{ \beta_1, \ldots, \beta_s \}$ be the irreducible Brauer characters.
Without loss of generality let $\beta_1$ be the one 
corresponding to $\Phi$.
For all $1 \leq i \leq s$ with
$\langle \varphi_i, \Phi \rangle > 0$ let
$z_i = m(\Phi, \Sigma, \varphi_i)$. As
$$\beta_1 \in \{ \varphi' \mid \varphi' \mbox{ bit of } 
\varphi_i \mbox{ with respect to }
\Phi \}$$
we have $\langle \beta_1, \Sigma \rangle \geq z_i$. Thus~$\Phi$ is
contained at least~$z_i$ times in~$\Sigma$\@. 
As $z = \max z_i$, also $\Sigma - z \cdot \Phi$
is projective.

Suppose now, that~$\Phi$ is multiplicity free.
Let $\Phi = v_1 \Pi_1 + \ldots + v_s \Pi_s$ be the decomposition of $\Phi$
into the PIMs. Then $v_i \in \{ 0, 1 \}$ by Lemma~\ref{mulfreebit}(a)\@.
Now let~$k$ be such that $\Pi_k$ occurs in~$\Phi$, i.e., $v_k = 1$\@.
Let~$\beta_k$ be the irreducible Brauer character corresponding
to~$\Pi_k$\@. We shall show that $\langle \beta_k , \Sigma \rangle
\geq z$\@. This will imply the result since~$\Phi$ is multiplicity free.
 
Since ${\bf BS}$ is a basic set, there is some $\varphi \in {\bf BS}$
with $\langle \varphi, \Pi_k \rangle > 0$\@. Then
$\beta_k$ is a bit of~$\varphi$
corresponding to~$\Phi$ by Lemma~\ref{mulfreebit}(b)\@. 
Hence $m(\Phi,\Sigma,\varphi) \leq
\langle \beta_k, \Sigma \rangle$\@. Since obviously
$\langle \varphi, \Phi \rangle > 0$, we have $z \leq m(\Phi,\Sigma,\varphi)$\@.
\end{bew}

\noindent  In \S~\ref{gomoryint} we shall show how to solve 
the problem of calculating $z$\@.
Changing projectives and Brauer characters give analogous results for
Brauer characters.

\begin{beisp}
\rm
Let us now return to our example $\mbox{\em Co}_2$ mod $5$. As mentioned above
we still have to find three 
more PIMs. First we want to show that $\Psi_{37}$ is
contained in $\Psi_{51}$. For this purpose we have to find the irreducible
Brauer character corresponding to $\Psi_{37}$. But this is a bit $\varphi'$
of $\varphi = 1291059$.
Its decomposition into the atoms is 
$$[\varphi] = (1, n_2', 0, 0, n_5', 0, n_7', n_8', 0, 0,
n_{11}', 0, 0, 0, n_{15}', n_{16}') \cdot [{\bf BA}]$$
where for all $i$ we have $0 \leq n_i' \leq
1$. As $\Psi_{43}, \Psi_{38}, \Psi_{34}$ and $\Psi_{11}$ are known to
be indecomposable we have $n_5' = n_7' = n_8' = n_{11}' = 0$. So it remains
to solve the following system of linear inequalities (remember that $\langle \varphi',
\Psi_{51} \rangle = n_2'$):

$$\mbox{Minimize \ } n_2' \mbox{ subject to}$$
$$\left( \begin{array}{rrr} 1 & 0 & 0 \\ 0 & 1 & 0 \\ 0 & 0 & 1 \\ 0 & 0 & 1 \\
1 & 0 & 0 \\ 0 & 1 & 0 \\ 0 & 0 & -1 \\ -1 & 0 & 0 \\ 0 & -1 & 0 \end{array}
\right) \cdot \left( \begin{array}{c} n_2' \\ n_{15}' \\ n_{16}' \end{array}
\right) \leq \left( \begin{array}{r} 1 \\ 1 \\ 1 \\ 2 \\ 1 \\ 1 \\ -1 \\ -1 \\
-1 \end{array} \right)$$
where all $n_i' \in \N_0$. As this minimum equals $1$ we find that
$\Psi_{51} - \Psi_{37}$ is a genuine projective character. Since
this projective is an atom
it is indecomposable. Similarly, we show that $\Psi_8 - \Psi_{37}$ and
$\Psi_4 - \Psi_{37}$ are projective. The two Brauer characters of degree
$245916$ and $2072576$ prove that the resulting projectives are indecomposable.
So we have found all PIMs. The resulting decomposition matrix as well
as the table of irreducible Brauer characters is printed in
the appendix. 
\end{beisp}

\subsection{Triangular matrix of scalar products}
\label{triangularalg}

In \S~\ref{subtractalg} we presented an algorithm for subtracting irreducible
or multiplicity free characters. In some cases, it is possible to improve
characters without knowing whether they are irreducible or not.
For this algorithm we assume that the matrix~$U$ of scalar products is
lower uni-triangular, i.e., $a_{ij} = \langle \varphi_i, \Phi_j \rangle = 0$ 
if $j > i$. Let $\{\Pi_1, \ldots,
\Pi_s\}$ denote the PIMs 
and $\{\beta_1, \ldots, \beta_s\}$ the irreducible Brauer characters.
We assume that the PIMs and ${\bf PS}$ are ordered consistently, i.e., 
that $\Phi_j$ consists of $\Pi_j$ and a sum of $\Pi_r$'s with $r > j$. 
Here we are going to describe a method which allows to 
subtract $\Pi_i$ from some
$\Phi_{j_0}$ without knowing what $\Pi_i$ looks like. 
Suppose we know for some reason that $\Pi_i$ is contained in $\Phi_{j_0}$
(see Lemma~\ref{lemmaB} below)\@. The idea now is to subtract $\Phi_i$
from $\Phi_{j_0}$\@. The possible errors due to the fact that $\Phi_i$
might not be indecomposable can easily be corrected by adding some
projectives $\Phi_l$ for $l > i$ (see Theorem~\ref{dreiecksstrip})\@.
This method produces
new characters which are smaller than the older ones in the following sense.
If we identify a projective~$\Phi$ with the column 
$(\langle \varphi_1, \Phi \rangle, \ldots , \langle \varphi_s, \Phi \rangle)^t$,
we can order the set of projectives lexicographically. The new characters
obtained are then smaller in this lexicographic ordering.

\begin{lem}
\label{lemmaA}
Let $\Phi_j = \Pi_j + \sum_{r=j+1}^s b_{rj} \Pi_r$ be the decomposition of
$\Phi_j$ into the PIMs. 
Then $b_{lj} \leq a_{lj} = \langle \varphi_l, \Phi_j \rangle$ 
for all $l = j+1, \ldots, s$\@. 
\end{lem}
\begin{bew}
Let $\varphi_l = \beta_l + \mu$ for a certain Brauer character $\mu$. 
Then
$a_{lj} = \langle \varphi_l, \Phi_j \rangle = 
\langle \beta_l, \Phi_j \rangle +
\langle \mu, \Phi_j \rangle  = 
b_{lj} + \langle \mu, \Phi_j \rangle \geq b_{lj}.$
\end{bew}

\begin{lem}
\label{lemmaB}
Let $\Phi$ be a projective character, $\Phi = \sum_{j=1}^s v_j \Phi_j$.
Suppose that $v_i < 0$ for some $i$. If there
exists $j_0 < i$ with $v_{j_0} > 0$ and 
$$| v_i | > \sum_{j=1, j \neq j_0, v_j > 0}^{i-1} v_j a_{ij}$$
then $\Pi_i$ is contained in $\Phi_{j_0}$ for a least $z$ times, where
\begin{equation}
\label{zdefinit}
z = \lceil \frac{1}{v_{j_0}} (| v_i | - \sum_{j=1, j \neq j_0, v_j > 0
} ^{i-1} v_j
a_{ij}) \rceil.
\end{equation}
(Here $\lceil x \rceil$ for $x \in \R$ denotes the integer $y$ with $y - 1
< x \leq y$.)
\end{lem}
\begin{bew}
Let $\Lambda = \sum_{j=1,j \neq i, v_j \leq 0}^s v_j \Phi_j$. Then 
$$\Phi - \Lambda + | v_i | \Phi_i = \sum_{j=1,v_j > 0}^s v_j \Phi_j$$
(remember that $v_i < 0$). Therefore $| v_i | \Phi_i$ is contained in the right
hand of this equation, hence so is $| v_i | \Pi_i$\@. By Lemma~\ref{lemmaA} $\Pi_i$ is
contained in $\Phi_j$ only for $j \leq i$ and in this case for at most $a_{ij}$
times. So $\Pi_i$ is contained in $v_{j_0} \Phi_{j_0}$ for at least 
$$| v_i | - \sum_{j=1,j \neq j_0, v_j > 0}^{i-1} v_j a_{ij}$$
times and therefore $\Pi_i$ is contained in $\Phi_{j_0}$ at least $z$ times.
\end{bew}

\begin{satz}
\label{dreiecksstrip}
Let the notation be as in {\rm Lemma~\ref{lemmaA}}\@. 
Suppose that~$\Pi_i$ is contained~$z$ times 
in~$\Phi_{j_0}$ for some $i > j_0$\@.
Then 
\begin{equation}
\label{PhiLexSmall}
\tilde{\Phi}_{j_0} = \Phi_{j_0} - z \Phi_i + z \sum_{l=i+1}^s a_{li} \Phi_l
\end{equation}
is a genuine projective character.
\end{satz}
\begin{bew}
We have
$$\Phi_{j_0} - z \Phi_i + z \sum_{l=i+1}^s a_{li} \Phi_l$$
$$= \Phi_{j_0} - (z \Pi_i + z \sum_{l=i+1}^s b_{li} \Pi_l)
+ z \sum_{l=i+1}^s a_{li} \Phi_l$$
$$= {(\Phi_{j_0} - z \Pi_i)} + z \sum_{l=i+1}^s {(a_{li} - b_{li})} \Phi_l
+ z \sum_{l=i+1}^s b_{li} {(\Phi_l - \Pi_l)},$$
which by Lemma \ref{lemmaA} 
is a non-negative linear combination of genuine projectives and therefore
is a genuine projective character.
\end{bew}

\noindent The following algorithm derived from Lemma \ref{lemmaB} and
Theorem~\ref{dreiecksstrip} is implemented in the {\sf MOC}-system.
Suppose that ${\bf P}$ is a set of projective characters.
For $\Phi \in {\bf P}$ write $\Phi = \sum_{j=1}^s v_j(\Phi) \Phi_j$\@.

\begin{algorithm}
\rm
\mbox{\ \ }

\begin{itemize}
\item[] For all $i$ in $1, \ldots, s$ do:
\begin{itemize}
\item[] For all $j_0$ in $1, \ldots, i-1$ do:
\begin{itemize}
\item[] $z := 0$;
\item[] For all $\Phi \in$ {\bf P} do:
\begin{itemize}
\item[] If $v_{j_0}(\Phi) \leq 0$ or $v_i(\Phi) \geq 0$
\item[] \hspace*{2em} $z_\Phi := 0$
\item[] else
\item[] \hspace*{2em} Calculate $z_{\Phi}$ according to (\ref{zdefinit});
\item[] fi
\item[] $z := \max \{ z, z_{\Phi} \};$
\end{itemize}
\item[] od;
\item[] If $z > 0$ calculate $\tilde{\Phi}_{j_0}$ 
according to (\ref{PhiLexSmall});
\end{itemize}
\item[] od;
\end{itemize}
\item[] od.
\end{itemize}
\end{algorithm}

\noindent If, in the above algorithm, $z > 0$ for some $j_0$, $i > j_0$,
then $\tilde{\Phi}_{j_0}$ is lexicographically smaller than ${\Phi}_{j_0}$
in the sense described at the beginning of this paragraph.

\subsection{Proving decomposability}
\label{decomposablealg}

Sometimes one can show that a projective or Brauer character has to
be divided into a sum of two or more characters. Let us assume that there is
a projective $\Phi \in {\bf P}$ which has at least one negative coefficient
in its decomposition into the basic set. Let $\Phi = v_1 \Phi_1 + \ldots +
v_s \Phi_s$ be this decomposition and assume that $v_i < 0$ for some $i$.
Hence, $\Phi_i$ is contained in $v_1 \Phi_1 + \ldots + v_{i-1} \Phi_{i-1} +
v_{i+1} \Phi_{i+1} + \ldots + v_s \Phi_s$. 
If $\Phi_i$ is indecomposable it has
to be contained in at least one $v_j \Phi_j$ for $j \neq i$. This can be
checked by using an algorithm similar to the one used for calculating the
$m_{ij}$ (see Definition~\ref{maxmult}). If there is no such
$j$ then $\Phi_i$ has to be decomposable. 

The easiest and most important case
is that $\Phi_i$ is the sum of two indecomposable projective characters. 
To check this we find the two lexicographically 
smallest solutions $(n_1', \ldots , n_s')$ and $(n_1'', \ldots ,n_s'')$
fulfilling the system of inequalities of Remark~\ref{gomoryhypo1}\@.
We put $\Psi_1 = \sum n_j' \varphi_j^*$ and 
$\Psi_2 = \sum n_j'' \varphi_j^*$\@. If $\Psi_1 + \Psi_2 = \Phi_i$, we
are done. In two special cases it is now easy to get a better basic set:
\begin{itemize}
\item[(i)] Assume that $\Psi_1 \in {\bf PS}$. Then $\{\Phi_1,\ldots,\Phi_{i-1},
\Psi_2, \Phi_{i+1}, \ldots, \Phi_s\}$ is a better basic set.
\item[(ii)] Assume that the matrix of scalar products is triangular. Let
$\Psi_1 = \sum_{k=1}^s w_k \varphi_k^*$ and $\Psi_2 = \sum_{k=1}^s x_k \varphi_k^*$
be the decompositions into the projective atoms. Then $w_k = x_k = 0$ for all
$k > i$. Without loss of generality
let $w_i = 1$. As $\Phi_i = \Psi_1 + \Psi_2$ and $U$ is
unitriangular we get $x_i = 0$. In ${\bf PS}$ 
we substitute $\Phi_i$ by $\Psi_1$.
Obviously, $\Psi_1$ is smaller than $\Phi_i$. Let $l = \max \{ k \mid x_k
> 0 \}$. Then $1 \leq l \leq i-1$. Now we substitute $\Phi_l$ by $\Psi_2$ and
we get the basic set required.
\end{itemize}

\subsection{Solving the second fundamental problem}

As we have seen in the last three paragraphs all algorithms depend on the
size of ${\bf P}$ and ${\bf B}$. So it becomes clear that these two sets
have to be kept small. First of all we remember 
from Section~\ref{fundamentalproblemII} that we only need to
consider essential characters. This is made precise below. 

\begin{rem}
\label{throwalg}
\rm
We assume that there is a set $E$ of projectives satisfying the following two
conditions:
\begin{itemize}
\item[(i)] ${\bf PS} \subseteq E \subseteq {\bf P}$,
\item[(ii)] For all $\Phi \in {\bf P}$ there exist some $v_1, \ldots, v_s, w_1,
\ldots, w_r \in \N_0$ such that
\begin{equation}
\label{throw}
\Phi = \sum_{i=1}^s v_i \Phi_i + \sum_{j=1}^r w_j \Sigma_j
\end{equation}
where 
$E \cap ( {\bf P} \setminus
{\bf PS}) = \{ \Sigma_1, \ldots, \Sigma_r \}$.
\end{itemize}
Then $E$ contains ``the same information'' as ${\bf P}$ in the following sense.
Suppose there is a
Brauer character $\varphi$ and a bit $\varphi'$ of $\varphi$
and a $\Sigma \in {\bf P} \setminus E$ with $\langle \varphi', \Sigma \rangle
< 0$ or $\langle \varphi
-\varphi', \Sigma \rangle < 0$\@. Then we have (without loss of generality 
$\langle \varphi', \Sigma
\rangle < 0)$ :
$$ \Sigma = \sum_{i=1}^s v_i \Phi_i + \sum_{j=1}^r w_j \Sigma_j$$
and therefore
$$0 > \langle \varphi', \Sigma \rangle = \sum_{i=1}^s v_i \langle \varphi',
\Phi_i \rangle + \sum_{j=1}^r
w_j \langle \varphi', \Sigma_j \rangle.$$
This implies $\langle \varphi', \Sigma_j \rangle < 0$ for some $j$.
So we do not need $\Sigma$.
\end{rem}

So again we have to solve a system of linear inequalities over $\N_0$ which we 
get by
subtracting $\sum_{i=1}^s v_i \Phi_i$ on the right hand side 
of equation (\ref{throw})\@.
The problem of integer linear programming is 
${\cal NP}$-complete. For a proof see \cite{burka}.
But looking at the proof of Remark \ref{throwalg}
we find that the conditions $v_1, \ldots, v_s, w_1, \ldots, w_r \geq 0$
are sufficient. So we are able to state the problem again, now over the field
of real numbers. In \cite{Karmarkar} one can find an polynomial algorithm for
solving this problem. 
Nevertheless, we use an ordinary simplex algorithm for the 
following reasons:
The solution can be constructed explicitly. So, by checking
whether it satisfies the inequalities, errors due to truncation are avoided.
Finally, it is quite fast. Our task is now:

Find a set $E$ of projectives fulfilling:
\begin{itemize}
\item[(i)] ${\bf PS} \subseteq E \subseteq {\bf P}$,
\item[(ii)] For all $\Phi \in {\bf P}$ there exists some $v_1, \ldots, v_s, w_1,
\ldots, w_r \in \R_{\geq 0}$ with
$$\Phi = \sum_{i=1}^s v_i \Phi_i + \sum_{j=1}^r w_j \Sigma_j$$
where $E \cap ({\bf P} \setminus {\bf PS}) = \{ \Sigma_1, \ldots, \Sigma_r \}$.
\end{itemize}

\noindent Here is our algorithm for solving this problem:

\begin{algorithm}
\label{algo5221}
\rm
\begin{description}
\item{{\bf Step 1.}} Let $E := {\bf PS}, F := {\bf P}$.
\item{{\bf Step 2.}} For all $\Phi \in F \setminus E$ do:
Try to solve
$$ \Phi = \sum_{i=1}^s v_i \Phi_i + \sum_{j=1}^r w_j \Sigma_j \mbox{ for } v_i,
w_j \geq 0,$$
where $E \setminus {\bf PS} = \{\Sigma_1, \ldots , \Sigma_r\}$\@. 
For this purpose we decompose $\Phi$ and $\Sigma_j$ into $\{ \Phi_1, \ldots, \Phi_s \}$
and denote by $b_1, \ldots, b_s$ and $a_{j1}, \ldots, a_{js}$ the resulting
coefficients. Then we have to find a solution $w = (w_1, \ldots ,w_r)$, 
$v = (v_1, \ldots, v_s)$ and $v_i, w_j
\in \R_{\geq 0}$
for the system $C \cdot x = b$ where $b = (b_1, \ldots, b_s)^t$, $A =
(a_{ij})_{ij}$, $C = [A,E_s]$ and $x = [w^t, v^t]$.
We multiply each equation by $-1$ or $+1$ so that the right hand side of any
equation is non-negative.
Let $\bar{C}$ and $\bar{b}$ denote the resulting coefficients.
Now we introduce some artificial variables $y_1, \ldots, y_s$ which have
to be greater or equal to zero and augment the system as follows:
$$\begin{array}{cccccccccccc}
\bar{a}_{11} x_1 & + & \ldots & + & \bar{a}_{1,r+s} x_{r+s} & + & y_1 & & & & 
 = & \bar{b}_1 \\
\vdots & & & & \vdots & & & \ddots & & & & \vdots \\
\bar{a}_{s1} x_1 & + & \ldots & + & \bar{a}_{s,r+s} x_{r+s} & & & & + & y_s & 
 = & \bar{b}_s \end{array}$$
Then the system is feasible by setting $x_i = 0$ and $y_i = \bar{b}_i$.
Now we minimize the function
$$f = y_1 + \ldots + y_s.$$
This function counts how many artificial variables we do need to make the
system solvable. If $f$ has a minimum equal to $0$, the original system is
solvable otherwise it is not.
If $\Phi$ is a non-negative linear combination of elements in $E$ throw it away,
i.e. put $F := F \setminus \{ \Phi \}$.
\item{{\bf Step 3.}} If $F \neq E$ take a new vector $\Psi \in F \setminus E$ 
and put $E := E \cup \{ \Psi \}$.
For example, choose a vector $\Psi \in F \setminus E$ with minimal sum of
coefficients. Go to 2.
\end{description}
\end{algorithm}

\begin{beisp}
\rm
Let us go back to our example {\em Co}$_2$ mod $5$. 
One useful method for generating
projectives is tensoring characters of defect $0$ with Brauer characters.
Here we tensor $91125_1$ (it is of defect $0$) with $275$ and restrict the
resulting projective to the principal block. We call this projective $\Phi$.
Its decomposition into the first ${\bf PS}$ we took for {\em Co}$_2$ mod $5$ is
$$\Phi = - \Psi_{37} + \Psi_{51}.$$
As (using the same notation as in Example \ref{c2mod5}) we have
$$\Phi_5 = \Phi + \Psi_{46} + \Psi_{42},$$
we could have thrown away $\Phi_5$. Of course, the final result would
have been the same.
\end{beisp}

\subsection{An all-integer integer simplex algorithm}
\label{gomoryint}

As seen in the earlier paragraphs we are given the following problem in linear
programming: Let $A \in
\Z^{m \times n}, b \in \Z^m$ and $c \in \Z^n$. Then
we have to find $\min \{ c^t x \}$ subject to $x \in \N_0^n$ and
$A x \leq b$ where we take $\leq$ in each component. So this problem
is similar to that which is solved by
a simplex algorithm. The main difference is the fact
that we have to find integral solutions~$x$. 
So we are interested in an algorithm only
using integers in order to avoid errors due to truncation. This problem
was solved by R.E.Gomory in 1963 \cite{gomory}. We use this Gomory algorithm
to calculate the minimum taken over all bits and to check whether a character
is irreducible (in this case we minimize the $0$-function). The Gomory algorithm
is based on the dual simplex algorithm. Here we are not going to
explain the whole theory of duality and linear programming. For details see 
\cite{colwet}. Our implementation is based on the book of Burkard \cite{burka}.
Similarly to Algorithm~\ref{algo5221}, Step~$2$, we first augment the
system of inequations by introducing some artificial variables
in order to get a system of equations.
We are going to explain only how we get integral solutions 
for the problem.
We put
$$a_{00} = 0,$$
$$a_{0j} = c_j \mbox{ for } j = 1, \ldots, n,$$
$$a_{i0} = b_i \mbox{ for } i = 1, \ldots, m$$
and write these numbers in the following table:
$$\begin{array}{c|ccccc} & & x_1 & x_2 & \ldots & x_n \\ \hline & a_{00} &
a_{01} & a_{02} & \ldots & a_{0n} \\ x_{n+1} & a_{10} & a_{11} & a_{12} & \ldots
& a_{1n} \\ \vdots & \vdots & \vdots & \vdots & & \vdots \\ x_{n+m} & a_{m0}
& a_{m1} & a_{m2} & \ldots & a_{mn} \end{array}$$
This a shortened way of writing 
$$z = \min ! = \sum_{j=1}^n a_{0j} x_j$$
subject to
$$\sum_{j=1}^n a_{ij} x_j + x_{n+i} = a_{i0} \quad \mbox{ for } \quad i = 1, \ldots, m.$$
As mentioned above, this algorithm is based on the dual simplex algorithm.
We assume now that there is a solution for the
dual problem such that the first non-zero element in each column $a_j$ for
$1 \leq j \leq n$ is positive. In the two algorithms presented in
\S\S~\ref{pimtestalg} and \ref{subtractalg}
these conditions are satisfied as remarked in \ref{gomoryhypo1} and
\ref{gomoryhypo2}\@.
Suppose now that the last solution the algorithm has found while optimizing
has a non-integral coefficient. Then we
have to exclude this point from the simplex. 
So we assure that all vertices
of the simplex have integral coefficients. The optimal
solution is found if all $b_i$ for $i = 1,\ldots,m$ are non-negative.
Otherwise one of the $b_i$ is negative. From now on we fix this $i$.
We are going to construct a new restriction out of the $i$-th row. To make the
algorithm becoming all-integer the pivot element 
will become $-1$. We will repeat
this until we get the optimal solution or we get a contradiction so that there
is no solution. Now we are going to make the new restriction. Let 
$a_{ij}$ as defined above. If $\lambda > 0$ and $a$ is a real number then $a$
can be written as
$$a = \lambda [\frac{a}{\lambda} ] + r_a \mbox{ with } 0 \leq r_a < \lambda$$
where $[ \cdot ]$ denotes the Gaussian function. Now let $a_{i0} < 0$ for some
$1 \leq i \leq m$. The corresponding restriction is
$$a_{i0} = a_{i1} x_1 + \ldots + a_{in} x_n + 1 \cdot x_{n+i}.$$
Therefore we get
$$\lambda \cdot [ \frac{a_{i0}}{\lambda} ] + r_0 = ( \lambda \cdot [
\frac{a_{i1}}{\lambda} ] + r_1) x_1 + \ldots + ( \lambda \cdot [
\frac{a_{in}}{\lambda} ] + r_n) x_n + ( \lambda \cdot [ \frac{1}{\lambda}] +
r_{n+i}) x_{n+i}$$
and then
$$\sum_{j=1}^n r_j x_j + r_{n+i}x_{n+i} = r_0 + \lambda \{ [ \frac{a_{i0}}
{\lambda} ] - \sum_{j=1}^n [ \frac{a_{ij}}{\lambda} ] x_j - [ \frac{1}{\lambda}]
x_{n+i} \}.$$
Since $r_j \geq 0$ and $x_j \geq 0$ for $j=1,\ldots,n,n+i$ the left hand of
this equation is non-negative. Also, the term in brackets is integral. Let
$$\bar{x} = [ \frac{a_{i0}}{\lambda} ] - \sum_{j=1}^n [ \frac{a_{ij}}{\lambda}]
x_j - [ \frac{1}{\lambda} ] x_{n+i}.$$
If $\bar{x} < 0$, then $\bar{x} \leq -1$, hence $\lambda \bar{x} \leq - \lambda
$. As $r_0 < \lambda$ we would have $r_0 + \lambda \bar{x}$ being negative, a
contradiction. In this algorithm we are going to choose $\lambda > 1$. Then
$[ \frac{1}{\lambda} ] = 0$ and we get the restriction
$$\bar{x} = [ \frac{a_{i0}}{\lambda} ] - \sum_{j=1}^n [ \frac{a_{ij}}{\lambda}
] x_j \geq 0.$$
From this construction we get that every integral point of the given problem
satisfies this restriction. We are going to choose a $\lambda$ such that the
pivot element
will become $-1$. Hence the system will remain integral after pivoting.
In the dual simplex algorithm the column $s$ containing the pivot element is
determined by
$$\frac{a_{0s}}{[ \frac{a_{is}}{\lambda}]} = \max \{ \frac{a_{0j}}{[ \frac 
{a_{ij}}{\lambda} ]} \mid a_{ij} < 0 , j = 1,\ldots,n \}.$$
If there is no such $a_{ij} < 0$ then the system is not solvable. As we are
going to choose $\lambda$ such that $[ \frac{a_{is}}{\lambda} ] = - 1$, we get
$$- a_{0s} = \max \{ \frac{a_{0j}}{[ \frac{a_{ij}}{\lambda}]} \mid a_{ij} < 0,
j=1,\ldots,n \}.$$
But $[ \frac{a_{ij}}{\lambda}]$ is an integer, less or equal to $-1$. If we put
$$\mu_j = [ \frac{a_{ij}}{\lambda} ]$$
we have
$$- a_{0s} \geq \frac{a_{0j}}{\mu_j} \geq - a_{0j},$$
hence $s$ is determined by
$$- a_{0s} = \max \{ - a_{0j} \mid a_{ij} < 0, j = 1,\ldots,n \}.$$
This is equivalent to
$$a_{0s} = \min \{ a_{0j} \mid a_{ij} < 0,  j=1,\ldots,n \}.$$
If $s$ is not determined uniquely consider the next rows of the table to make
a decision. We determine for $k=1,2 \ldots$
$$a_{ks} = \min \{ a_{kj} \mid a_{ij} < 0 , a_{k-1,j} \mbox{ minimal for }
a_{ij} < 0 \}$$
until the minimum is unique. So we achieve that for all $j$ with $a_{ij} < 0$
the first element of the vector $a_s - a_j$ is negative. We are now going to
examine how to choose $\lambda$ best. After pivoting we get a new value
$$\bar{a}_{00} = a_{00} + [ \frac{a_{i0}}{\lambda} ] a_{0s}$$
for the function to be minimized. In order to decrease the value as much as
possible we have to choose $\lambda$ as minimal as possible. But we have to
make sure that the first non-zero element $\bar{a}_{kj}$ of the columns
$j= 1,\ldots,n$ with $j \neq s$ is
$$\bar{a}_{kj} = a_{kj} + [ \frac{a_{ij}}{\lambda} ] a_{ks} > 0.$$
Let $i_0$ be the minimum of $i(j)$ and $i(s)$ where $i(j)$ is first index
for which $a_{i(j),j} \neq 0$ (and therefore $a_{i(j),j} > 0$). Similarly, let
$i(s)$ be the first index of rows with $a_{i(s),s} > 0$. Then $\bar{a}_{i_0j}$
becomes the first non-zero element of column $j$. If $a_{ij} \geq 0$ then
we have $\bar{a}_{i_0j} > 0$ for all possible $\lambda$ as we either have
$a_{i_0j} = 0$ and $a_{i_0s} > 0$ or $a_{i_0j} > 0$ and $a_{i_0s} \geq 0$.
But if we have $a_{ij} < 0$ and $a_{i_0j} > 0$ then we get a new restriction
for $\lambda$ resulting from the equation determing $\bar{a}_{kj}$. 
As $0 \leq a_{is} \leq a_{i0}$ we get
$$ - 1 \geq \mu_j > - \frac{a_{i_0j}}{a_{i_0s}}.$$
As $\lambda$ has to be chosen as small as possible and $a_{ij} < 0$ we have
to choose $\mu_j$ as small as possible, too. Therefore we take $\mu_j$ to be
the smallest integer greater than $- \frac{a_{i_0j}}{a_{i_0s}}$. Now 
for $j \in J = \{ j \mid a_{ij} < 0 \}$ we put:
$$\lambda_j = \left\{ \begin{array}{l} \frac{a_{ij}}{\mu_j} \mbox{ if } 
a_{i_0s} > 0 \\ 0 \mbox{ otherwise } \end{array} \right..$$
Then 
$$\lambda = \max \{ \lambda_j \mid j \in J \}$$
satisfies all conditions, i.e. $[ \frac{a_{is}}{\lambda} ] = -1$, the value
of the function being optimized decreases as much as possible and after 
pivoting in all columns the first non-zero element is positive.

\begin{algorithm}
\rm
Gomory's all-integer integer algorithm (see \cite{burka}) 
for solving the problem
$$\sum_{j=1}^n a_{0j} x_j = \min !$$
subject to $Ax \leq a_0$ for $x_j \in \N_0$ and $1 \leq j \leq n$.
We assume that the problem is 
dual feasible. Then:
\begin{enumerate}
\item If for all $i = 1,\ldots, m$
$$a_{i0} \in \N_0$$
then the optimal solution is found. Stop. Otherwise go to Step~2.
\item Let
$$\tilde{r} := \min \{ i \mid a_{i0} \notin \N_0, i = 1,\ldots, m \}.$$
\item If $a_{\tilde{r}j} \geq 0$ for all $j = 1,\ldots, n$ then the problem
is not solvable. Stop. Otherwise let
$$J := \{ j \mid j > 0 \mbox{\rm\ and } a_{\tilde{r}j} < 0 \}.$$
\item Determine $s$ to be the index of the pivoting column by
$$a_{0s} := \min \{ a_{0j} \mid j \in J \}.$$
If $s$ is not determined uniquely repeat for $i = 1,2,\ldots$
$$a_{is} := \min \{ a_{ij} \mid j \in J\mbox{ and } a_{i-1,j}\mbox{ minimal }\}$$
until the minimum is determined uniquely.
\item For $j \in J$ let $i(j)$ be the first index of rows such that $a_{i(j),j}
\neq 0$. Determine $\mu_j$ to be the smallest integer for which
$$a_{i(j),j} + \mu_j a_{i(j),s} > 0.$$
If $a_{i(j),s} = 0$ put $\mu_j := - \infty$.
\item For $j \in J$ let
$$\lambda_j := \left\{ \begin{array}{l} 0 \mbox{ if } \mu_j = - \infty \\
\frac{a_{\tilde{r}j}}{\mu_j}, \mbox{ otherwise } \end{array} \right.$$
and let
$$\lambda := \max \{ \lambda_j \mid j \in J \}.$$
If $\lambda = 1$, replace $\lambda$ by $1 + \varepsilon$ where $\varepsilon$
is some arbitrary little number.
\item For $j = 0,1,\ldots, n$ put
$$a_{m+1,j} := \lceil \frac{a_{\tilde{r}j}}{\lambda} \rceil$$
and add $a_{m+1,0}, a_{m+1,1}, \ldots, a_{m+1,n}$ as $(m+1)$st row to the table.
Let $r := m := m+1$.
\item Put
$$\bar{a}_{rj} := - a_{rj} \mbox{ for } j = 0,1, \ldots, n; j \neq s,$$
$$\bar{a}_{ij} := a_{ij} + a_{is} a_{rj}$$
$$\mbox{ for } i = 0,1,\ldots, m-1, j = 0,1,\ldots,s-1,s+1,\ldots,n.$$
\item Let for $i = 0,1,\ldots,m; j = 0,1,\ldots,n; j \neq s$
$$a_{ij} := \bar{a}_{ij}$$
and go to 1.
\end{enumerate}
\end{algorithm}

\section{Automatic proofs}

We once more emphasize the fact that {\sf MOC} does
{\bf not} provide an algorithm for calculating the irreducible
modular characters for a given finite group. 
No such algorithm is yet known. 
Rather {\sf MOC} provides a collection of algorithms
for various purposes to support the calculation of some
decomposition numbers and irreducible Brauer characters. 

Finding the modular characters of a finite
group is a kind of problem where a result can be checked
with much less effort than was needed to obtain it.
We have to produce a very large number of Brauer characters
and projective characters to begin with. 
If the irreducible respectively indecomposable
characters happen to be among those of our pool,
then we would find them by means of the scalar product
(see Proposition~\ref{idealfall})\@.

The methods described in Section~\ref{generation} 
for constructing characters do not 
allow to select good candidates beforehand. 
Usually most of the characters generated are not needed at all. 
Only after the construction can one recognise those which really 
contain the information wanted. 
The other characters constructed during the session can be thrown away.
For the final proof, which can usually be given in a conventional
form, it suffices, of course, to keep the ``good guys''\@. 

Every run of {\sf MOC} is documented on an info-file corresponding
to the pair $(G,p)$\@. Some additional information,
for example output of certain programs is stored on this
file, which is called {\tt G.p.info}. This method of automatic documentation
is very much supported by the UNIX operating system.
As already mentioned, {\sf MOC} consists of a collection of 
FORTRAN programs and UNIX Bourne-shell scripts. 
The FORTRAN programs serve some well defined purposes
such as, for example, matrix multiplication or solving integral
linear equations. Transfer of data is via files, and the
files are moved by shell scripts according to the
needs of the programs.

Let us give an example of such a shell script.
In {\sf MOC}, the actual basic set ${\bf PS}$ of projective
characters is stored on the file {\tt G.p} under the label
30700 as a matrix~$X_0$, expressing ${\bf PS}$ 
in terms of the projective atoms ${\bf PA}_0$ 
corresponding to the special basic set ${\bf BS}_0$\@.
By equation~(\ref{PinPA0}) we have
$X_0 = \langle {\bf PS}, {\bf BS}_0 \rangle$\@.
The following shell script, called {\tt bsinba} transposes
the matrix~$X_0$ stored under label 30700 on the file
{\tt G.p}. By equation~(\ref{BS0inBA}) this expresses the
special basic set ${\bf BS}_0$ in terms of the Brauer atoms ${\bf BA}$
corresponding to ${\bf PS}$\@.
Since the FORTRAN program {\tt transp} to transpose a matrix
expects the matrix under label 30900, we have to {\tt merge} it first from
30700 to 30900. The shell script {\tt get} copies the file {\tt G.p}
to the current working directory, which for safety reasons
should always be different from the directory containing the
data-files. The first argument of the following shell script
is the prime~$p$, the second the name of the group~$G$\@. 

     \begin{verbatim}
     set -x
     # bsinba
     rm -f t1 t2 t3 t4
     get $2.$1 || exit
     cp $2.$1 t1
     cp t1 t2
     merge <<%
     30900
     0
     30700
     %
     rm -f t2
     mv t3 t1
     transp
     mv t2 t3
     rm -f t1 t2
     \end{verbatim}
     The shell script expressing the Brauer atoms in terms of the 
     special basic set ${\bf BS}_0$ is called {\tt bainbs}.
     By equation~(\ref{BS0inBA}) it has to invert the matrix~$X_0^t$\@.
     \begin{verbatim}
     set -x
     # 
     # bainbs 
     #
     # Expresses the Brauer atoms corresponding to 
     # the actual  basic set of projectives, which 
     # are stored in $2.$1 under the label 30700, 
     # in terms of the special basic set.
     # The result is put onto file t3.
     #
     # ``bsinba'' is the program inverse to ``bainbs''
     #
     bsinba $1 $2
     mv t3 hilf
     #
     # ``invers'' calculates the inverse of a matrix
     #       
     invers hilf
     rm t1 t2
     mv hilf.inv t3
     \end{verbatim}
     The following shell script calculates the set of Brauer atoms ${\bf BA}$
     as class functions. For this purpose it has to multiply the 
     matrix~$X_0^{-t}$ with $[{\bf BS}_0]$, the latter being
     stored on {\tt G.p} under the label 30900.
     \begin{verbatim}
     set -x
     #
     # bratoms
     #
     # This program gives the $1-modular 
     # Brauer atoms of group $2 corresponding
     # to the actual basic set of projective
     # characters as class functions.
     # The answer is stored in $2.$1.bat.
     #
     bainbs $1 $2
     mv t3 t2
     get $2.$1 1
     matmul2
     mv t3 $2.$1.bat
     rm -f $2.$1 t1 t2 t4
     \end{verbatim}
     We now give an example of an info-file which documents the 
     session for determining the Brauer characters for the
     Mathieu group $M_{11}$ modulo~$5$\@.
     \begin{verbatim}
     Fri Mar 20 11:25:16 MET 1992
     /home/euterpe/hiss/mocha/proggy/prp
     pst run
     fct run
      
     projectives no     1 -    5 obtained by:
     defect 0 characters
     
     relations stored under 30550
      
     Brauer characters no     1 -   10 obtained by:
     ordinary characters restricted to p-regular classes
     
     Fri Mar 20 11:25:56 MET 1992
     /home/euterpe/hiss/mocha/proggy/defzo
      
     projectives no     6 -   55 obtained by:
     tensoring ordinaries with defect 0 characters
     option:1:all ordinaries tensor all defect 0 characters
     
     Fri Mar 20 11:26:17 MET 1992
     /home/euterpe/hiss/mocha/proggy/basepro
     input label (i5)
     insert depth of search (i3), e.g. 10
     input block-number (i3)
     number of basic characters found is   4
     numbas= 10011 10012 10010 10007
     
     Fri Mar 20 11:27:21 MET 1992
     /home/euterpe/hiss/mocha/proggy/protest
     projective  nr    11 in block   1 is indecomposable,
       because it is an atom
     projective  nr    12 in block   1 is indecomposable,
       because it is an atom
     projective  nr    10 in block   1 is indecomposable,
       because it is an atom
     projective  nr     7 in block   1 is indecomposable,
       because it is an atom
     \end{verbatim}

The info-file serves two principal purposes. It can be used
to repeat the calculations for a particular group and prime.
This may be necessary if some information is lost by accident
or if the files get too large and immovable.
In a second run one can usually proceed much faster, because one already
knows which of the steps have been redundant.

The second purpose is of course the reconstruction of those
characters which are used in the proof. 
Every character generated is supplied with an identification
number. 
In particular, every irreducible Brauer character or
indecomposable projective character has its generating number, 
which can be used
to trace back to its origin with the help of the info-file.
We even have programs which do the trace-back automatically, at
least if the character was obtained by some standard constructions.
This gives the method of automatic proofs.
In \cite{BBB} about two thirds of the proofs
were obtained in this fashion.
If such a proof is written down in conventional form
everybody can reconstruct the characters
by the description we have given and thus check the proof 
by hand. An example of such a proof is given in Chapter~\ref{beweise},
where we have collected all the information necessary
to determine the decomposition numbers
modulo~$7$ of Conway's largest group $\mbox{\em Co}_1$\@.

\section{Extensions}

There are some features which we would like {\sf MOC} to have,
but which are still not implemented.
The first and most important one is the additional information one obtains from
duality and group automorphisms. For instance, a real valued
ordinary character must be contained in a PIM and its dual
with equal multiplicity. This simple observation and analogous
statements in the case of group automorphisms can be used
to improve the algorithms described in Section~\ref{Improvingbasicsets}
considerably. This is one of the next things we are planning
to add to {\sf MOC}\@.

The additional information available in characteristic~$2$ by
Fong's lemma (see \ref{fongslemma} below) has not been implemented yet.

There are examples such as the Rudvalis group in characteristic~$5$
which seem to be very hard to settle with the current methods
of {\sf MOC}\@. Here, no improvement of the characters as described
in Section~\ref{Improvingbasicsets} is possible, 
but still there seems to be enough
information to rule out most of the cases. 
Here, it really is desirable to have a program
which tries to solve Fundamental Problem~I in its full generality. 
We plan to do some experiments with these factorization
problems.

Finally, we would very much like to have a general
trace-back program for writing down more complicated
proofs. The current versions cannot cope with the numerous
possibilities for generating and improving characters.

\section{Advanced methods}
\label{AdvancedMethods}

As remarked in the introduction, {\sf MOC} uses only elementary
methods to calculate decomposition numbers. In large examples,
however, these methods are not sufficient to find all of the
irreducible Brauer characters. 

One of the reasons for this is the fact that there may be
proper ordinary characters which vanish on all $p$-singular
classes and thus are virtual projective characters, but which are not 
characters of projective modules.
This means that in the decomposition matrix there are 
two columns, corresponding to PIMs $\Phi$ and $\Psi$,
such that one is contained in the other, in other words,
$\Phi - \Psi$, say, is a proper ordinary character.
This happens for example in the case of the principal $3$-block of 
the second Janko group~$J_2$, which has the 
decomposition matrix displayed in Table~\ref{j2mod3dec} (for a proof see \cite{hisslux})\@.
\begin{table}[tb]
\caption{\label{j2mod3dec} The $3$-modular decomposition of $J_2$}
$$
\begin{array}{rcccccccc}\hline
        & \Phi_{ 1} & \Phi_{ 2} & \Phi_{ 3}
        & \Phi_{ 4} & \Phi_{ 5} & \Phi_{ 6}
        & \Phi_{ 7} & \Phi_{ 8}  \rule[- 7pt]{0pt}{ 20pt} \\ \hline
   1 &  1 &  . &  . &  . &  . &  . &  . & . \rule[ 0pt]{0pt}{ 13pt}\\
  14 &  1 &  1 &  . &  . &  . &  . &  . & . \\
  14 &  1 &  . &  1 &  . &  . &  . &  . & . \\
  21 &  . &  . &  . &  1 &  . &  . &  . & . \\
  21 &  . &  . &  . &  . &  1 &  . &  . & . \\
  70 &  . &  1 &  . &  . &  . &  1 &  . & . \\
  70 &  . &  . &  1 &  . &  . &  . &  1 & . \\
 160 &  1 &  1 &  1 &  . &  . &  . &  . & 1 \\
 175 &  . &  . &  . &  1 &  1 &  . &  . & 1 \\
 224 &  . &  1 &  . &  . &  1 &  1 &  . & 1 \\
 224 &  . &  . &  1 &  1 &  . &  . &  1 & 1 \\
 300 &  1 &  2 &  2 &  . &  . &  1 &  1 & 1 \\
 336 &  2 &  1 &  1 &  1 &  1 &  . &  . & 2 \rule[- 7pt]{0pt}{  5pt} \\ \hline
\end{array}
$$
\end{table}
In this example $\Phi_2 - \Phi_6$ is a proper ordinary character.
Table~\ref{tabbaschar} gives the character table of a basic set of Brauer characters.
The table is given in the {\sf MOC}-character table format.
We use the integral bases $\{1, (-1 + \sqrt{5})/2\}$ for the
occuring irrationalities.
\begin{table}[tb]
\caption{\label{tabbaschar} Some $3$-modular characters of $J_2$}
$$
\begin{array}{r@{\hspace{ 6pt}}*{13}{r@{\hspace{ 6pt}}}}\hline
   1A &  2A & 2A & 4A & 5A & 5A & 5C & 5C & 7A & 8A &  10A & 10A & 10C & 10C 
\rule[- 7pt]{0pt}{ 20pt} \\ \hline
   1 &  1 & 1 & 1 & 1 & 0 & 1 & 0 & 1 & 1 & 1 & 0 & 1 & 0 \rule[ 0pt]{0pt}{ 13pt} \\
  13 & -3 & 1 & 1 & 2 & 3 & 1 & 1 & -1 & -1 & -1 & 1 & 0 & 1 \\
  13 & -3 & 1 & 1 & -1 & -3 & 0 & -1 & -1 & -1 & -2 & -1 & -1 & -1 \\
  21 &  5 & -3 & 1 & 4 & 1 & 2 & 2 & 0 & -1 & -1 & -1 & 0 & 0 \\
  21 &  5 & -3 & 1 & 3 & -1 & 0 & -2 & 0 & -1 & 0 & 1 & 0 & 0 \\
  70 & -10 & -2 & 2 & 5 & 5 & 0 & 0 & 0 & 0 & 0 & -1 & 0 & 0 \\
  70 & -10 & -2 & 2 & 0 & -5 & 0 & 0 & 0 & 0 & 1 & 1 & 0 & 0 \\
 133 &  5 & 1 & -3 & -7 & 0 & -2 & 0 & 0 & 1 & 1 & 0 & 0 & 0 \\
  36 &  4 & 0 & 4 & -4 & 0 & 1 & 0 & 1 & 0 & 0 & 0 & -1 & 0 \\
  90 & 10 & 6 & -2 & 5 & 0 & 0 & 0 & -1 & 0 & 1 & 0 & 0 & 0 \\
  63 & 15 & -1 & 3 & 3 & 0 & -2 & 0 & 0 & 1 & -1 & 0 & 0 & 0 \\
 225 & -15 & 5 & -3 & 0 & 0 & 0 & 0 & 1 & -1 & 0 & 0 & 0 & 0 \\
 189 &  -3 & -3 & -3 & 3 & 3 & 2 & 1 & 0 & 1 & 0 & 1 & -1 & -1 \\
 189 &  -3 & -3 & -3 & 0 & -3 & 1 & -1 & 0 & 1 & -1 & -1 & 0 & 1 \rule[- 7pt]{0pt}{ 20pt} \\ \hline 
 169 &  9 & 1 & 1 & 13 & 3 & 2 & 1 & 1 & 1 & 2 & -3 & 1 & -1
\rule[- 7pt]{0pt}{  5pt} \\ \hline
\end{array}
$$
\end{table}

Let~$\varphi$ denote the irreducible Brauer character
corresponding to $\Phi$\@. The proof that $\Psi$ is not a 
direct summand of
the projective with character $\Phi$ is usually achieved
by giving a lower bound for the multiplicity of $\varphi$
in the reduction modulo~$p$ of some ordinary character. 

In the above example, $\varphi$ is an irreducible
Brauer character of degree~$13$, and to show that $\Phi_2 - \Phi_6$ 
is not projective, is the same as showing that~$\varphi$
is contained in the reduction modulo~$3$ of the first
ordinary character of degree~$70$\@. By inspecting the
table of Brauer characters~\ref{tabbaschar} we see that
$$ 169 = 13_1 \otimes 13_1 = 1 - 13_1 + 21_1 + 70_1 + 90,$$
and thus~$13_1$ must be contained in $70_1$\@.
This example is typical insofar as we have used tensor 
products of Brauer characters to obtain the bounds.
The larger the group and hence the
degrees of the irreducible Brauer characters, 
the weaker gets this method.
Sometimes it helps to distinguish cases as 
indicated in Section~\ref{fundamentalproblemI}
with the help of Brauer characters
obtained by tensor products. 

If all this fails, there are more advanced methods available, 
which make use of structure theory of $G$-modules, 
which we are now going to describe.
For every such method we also give an example where it can be
applied. In fact, all the tricks to be described in the following
have come to mind only when we were trying to
solve the corresponding problems.

Let $S$ be one of the rings $\{K,R,F\}$ of our $p$-modular
system for~$G$\@. We only consider $SG$-modules~$X$, which are free
and finitely generated as $S$-modules. Let~$X$ and~$Y$ be two such
$SG$-modules. We set
$$[X,Y] = \rk_S\,\Hom_{SG}(X,Y).$$
The $SG$-module dual to~$X$ is denoted by $X^*$\@. 

\subsection{Associativity of the tensor product}

Let~$X$,~$Y$ and~$Z$ denote three $SG$-modules. Then there
is a canonical isomorphism of $SG$-modules:
$$(X \otimes Y) \otimes Z \cong X \otimes (Y \otimes Z).$$
This can sometimes be used if the dimensions of the 
modules involved in the tensor product are small.

\subsection{Adjointness}

Let $X$, $Y$ and~$Z$ denote three $SG$-modules. We then have
(see \cite[Corollary~II.6.9]{laro})
$$[X \otimes Y, Z] = [X, Y^* \otimes Z].$$
This has been used in the final state of our example
of the Conway group modulo~$7$ 
(see Section~\ref{Co1mod7PrincBloPro}), and also 
in the proof for the Janko group
$J_2$ modulo~$5$ \cite{hisslux}\@.
This method can be used either to show that some composition
factor must occur in a certain tensor product ($J_2$ modulo~$5$)
or to show that it cannot occur ({\em Co}$_1$ modulo~$7$)\@.

It is worth noting that via tensor products
one can transfer information between different blocks
with the help of this adjointness property.
This is particularly helpful if some of the blocks
involved are of defect~$1$, since for those one can
write down all the indecomposables, if the Brauer tree
is known. An example of a successful application is provided
by the second Janko group $J_2$ modulo~$3$\@.

\subsection{Nakayama relations}

Let~$H$ be a subgroup of~$G$ and let $X$ be an $SG$-module,
$Y$ an $SH$ module. Then we have the Nakayama relations
(or Frobenius reciprocity) (see \cite[Corollary~II.1.4]{laro})
$$[X, Y^G] = [X_H, Y], \quad \mbox{ and } \quad [Y^G,X] = [Y, X_H].$$
This has been applied for instance in the calculation
of the modular characters for the Tits group (\cite{hisstits})\@.

\subsection{Self duality}

If an $FG$-module~$X$ is self dual, i.e., $X \cong X^*$,
then its socle series is the dual of its Loewy series
(see \cite[Lemma~I.8.4(i)]{laro})\@. 
In particular, any self dual composition factor $M$ in the head of~$X$
must also occur in the socle. 
Hence such an~$M$ either is twice a composition factor of~$X$ or 
else it is a direct summand.
If the latter is not the case, $X$ has a
non-trivial endomorphism.

The following lemma, which follows from the above remarks,
is often quite useful.
\begin{lem}
\label{lemma551}
Let~$X$ be a self dual $FG$-module with $[X,X] = 1$\@.
If every composition factor of~$X$ is self dual,
then~$X$ is irreducible. $\Box$
\end{lem}

\subsection{Trivial source modules}

The methods described above are very powerful in connection
with the theory of modules with a trivial source. 
These are, by definition, indecomposable direct summands
of permutation modules over~$F$ or~$R$\@.
Every indecomposable $FG$-module~$M$ with trivial source
is liftable to a uniquely determined trivial source $RG$-module~$X$
(see \cite[Section~II.12]{laro})\@. 
This means that $M$ is isomorphic to the
$FG$-module $X \otimes_{R} FG$\@. Furthermore, $[M,M] = [X,X]$\@.

Any direct summand of the tensor product of two  trivial source modules
is again a trivial source module. 
If~$H$ is a subgroup of~$G$ and~$Y$ a trivial source module of~$H$,
than every direct summand of $Y^G$ is a trivial source module.

If~$B$ is a block of~$G$ containing only real valued ordinary
characters, then~$B$ is called {\sl strongly real}\@.
\index{strongly real}
In a strongly real block every simple $FG$-module is self dual.
This follows from the surjectivity of the decomposition 
homomorphism: every irreducible Brauer character in such a block
is real valued and thus corresponds to a self dual simple module.

Let~$M$ be a trivial source $FG$-module with trivial source
lift~$X$\@. Suppose that the ordinary character of~$X$ is
irreducible and that~$M$ is contained in a strongly real block.
Then~$M$ is irreducible by Lemma~\ref{lemma551}\@. This applies in particular to 
permutation modules of doubly transitive permutation 
representations on~$n$ points if~$p$ does not divide~$n$\@.
If $1 + \chi$ is the ordinary character of a doubly transitive
permutation representation, then $\chi$ is irreducible modulo~$p$
if $\chi$ is contained in a strongly real block
and if $p$ does not divide $1 + \chi(1)$\@.
The symplectic group $S_6(2)$ has a doubly transitive
permutation representation of degree~$28$, and, since the
ordinary character table of $S_6(2)$ has only rational entries,
the non-trivial constituent of the permutation character remains
irreducible modulo~$3$\@.

\subsection{Lattices}

An $RG$-module~$X$, which is free as $R$-module, is also called
an $RG$-{\sl lattice}\@. A sublattice~$Y$ of~$X$ is called {\sl pure},
\index{lattice}
if the factor module $X/Y$ again is a lattice. The following result
is more or less folclore.

\begin{lem}
Let~$X$ be an $RG$-lattice with
ordinary character $\chi + \psi$\@. Then there exists
a pure sublattice $Y \leq X$ such that the character
of~$Y$ is $\chi$ and the character of $X/Y$ is $\psi$\@.
\end{lem}
\begin{bew} See \cite[Theorem~I.17.3]{laro}\@.
\end{bew}

\noindent This has turned out to very powerful in the theory
of blocks with a cyclic defect group.
But it has also found applications for non-cyclic
defect groups, for example in the proof of Theorem~$1$ of \cite{hisstits}\@.

\subsection{Fong's lemma}

\begin{lem}
\label{fongslemma}
{\rm (Fong's lemma:)} Let $p = 2$\@. Then every self dual non-trivial
simple $FG$-module has even dimension.
\end{lem}
\begin{bew}
Such a module carries a non-degenerate $G$-invariant symplectic
form (see \cite[Theorem~VII.8.13]{HuBlII})\@.
\end{bew}

\noindent Let $\chi \in \Irr(G)$ be real-valued.
Then it easily follows from Fong's lemma that 
$d_{\chi 1} \equiv \chi(1)(\mod\,2)$\@.

For example, consider the Mathieu group $M_{11}$ modulo~$2$\@.
There are two characters of defect~$0$, namely $\chi_6$ and $\chi_7$
(in Atlas notation)\@. If we define
$$ \Psi   _{1  } :=  \chi   _{7} \otimes        \chi   _{6},$$
$$ \Psi   _{2  } :=  \chi   _{6} \otimes        \chi   _{6},$$
$$ \Psi   _{3  } :=  \chi   _{2} \otimes        \chi   _{6},$$
we obtain the following basic set of projective characters for the 
principal $2$-block of $M_{11}$\@.
$$
\begin{array}{rrrr}\hline
        & \Psi_{ 1} & \Psi_{ 2} & \Psi_{ 3}
        \rule[- 7pt]{0pt}{ 20pt} \\ \hline
  1 &  1 &  . & . \rule[ 0pt]{0pt}{ 13pt} \\
 10 &  . &  1 & . \\
 10 &  . &  1 & . \\
 10 &  . &  1 & . \\
 11 &  1 &  1 & . \\
 44 &  1 &  1 & 1 \\
 45 &  2 &  1 & 1 \\
 55 &  2 &  2 & 1 
\rule[- 7pt]{0pt}{  5pt} \\ \hline
\end{array}
$$
Fong's lemma shows that $\Psi_3$ is contained in $\Psi_1$\@.

\subsection{Theorem of Benson and Carlson}

\begin{satz}
{\rm (Benson-Carlson:)} Let~$M$,~$N$ be indecomposable $FG$-mod\-ul\-es
such that $p \mid \dim_k\,M$\@. Then every indecomposable direct
summand of $M \otimes N$ has dimension divisible by~$p$\@.
\end{satz}
\begin{bew}
See \cite[Vol.~I,Theorem~3.1.9]{benbook}\@,
\end{bew}

\noindent This is applied, for example, in
\cite[p.~$111$]{hilupa}\@. 

\subsection{Blocks with cyclic defect groups}

In a block with a cyclic defect group one can describe
all indecomposable $FG$-modules, once the planar embedded
Brauer tree is known (See \cite[Chapter~V]{Alperin})\@. 

This information is usually easier to obtain than that for
a block with a non-cyclic defect group. It is very powerful
in connection with the Nakayama relations (where we assume
that some block for a subgroup with a cyclic defect group
is known) and, if some blocks with a cyclic defect group
for~$G$ are known, in connection with the adjointness or
the theorem of Benson-Carlson.

This has been applied,
for example, for some blocks of the Monster group
in characteristic~$5$ (see~\cite[p.~454f]{BBB}),
and also in the Rudvalis group modulo~$3$ \cite{hissru}\@.

\subsection{Webb's theorem}

Let~$M$ be an indecomposable $FG$-module. The heart of~$M$
is the $FG$-module $H(M) = \mbox{\rm rad}(M)/\mbox{\rm soc}(M)$,
where $\mbox{\rm rad}(M)$ denotes the intersection of all maximal
submodules of~$M$ and $\mbox{\rm soc}(M)$ is the sum of all simple
submodules. It is not difficult to see that if~$M$ is self dual,
the same is true for the heart of~$M$\@. By classifying the possible
Auslander-Reiten quivers for a finite-dimensional group algebra,
Webb has shown that the heart of a projective indecomposable
module of~$FG$ has at most~$4$ indecomposable direct summands.
This result has been extended by Bessenrodt.

\begin{satz}
{\rm (Webb, Bessenrodt:)} Let~$P$ be a projective indecomposable
module of~$FG$\@. Then $H(P)$ has at most~$3$
indecomposable direct summands.
\end{satz}
\begin{bew}
See \cite[Corollary~1.3]{bessenrodt}\@.
\end{bew}

We can apply this in some situations as follows.

\begin{lem}
Let~$B$ be a strongly real $p$-block of~$G$\@. Let~$\Phi$
be a PIM of~$B$ such that $\langle \Phi, \Phi \rangle < 4$ and
\begin{equation}
\label{grossesherz}
\sum_{\Psi \in {\rm IPr}(B)} \langle \Phi, \Psi \rangle \geq 6.
\end{equation}
Then there exists some $\Theta \in \IPr(B)$, $\Theta \neq \Phi$
such that $\langle \Phi, \Theta \rangle \geq 2$\@.
\end{lem}
\begin{bew}
Let~$P$ denote the projective indecomposable $FG$-module with
character~$\Phi$\@.
By assumption~(\ref{grossesherz}), the heart $H(P)$ of~$P$
has at least~$4$ composition factors. 
By Bessenrodt's extension of Webb's theorem, $H(P)$ has at
most~$3$ indecomposable direct summands, hence it cannot be
semi-simple. It follows that $H(P)$ contains some
composition factor with multiplicity at least~$2$\@. 
Since $\langle \Phi, \Phi \rangle < 4$, this must be a
simple module not isomorphic to $P/\mbox{\rm rad}(P)$
and the assertion follows.
\end{bew}

For example, this lemma can be applied in the case of the
Rudvalis group modulo~$3$ \cite{hissru}\@.
We finally proof a lemma in the same spirit as the one above.
\begin{lem}
Suppose that~$G$ has no factor group of order~$p$\@.
Suppose also, that the principal block of~$G$ is strongly real.
Let~$\Phi$ denote the character of the projective cover~$P$ of 
the trivial $FG$-module. If $\langle \Phi, \Phi \rangle > 2$,
then there exists some $\Theta \in \IPr(B)$, $\Theta \neq \Phi$ 
such that $\langle \Phi, \Theta \rangle \geq 2$\@. 
\end{lem}
\begin{bew} 
If the trivial $FG$-module is in the second Loewy layer of~$P$,
then~$G$ has some factor group of order~$p$ 
(see \cite[Corollary~I.10.13]{laro})\@.
It follows that $H(P)$ is not semi-simple and that it must
contain some non-trivial composition factor with multiplicity
at least~$2$\@.
\end{bew}


\clearpage
\thispagestyle{empty}

\chapter{Calculating the $7$-modular decomposition matrices of the Conway group}
\markboth{DECOMPOSITION MATRICES OF CONWAY'S GROUP}{THE FAITHFUL BLOCK OF MAXIMAL DEFECT}
\label{beweise}

The double covering group of the Conway group {\em Co}$_1$ has nine
$7$-blocks of defect larger than~$0$\@. Seven of these are of defect~$1$
and their decomposition matrices are given in~\cite{BBB}\@. 
The remaining two blocks are of maximal defect~$2$\@. 
Each contains $29$ ordinary and $21$ irreducible 
Brauer characters.

\section{The faithful block of maximal defect}
\markboth{DECOMPOSITION MATRICES OF CONWAY'S GROUP}{THE FAITHFUL BLOCK OF MAXIMAL DEFECT}

We start with the block of maximal defect containing faithful characters,
shortly called the faithful block in the following. 
The proof for the faithful block is easier than that for the principal block, 
and the result is used in the 
proof for the principal block.

We start with the set of projective characters displayed in Table~\ref{tab1}\@.
The origin of the projective characters is documented in table~\ref{origin}\@.
The symbol $\theta_i$ denotes the $i$th character of {\em Co}$_2$
(in ATLAS ordering), and a bar denotes the character times
the sign character of $2 \times \mbox{\em Co}_2$\@.

The two characters of degree $9\,152\,000$ are complex conjugates of
each other. All other characters in the block are real.
Since the number of irreducible Brauer characters in
the block is~$21$, there are exactly~$19$ real valued projective
indecomposable characters. Each of $\Lambda_{ 1}$--$\Lambda_{20}$,
except~$\Lambda_{12}$, contains one of these. Therefore,~$\Lambda_{12}$
contains a pair of complex conjugate projective indecomposable
characters. In particular, the decomposition matrix has wedge shape.
From this it immediately follows, that $\Lambda_{ 2}$--$\Lambda_{ 4}$,
$\Lambda_{ 7}$, $\Lambda_{ 9}$, $\Lambda_{11}$ and
$\Lambda_{14}$--$\Lambda_{20}$ are projective indecomposable
characters.
 
We now consider the projectives given in Table~\ref{tab2}\@.
Table~\ref{origin} gives their origin and
Table~\ref{expr1} their decompositions in terms of the projectives of
Table~\ref{tab1}, the relations.

Using these relations, one can reduce the first projectives to obtain a new
set of projectives, given in Table~\ref{tab3}\@.
For example, $\Lambda_{30}$ shows, that the PIM $\Lambda_{18}$ is contained
in $\Lambda_{1}$, and $\Lambda_{31}$ shows, that the PIM $\Lambda_{3}$
is contained in $\Lambda_{1}$\@. This gives~$\Phi_1$ of Table~\ref{tab3}\@.
All the other new projectives are obtained by similar
arguments.
Since the two characters of degree 9\,152\,000 are complex conjugates 
of each other, $\Phi_{12}$
contains a pair of complex conjugate PIMs. 
This immediately determines the decomposition matrix.
 
\section{The principal block}
\markboth{DECOMPOSITION MATRICES OF CONWAY'S GROUP}{THE PRINCIPAL BLOCK}
 
\label{Co1mod7PrincBloPro}
 
We start with the set of projectives displayed in Table~\ref{tab4}\@.
The origin of these is documented in Table~\ref{orig2}\@.
It is not too difficult to check that
$\Lambda_{ 2}$, $\Lambda_{ 4}$,
$\Lambda_{ 6}$--$\Lambda_{ 9}$, $\Lambda_{12}$,
$\Lambda_{15}$ and $\Lambda_{19}$ are projective indecomposable
characters. (For example, this can be done by checking
that no subsum of these characters is zero on all $7$-singular
classes.)
 
We have displayed further projectives in Table~\ref{tab5},
their origin in Table~\ref{orig2}, where the notation
$\Phi_{i,j}$ stands for the $i$-th projective indecomposable
character in Block~$j$\@. The numbering of the blocks follows
Appendix~\ref{Co1mod7}\@.
The expression of the new projectives in 
terms of the first set of projectives is given in Table~\ref{expr2}\@.
Using these relations,
one can reduce the basic projectives to obtain a new
set of projectives, given in Table~\ref{tab6}\@.
The arguments are exactly the same as those used for the proof of Block~$7$\@.

One immediately checks that all of these, except possibly $\Phi_{10}$ and $\Phi_{14}$
are characters of projective indecomposable modules. Furthermore, either $\Phi_{10}$
is a PIM or $\Phi_{10} - \Phi_{17}$ is one. Similarly, either $\Phi_{14}$ is a PIM
or $\Phi_{14} - \Phi_{20}$ is one.
 
Thus we are left with four possible decomposition matrices. To decide which of these is
correct, we reformulate the problem in terms of Brauer characters. All irreducible
Brauer characters are known, except $\phi_{17}$ and $\phi_{20}$\@.
The degree of $\phi_{10}$ is 2\,038\,674, that of $\phi_{14}$ is 10\,140\,998\@.
If $\Phi_{10}$ is a PIM, then the degree of $\phi_{17}$ is 62\,725\,301, if
$\Phi_{10} - \Phi_{17}$ is a PIM we have $\phi_{17}(1) = 62\,725\,301 + 2\,038\,674 =
64\,763\,975$\@.
Similarly, if $\Phi_{14}$ is a PIM, then the degree of $\phi_{20}$ is 124\,375\,559 if
$\Phi_{14} - \Phi_{20}$ is a PIM we have $\phi_{20}(1) = 124\,375\,559 + 10\,140\,998
= 134\,516\,557$\@.
 
Suppose first, that $\Phi_{10}$ is a PIM. Consider the irreducible Brauer character
of Block~$7$ of degree 50\,207\,872\@. Denoting the irreducible Brauer characters by their
degrees, we obtain
\begin{equation}
\label{gl1}
 24 \otimes 2\,038\,674 = 48\,928\,176,
\end{equation}
\begin{equation}
\label{gl2}
 24 \otimes 50\,207\,872 = 62\,725\,301 + 2\,038\,674 + \psi,
\end{equation}
where~$\psi$ is a sum of characters not in the principal block,
namely $\psi = 52\,465\,644 + 256\,168\,549 + 309\,429\,120 + 522\,161\,640$\@.
Let~$F$ denote an algebraically closed field of characteristic~$7$\@.
For two $FG$-modules~$M$ and~$N$, the dimension of the space of
$FG$-homomorphisms from~$M$ to~$N$ is denoted by $[M,N]$\@.
With this notation we have,
$$[24 \otimes 2\,038\,674, 50\,207\,872] = 0,$$
by equation~(\ref{gl1}), and so,
\begin{equation}
\label{gl4}
[2\,038\,674, 24 \otimes 50\,207\,872] = 0.
\end{equation}
Since all modules considered are self-dual, equation~(\ref{gl2}) implies
that the irreducible 2\,038\,674 is in the socle of $24 \otimes 50\,207\,872$,
contradicting equation~(\ref{gl4})\@.
We have proved that $\Phi_{10}$ is decomposable.
 
Now suppose that $\Phi_{14}$ is a PIM. This time we consider the irreducible
Brauer character of degree 10\,039\,568 of Block~$7$\@. We obtain:
$$ 24 \otimes 10\,140\,998 = 243\,383\,952,$$
\begin{equation}
\label{gl6}
24 \otimes 10\,039\,568 = 124\,375\,559 + 10\,140\,998 + \psi,
\end{equation}
where~$\psi$ is a sum of characters not in the principal block,
namely $\psi = 77\,700\,854 + 25\,892\,020 +  2\,840\,201$\@.
However, as above,
$$[10\,140\,998, 24 \otimes 10\,039\,568] = 0,$$
contradicting equation~(\ref{gl6})\@.
We have shown that $\Phi_{14}$ decomposes, and thus completed the proof of
the decomposition matrix.

\pagebreak

\renewcommand{\textfraction}{0.0}

{\samepage
\section{Tables}
\markboth{DECOMPOSITION MATRICES OF CONWAY'S GROUP}{TABLES}

Here we collect the tables of projective characters
which are needed in the previous sections.

{\small
\begin{table}[h]
\caption{\label{tab1} Some projectives in the faithful block}
$$\begin{array}{r@{\hspace{ 7pt}}*{10}{r@{\hspace{ 5pt}}}*{10}{r@{\hspace{ 4pt}}}} \hline
 \Lambda:  &  1 &  2 &  3 &  4 &  5 &  6 &  7 &  8 &  9 & 10 &
 11 & 12 & 13 & 14 & 15 & 16 & 17 & 18 & 19 & 20 \rule[- 7pt]{0pt}{ 20pt} \\ \hline
 24 &  1 &  . &  . &  . &  . &  . &  . &  . &  . &  . &  . &  . &  . &  . &  . &  . &
 . &  . &  . & . \rule[  0pt]{0pt}{ 13pt} \\
 2024 &  . &  1 &  . &  . &  . &  . &  . &  . &  . &  . &  . &  . &  . &  . &  . &  .
&  . &  . &  . & . \\
 4576 &  1 &  . &  1 &  . &  . &  . &  . &  . &  . &  . &  . &  . &  . &  . &  . &  .
&  . &  . &  . & . \\
 40480 &  . &  1 &  . &  1 &  . &  . &  . &  . &  . &  . &  . &  . &  . &  . &  . &  . &  . &  . &  . & . \\
 95680  &  1 &  . &  . &  . &  1 &  . &  . &  . &  . &  . &  . &  . &  . &  . &  . &  . &  . &  . &  . & . \\
 299000 &  1 &  . &  1 &  . &  . &  1 &  . &  . &  . &  . &  . &  . &  . &  . &  . &
. &  . &  . &  . & . \\
 315744 &  1 &  . &  1 &  . &  . &  . &  1 &  . &  . &  . &  . &  . &  . &  . &  . &
. &  . &  . &  . & . \\
 789360 &  . &  . &  . &  . &  . &  1 &  . &  1 &  . &  . &  . &  . &  . &  . &  . &
. &  . &  . &  . & . \\
 1937520 &  1 &  . &  . &  . &  . &  . &  . &  . &  1 &  . &  . &  . &  . &  . &  . &
 . &  . &  . &  . & . \\
 5051904 &  . &  . &  . &  1 &  . &  . &  . &  . &  1 &  1 &  . &  . &  . &  . &  . &
 . &  . &  . &  . & . \\
 7104240 &  . &  . &  . &  . &  1 &  . &  . &  . &  . &  . &  1 &  . &  . &  . &  . &
 . &  . &  . &  . & . \\
 9152000 &  . &  . &  . &  . &  . &  . &  . &  . &  . &  . &  . &  1 &  . &  . &  . &
 . &  . &  . &  . & . \\
 9152000 &  . &  . &  . &  . &  . &  . &  . &  . &  . &  . &  . &  1 &  . &  . &  . &
 . &  . &  . &  . & . \\
 13156000 &  . &  1 &  . &  1 &  . &  . &  . &  . &  . &  1 &  . &  . &  1 &  . &  . &  . &  . &  . &  . & . \\
 17050176  &  . &  1 &  . &  . &  . &  . &  . &  . &  . &  . &  1 &  . &  1 &  . &  . &  . &  . &  . &  . & . \\
 19734000 &  1 &  . &  . &  . &  . &  . &  . &  . &  1 &  2 &  . &  . &  . &  1 &  . &  . &  . &  . &  . & . \\
 34155000 &  1 &  . &  1 &  . &  . &  1 &  1 &  . &  . &  2 &  . &  . &  . &  . &  1 &  . &  . &  . &  . & . \\
 49335000 &  . &  . &  . &  . &  2 &  1 &  1 &  . &  . &  . &  . &  . &  . &  . &  . &  1 &  . &  . &  . & . \\
 50519040 &  . &  . &  . &  . &  . &  . &  1 &  . &  . &  . &  . &  1 &  1 &  . &  . &  . &  1 &  . &  . & . \\
 67358720 &  . &  . &  . &  . &  . &  1 &  . &  1 &  . &  3 &  . &  2 &  . &  1 &  1 &  . &  . &  . &  . & . \\
 210496000 &  1 &  . &  . &  . &  1 &  1 &  . &  . &  . &  . &  1 &  . &  1 &  . &  .
&  . &  . &  1 &  . & . \\
 215547904 &  . &  . &  . &  . &  . &  1 &  . &  3 &  . &  1 &  . &  6 &  2 &  . &  .
&  . &  . &  . &  1 & . \\
 313524224 &  1 &  . &  . &  . &  . &  . &  . &  2 &  . &  4 &  . &  4 &  . &  1 &  1
&  . &  . &  . &  . & 1 \\
 394680000 &  . &  . &  . &  . &  1 &  1 &  1 &  1 &  . &  2 &  . &  5 &  1 &  . &  1
&  1 &  1 &  . &  . & 1 \\
 485760000 &  2 &  . &  . &  . &  2 &  2 &  . &  1 &  . &  . &  . &  2 &  . &  . &  .
&  1 &  . &  1 &  . & 1 \\
 517899096 &  . &  . &  . &  . &  . &  1 &  . &  5 &  . &  1 &  . &  11 &  2 &  . &  . &  . &  1 &  . &  1 & 1 \\
 655360000 &  1 &  . &  . &  . &  . &  2 &  . &  4 &  . &  2 &  . &  8 &  2 &  . &  .
&  . &  . &  1 &  1 & 1
\rule[- 7pt]{0pt}{  5pt} \\ \hline
\end{array}
$$
\end{table}
}
}
\begin{table}[tb]
\caption{\label{tab2} More projectives in the faithful block}
$$\begin{array}{r@{\hspace{10pt}}*{13}{r@{\hspace{ 5pt}}}} \hline
 \Lambda:  & 21 & 22 & 23 & 24 & 25 & 26 & 27 & 28 & 29 & 30 &
 31 & 32  & \rule[- 7pt]{0pt}{ 20pt} \\ \hline
 24 &  . &  . &  . &  . &  . &  . &  . &  . &  . &  1 &  1 &  . &  \rule[  0pt]{0pt}{ 13pt} \\
 2024 &  . &  . &  . &  . &  . &  . &  . &  . &  . &  1 &  . &  . & \\
 4576 &  . &  . &  . &  . &  . &  . &  . &  . &  . &  2 &  . &  . & \\
 40480 &  . &  . &  . &  . &  . &  1 &  . &  . &  . &  1 &  . &  . & \\
 95680 &  . &  . &  . &  . &  1 &  . &  . &  . &  . &  1 &  1 &  . & \\
 299000 &  2 &  . &  . &  . &  . &  3 &  . &  . &  . &  4 &  . &  . & \\
 315744 &  . &  . &  . &  . &  1 &  . &  . &  . &  . &  2 &  . &  . & \\
 789360 &  2 &  . &  . &  . &  . &  4 &  . &  . &  . &  2 &  . &  1 & \\
 1937520 &  3 &  1 &  . &  . &  . &  7 &  . &  . &  1 &  5 &  1 &  . & \\
 5051904 &  3 &  2 &  . &  1 &  . &  9 &  1 &  1 &  1 &  4 &  . &  1 & \\
 7104240 &  . &  . &  2 &  . &  2 &  . &  . &  . &  . &  . &  . &  . & \\
 9152000 &  . &  . &  . &  . &  . &  . &  1 &  1 &  1 &  . &  . &  . & \\
 9152000 &  . &  . &  . &  . &  . &  . &  1 &  1 &  1 &  . &  . &  . & \\
 13156000 &  . &  1 &  1 &  2 &  . &  2 &  2 &  2 &  . &  1 &  . &  1 & \\
 17050176 &  . &  . &  3 &  1 &  1 &  . &  1 &  1 &  . &  1 &  . &  . & \\
 19734000 &  4 &  5 &  . &  3 &  . &  14 &  2 &  1 &  2 &  6 &  1 &  2 & \\
 34155000 &  2 &  4 &  . &  4 &  1 &  7 &  . &  1 &  1 &  6 &  1 &  1 & \\
 49335000 &  . &  . &  . &  . &  2 &  . &  . &  . &  . &  . &  2 &  . & \\
 50519040 &  . &  . &  . &  . &  1 &  . &  2 &  1 &  . &  . &  . &  . & \\
 67358720 &  3 &  7 &  . &  6 &  . &  14 &  3 &  3 &  4 &  5 &  1 &  3 & \\
 210496000 &  . &  . &  6 &  1 &  2 &  . &  4 &  3 &  . &  . &  4 &  . & \\
 215547904 &  . &  . &  2 &  1 &  . &  . &  11 &  6 &  1 &  . &  1 &  1 & \\
 313524224 &  1 &  10 &  1 &  10 &  . &  13 &  8 &  7 &  5 &  4 &  6 &  4 & \\
 394680000 &  . &  6 &  1 &  7 &  1 &  5 &  8 &  7 &  4 &  2 &  7 &  1 & \\
 485760000 &  . &  2 &  4 &  3 &  1 &  1 &  7 &  5 &  1 &  1 &  11 &  . & \\
 517899096 &  . &  2 &  2 &  3 &  . &  2 &  18 &  11 &  4 &  . &  5 &  2 & \\
 655360000 &  . &  3 &  6 &  5 &  . &  2 &  19 &  12 &  2 &  . &  9 &  2 &
\rule[- 7pt]{0pt}{  5pt} \\ \hline
\end{array}
$$
\end{table}
\begin{table}[t]
\caption{\label{expr1} Relations in the faithful block}
$$\begin{array}{r@{\hspace{ 6pt}}*{20}{r@{\hspace{ 4.5pt}}}} \hline
 \Lambda:  &  1 &  2 &  3 &  4 &  5 &  6 &  7 &  8 &  9 & 10
 & 11 & 12 & 13 & 14 & 15 & 16 & 17 & 18 & 19 & 20
 \rule[- 7pt]{0pt}{ 20pt} \\ \hline
 \Lambda_{21} & . &  . &  . &  . &  . &  2 &  . &  . &  3 &  . &  . &  . &  . &  1 &
 . &  -2 &  . &  -2 &  -2 & .\rule[  0pt]{0pt}{ 13pt} \\
 \Lambda_{22} & . &  . &  . &  . &  . &  . &  . &  . &  1 &  1 &  . &  . &  . &  2 &
 2 &  . &  . &  . &  -1 & 2 \\
 \Lambda_{23} & . &  . &  . &  . &  . &  . &  . &  . &  . &  . &  2 &  . &  1 &  . &
 . &  . &  -1 &  3 &  . & 1 \\
 \Lambda_{24} & . &  . &  . &  . &  . &  . &  . &  . &  . &  1 &  . &  . &  1 &  1 &
 2 &  . &  -1 &  . &  -2 & 3 \\
 \Lambda_{25} & . &  . &  . &  . &  1 &  . &  1 &  . &  . &  . &  1 &  . &  . &  . &
 . &  -1 &  . &  . &  . & . \\
 \Lambda_{26} & . &  . &  . &  1 &  . &  3 &  . &  1 &  7 &  1 &  . &  . &  . &  5 &
 2 &  -3 &  . &  -3 &  -7 & . \\
 \Lambda_{27} & . &  . &  . &  . &  . &  . &  . &  . &  . &  1 &  . &  1 &  1 &  . &
 -2 &  . &  . &  3 &  2 & 2 \\
 \Lambda_{28} & . &  . &  . &  . &  . &  . &  . &  . &  . &  1 &  . &  1 &  1 &  -1 &  -1 &  . &  -1 &  2 &  -3 & 1 \\
 \Lambda_{29} & . &  . &  . &  . &  . &  . &  . &  . &  1 &  . &  . &  1 &  . &  1 &
 1 &  . &  -1 &  . &  -5 & -1 \\
 \Lambda_{30} & 1 &  1 &  1 &  . &  . &  2 &  . &  . &  4 &  . &  . &  . &  . &  1 &
 2 &  -2 &  . &  -3 &  -2 & . \\
 \Lambda_{31} & 1 &  . &  -1 &  . &  . &  . &  . &  . &  . &  . &  . &  . &  . &  . &  1 &  2 &  . &  3 &  1 & 4 \\
 \Lambda_{32} & . &  . &  . &  . &  . &  . &  . &  1 &  . &  1 &  . &  . &  . &  . &
 -1 &  . &  . &  . &  -3 & -1
\rule[- 7pt]{0pt}{  5pt} \\ \hline
\end{array}
$$
\end{table}
\begin{table}[b]
\caption{\label{origin} Origin of projective characters of the faithful block}
$$\begin{array}[t]{l@{\hspace{ 6pt}}cc} \hline
\mbox{\rm Char.} & \multicolumn{1}{c}{\mbox{\rm Origin}}
\rule[- 7pt]{0pt}{ 20pt} \\ \hline
   \Lambda   _{ 1 }: & {\Ind}_{\rm 2xCo2} (\bar{\theta}_{ 2} + \bar{\theta}_{21}) & \rule[  0pt]{0pt}{ 13pt} \\ 
   \Lambda   _{ 2 }: & {\Ind}_{\rm 2xCo2} (\bar{\theta}_{ 5}) & \\ 
   \Lambda   _{ 3 }: &    \chi   _{3} \otimes        \chi   _{111} & \\ 
   \Lambda   _{ 4 }: & {\Ind}_{\rm 2xCo2} (\bar{\theta}_{ 9}) & \\ 
   \Lambda   _{ 5 }: & {\Ind}_{\rm 2xCo2} (\bar{\theta}_{18}) & \\ 
   \Lambda   _{ 6 }: & {\Ind}_{\rm 2xCo2} (\bar{\theta}_{15} + \bar{\theta}_{24}) & \\ 
   \Lambda   _{ 7 }: &    \chi   _{104} \otimes      \chi   _{16} & \\ 
   \Lambda   _{ 8 }: & {\Ind}_{\rm 2xCo2} (   \bar{\theta}_{22} + \bar{\theta}_{23} + \bar{\theta}_{24}) & \\ 
   \Lambda   _{ 9 }: &    \chi   _{102} \otimes      \chi   _{36} & \\ 
   \Lambda   _{10 }: &    \chi   _{2} \otimes        \chi   _{141} & \\ 
   \Lambda   _{11 }: &    \chi   _{102} \otimes      \chi   _{39} & \\ 
   \Lambda   _{12 }: &    \chi   _{5} \otimes        \chi   _{118} & \\ 
   \Lambda   _{13 }: &    \chi   _{104} \otimes      \chi   _{33} & \\ 
   \Lambda   _{14 }: &    \chi   _{102} \otimes      \chi   _{49} & \\ 
   \Lambda   _{15 }: &    \chi   _{2} \otimes        \chi   _{127} & \\ 
   \Lambda   _{16 }: &    \chi   _{3} \otimes        \chi   _{118} & 
\rule[- 7pt]{0pt}{  5pt} \\ \hline
\end{array}
\quad
\begin{array}[t]{l@{\hspace{ 6pt}}cc} \hline
\mbox{\rm Char.} & \multicolumn{1}{c}{\mbox{\rm Origin}}
\rule[- 7pt]{0pt}{ 20pt} \\ \hline
   \Lambda   _{17 }: &    \chi   _{2} \otimes        \chi   _{118} & \rule[  0pt]{0pt}{ 13pt} \\
   \Lambda   _{18 }: &    \chi   _{2} \otimes        \chi   _{128} & \\
   \Lambda   _{19 }: &    \chi   _{2} \otimes        \chi   _{125} & \\
   \Lambda   _{20 }: &    \chi   _{102} \otimes      \chi   _{75}  & \\
   \Lambda   _{21 }: &    \chi   _{2} \otimes        \chi   _{119}     & \\ 
   \Lambda   _{22 }: &    \chi   _{2} \otimes        \chi   _{147}     & \\
   \Lambda   _{23 }: &    \chi   _{2} \otimes        \chi   _{148}     & \\
   \Lambda   _{24 }: &    \chi   _{2} \otimes        \chi   _{149}     & \\
   \Lambda   _{25 }: &    \chi   _{3} \otimes        \chi   _{131}     & \\
   \Lambda   _{26 }: &    \chi   _{4} \otimes        \chi   _{132}     & \\
   \Lambda   _{27 }: & {\Ind}_{\rm 2xCo2} (\bar{\theta}_{50})          & \\
   \Lambda   _{28 }: &    \chi   _{3} \otimes        \chi   _{156}     & \\
   \Lambda   _{29 }: &    \chi   _{4} \otimes        \chi   _{127}     & \\
   \Lambda   _{30 }: &    \chi   _{7} \otimes        \chi   _{111}     & \\
   \Lambda   _{31 }: & {\Ind}_{\rm 2xCo2} (\bar{\theta}_{ 1} + \bar{\theta}_{44})   & \\
   \Lambda   _{32 }: &    \chi   _{3} \otimes        \chi   _{141}     &
\rule[- 7pt]{0pt}{  5pt} \\ \hline
\end{array}
$$
\end{table}
\begin{table}[tb]
\caption{\label{tab3} Refined projectives in the faithful block}
$$\begin{array}{r@{\hspace{ 7pt}}*{10}{r@{\hspace{ 5pt}}}*{10}{r@{\hspace{ 4pt}}}} \hline
  \Phi:  &  1 &  2 &  3 &  4 &  5 &  6 &  7 &  8 &  9 & 10 &
 11 & 12 & 13 & 14 & 15 & 16 & 17 & 18 & 19 & 20 \rule[- 7pt]{0pt}{ 20pt} \\ \hline
 24 &  1 &  . &  . &  . &  . &  . &  . &  . &  . &  . &  . &  . &   . &  . &  . &  . &  . &  . &  . & . \rule[  0pt]{0pt}{ 13pt} \\
 2024 &  . &  1 &  . &  . &  . &  . &  . &  . &  . &  . &  . &  . &   . &  . &  . &  . &  . &  . &  . & . \\
 4576 &  . &  . &  1 &  . &  . &  . &  . &  . &  . &  . &  . &  . &   . &  . &  . &  . &  . &  . &  . & . \\
 40480 &  . &  1 &  . &  1 &  . &  . &  . &  . &  . &  . &  . &  . &   . &  . &  . &  . &  . &  . &  . & . \\
 95680 &  1 &  . &  . &  . &  1 &  . &  . &  . &  . &  . &  . &  . &   . &  . &  . &  . &  . &  . &  . & . \\
 299000 &  . &  . &  1 &  . &  . &  1 &  . &  . &  . &  . &  . &  . &   . &  . &  . &  . &  . &  . &  . & . \\
 315744 &  . &  . &  1 &  . &  . &  . &  1 &  . &  . &  . &  . &  . &   . &  . &  . &  . &  . &  . &  . & . \\
 789360 &  . &  . &  . &  . &  . &  1 &  . &  1 &  . &  . &  . &  . &   . &  . &  . &  . &  . &  . &  . & . \\
 1937520 &  1 &  . &  . &  . &  . &  . &  . &  . &  1 &  . &  . &  . &   . &  . &  . &  . &  . &  . &  . & . \\
 5051904 &  . &  . &  . &  1 &  . &  . &  . &  . &  1 &  1 &  . &  . &   . &  . &  . &  . &  . &  . &  . & . \\
 7104240 &  . &  . &  . &  . &  1 &  . &  . &  . &  . &  . &  1 &  . &   . &  . &  . &  . &  . &  . &  . & . \\
 9152000 &  . &  . &  . &  . &  . &  . &  . &  . &  . &  . &  . &  1 &   . &  . &  . &  . &  . &  . &  . & . \\
 9152000 &  . &  . &  . &  . &  . &  . &  . &  . &  . &  . &  . &  1 &  . &  . &  . &  . &  . &  . &  . & . \\
 13156000 &  . &  1 &  . &  1 &  . &  . &  . &  . &  . &  1 &  . &  . &  1 &  . &  . &  . &  . &  . &  . & . \\
 17050176 &  . &  1 &  . &  . &  . &  . &  . &  . &  . &  . &  1 &  . &  1 &  . &  . &  . &  . &  . &  . & . \\
 19734000 &  1 &  . &  . &  . &  . &  . &  . &  . &  1 &  1 &  . &  . &  . &  1 &  . &  . &  . &  . &  . & . \\
 34155000 &  . &  . &  1 &  . &  . &  1 &  1 &  . &  . &  . &  . &  . &  . &  . &  1 &  . &  . &  . &  . & . \\
 49335000 &  . &  . &  . &  . &  1 &  . &  1 &  . &  . &  . &  . &  . &  . &  . &  . &  1 &  . &  . &  . & . \\
 50519040 &  . &  . &  . &  . &  . &  . &  1 &  . &  . &  . &  . &  . &  . &  . &  . &  . &  1 &  . &  . & . \\
 67358720 &  . &  . &  . &  . &  . &  1 &  . &  1 &  . &  . &  . &  2 &  . &  1 &  1 &  . &  . &  . &  . & . \\
 210496000 &  . &  . &  . &  . &  1 &  . &  . &  . &  . &  . &  1 &  . &  1 &  . &  . &  . &  . &  1 &  . & . \\
 215547904 &  . &  . &  . &  . &  . &  . &  . &  . &  . &  . &  . &  1 &  1 &  . &  . &  . &  . &  . &  1 & . \\
 313524224 &  1 &  . &  . &  . &  . &  . &  . &  1 &  . &  1 &  . &  3 &  . &  1 &  1 &  . &  . &  . &  . & 1 \\
 394680000 &  . &  . &  . &  . &  . &  . &  1 &  . &  . &  . &  . &  3 &  . &  . &  1 &  1 &  1 &  . &  . & 1 \\
 485760000 &  1 &  . &  . &  . &  1 &  . &  . &  . &  . &  . &  . &  1 &  . &  . &  . &  1 &  . &  1 &  . & 1 \\
 517899096 &  . &  . &  . &  . &  . &  . &  . &  1 &  . &  . &  . &  4 &  . &  . &  . &  . &  1 &  . &  1 & 1 \\
 655360000 &  . &  . &  . &  . &  . &  . &  . &  . &  . &  1 &  . &  2 &  1 &  . &  . &  . &  . &  1 &  1 & 1 
\rule[- 7pt]{0pt}{  5pt} \\ \hline
\end{array}
$$
\end{table}

\begin{table}[tb]
\caption{\label{tab4} Basic set of projectives of the principal block}
$$\begin{array}{r@{\hspace{ 7pt}}*{ 9}{r@{\hspace{ 5pt}}}*{12}{r@{\hspace{ 3.5pt}}}} \hline
 \Lambda:  &  1 &  2 &  3 &  4 &  5 &  6 &  7 &  8 &  9 & 10 &
 11 & 12 & 13 & 14 & 15 & 16 & 17 & 18 & 19 & 20 & 21 \rule[- 7pt]{0pt}{ 20pt} \\ \hline
 1 &  1 &  . &  . &  . &  . &  . &  . &  . &  . &  . &  . &  . &  . &  . &  . &  . &  . &  . &  . &  . & . \rule[  0pt]{0pt}{ 13pt} \\
 276 &  1 &  1 &  . &  . &  . &  . &  . &  . &  . &  . &  . &  . &  . &  . &  . &  . &  . &  . &  . &  . & . \\
 299 &  1 &  . &  1 &  . &  . &  . &  . &  . &  . &  . &  . &  . &  . &  . &  . &  . &  . &  . &  . &  . & . \\
 17250 &  1 &  . &  1 &  1 &  . &  . &  . &  . &  . &  . &  . &  . &  . &  . &  . &  . &  . &  . &  . &  . & . \\
 80730 &  1 &  . &  . &  . &  1 &  . &  . &  . &  . &  . &  . &  . &  . &  . &  . &  . &  . &  . &  . &  . & . \\
 94875 &  1 &  1 &  . &  . &  . &  1 &  . &  . &  . &  . &  . &  . &  . &  . &  . &  . &  . &  . &  . &  . & . \\
 822250 &  1 &  . &  1 &  . &  . &  . &  1 &  . &  . &  . &  . &  . &  . &  . &  . &  . &  . &  . &  . &  . & . \\
 871884 &  1 &  . &  . &  . &  . &  . &  . &  1 &  . &  . &  . &  . &  . &  . &  . &  . &  . &  . &  . &  . & . \\
 1821600 &  1 &  1 &  . &  1 &  1 &  . &  . &  . &  1 &  . &  . &  . &  . &  . &  . &  . &  . &  . &  . &  . & . \\
 2055625 &  . &  . &  . &  1 &  . &  . &  . &  . &  . &  1 &  . &  . &  . &  . &  . &  . &  . &  . &  . &  . & . \\
 9221850 &  1 &  . &  . &  . &  . &  1 &  1 &  . &  . &  . &  1 &  . &  . &  . &  . &  . &  . &  . &  . &  . & . \\
 16347825 &  . &  . &  . &  1 &  2 &  . &  . &  . &  1 &  . &  . &  1 &  . &  . &  . &  . &  . &  . &  . &  . & . \\
 21528000 &  1 &  . &  . &  . &  . &  1 &  . &  1 &  . &  . &  . &  . &  1 &  . &  . &  . &  . &  . &  . &  . & . \\
 21579129 &  . &  . &  . &  . &  . &  . &  . &  . &  . &  . &  . &  . &  . &  1 &  1 &  . &  . &  . &  . &  . & . \\
 24667500 &  1 &  . &  . &  . &  1 &  . &  . &  . &  . &  . &  . &  1 &  . &  1 &  . &  . &  . &  . &  . &  . & . \\
 31574400 &  2 &  . &  . &  . &  . &  . &  . &  1 &  . &  . &  . &  . &  1 &  1 &  . &  . &  . &  . &  . &  . & . \\
 57544344 &  2 &  . &  1 &  1 &  . &  . &  1 &  . &  . &  . &  1 &  . &  . &  1 &  . &  1 &  . &  . &  . &  . & . \\
 66602250 &  . &  . &  . &  2 &  1 &  . &  . &  . &  1 &  1 &  . &  . &  . &  1 &  . &  . &  2 &  1 &  1 &  . & . \\
 85250880 &  2 &  . &  1 &  2 &  . &  . &  . &  . &  . &  1 &  . &  . &  . &  2 &  . &  1 &  1 &  . &  . &  . & . \\
 150732800 &  . &  . &  . &  1 &  . &  . &  . &  . &  . &  1 &  . &  . &  1 &  1 &  . &  1 &  1 &  1 &  . &  1 & 1 \\
 163478250 &  3 &  1 &  . &  1 &  . &  1 &  . &  . &  . &  . &  1 &  . &  1 &  1 &  . &  1 &  1 &  . &  . &  1 & 1 \\
 191102976 &  1 &  . &  . &  1 &  . &  . &  . &  . &  . &  . &  . &  1 &  1 &  2 &  1 &  1 &  1 &  . &  . &  1 & 1 \\
 207491625 &  1 &  . &  . &  1 &  . &  . &  . &  . &  . &  1 &  . &  . &  1 &  3 &  . &  . &  3 &  1 &  1 &  1 & . \\
 215547904 &  1 &  1 &  . &  2 &  1 &  . &  . &  . &  1 &  1 &  . &  1 &  . &  2 &  . &  1 &  3 &  1 &  1 &  1 & 1 \\
 219648000 &  . &  . &  . &  . &  . &  . &  . &  . &  . &  . &  . &  . &  2 &  1 &  1 &  . &  1 &  1 &  . &  2 & 1 \\
 299710125 &  2 &  . &  . &  1 &  . &  . &  . &  . &  . &  . &  1 &  . &  2 &  3 &  . &  1 &  2 &  1 &  . &  2 & 1 \\
 326956500 &  1 &  . &  . &  1 &  . &  . &  . &  . &  . &  1 &  1 &  . &  1 &  3 &  1 &  1 &  3 &  1 &  1 &  2 & 1 
\rule[- 7pt]{0pt}{  5pt} \\ \hline
\end{array}
$$
\end{table}
\begin{table}[tb]
\caption{\label{tab5} More projectives in the principal block}
$$\begin{array}{r@{\hspace{10pt}}*{17}{r@{\hspace{ 5pt}}}} \hline
 \Lambda: & 22 & 23 & 24 & 25 & 26 & 27 & 28 & 29 & 30 &
 31 & 32 & 33 & 34 & 35 & 36 & 37 & \rule[- 7pt]{0pt}{ 20pt} \\ \hline
 1 &  1 &  1 &  . &  . &  1 &  . &  . &  . &  . &  . &  . &  . &  . &  . &  . &  . & \rule[  0pt]{0pt}{ 13pt} \\
 276 &  . &  . &  . &  . &  . &  . &  . &  . &  . &  . &  . &  . &  . &  . &  . &  . & \\
 299 &  . &  2 &  . &  . &  1 &  . &  . &  . &  . &  . &  . &  . &  . &  . &  . &  . & \\
 17250 &  1 &  2 &  . &  . &  1 &  . &  1 &  . &  . &  . &  . &  . &  . &  . &  . &  . & \\
 80730 &  2 &  1 &  . &  . &  1 &  . &  2 &  . &  . &  . &  . &  . &  . &  . &  . &  . & \\
 94875 &  . &  2 &  . &  . &  . &  . &  . &  . &  . &  . &  . &  . &  . &  . &  . &  . & \\
 822250 &  . &  8 &  . &  . &  1 &  . &  . &  . &  . &  . &  . &  . &  . &  . &  . &  1 & \\
 871884 &  4 &  6 &  . &  . &  1 &  . &  . &  . &  . &  . &  . &  . &  . &  . &  . &  . & \\
 1821600 &  . &  1 &  . &  . &  . &  . &  1 &  . &  . &  . &  . &  . &  1 &  . &  . &  . & \\
 2055625 &  3 &  . &  . &  . &  . &  . &  1 &  . &  . &  . &  . &  . &  . &  . &  . &  . & \\
 9221850 &  . &  11 &  . &  . &  . &  . &  . &  . &  . &  . &  1 &  2 &  . &  1 &  2 &  1 & \\
 16347825 &  4 &  1 &  1 &  1 &  1 &  . &  3 &  . &  . &  1 &  . &  . &  1 &  . &  . &  . & \\
 21528000 &  5 &  13 &  . &  . &  . &  . &  . &  . &  . &  . &  1 &  2 &  . &  1 &  2 &  1 & \\
 21579129 &  1 &  . &  1 &  1 &  . &  1 &  . &  1 &  1 &  . &  . &  . &  . &  1 &  . &  . & \\
 24667500 &  6 &  1 &  2 &  2 &  2 &  1 &  2 &  1 &  . &  1 &  . &  . &  . &  . &  . &  . & \\
 31574400 &  7 &  12 &  1 &  1 &  1 &  1 &  . &  1 &  . &  . &  1 &  2 &  . &  1 &  2 &  1 & \\
 57544344 &  5 &  19 &  . &  1 &  2 &  1 &  . &  1 &  . &  . &  3 &  3 &  . &  2 &  2 &  2 & \\
 66602250 &  10 &  1 &  2 &  3 &  2 &  2 &  2 &  . &  . &  1 &  . &  . &  1 &  . &  . &  . & \\
 85250880 &  14 &  10 &  1 &  3 &  4 &  1 &  1 &  3 &  . &  1 &  3 &  1 &  . &  2 &  1 &  1 & \\
 150732800 &  13 &  12 &  3 &  2 &  1 &  2 &  . &  1 &  . &  . &  2 &  1 &  . &  2 &  . &  1 & \\
 163478250 &  5 &  18 &  3 &  1 &  1 &  1 &  . &  2 &  . &  . &  4 &  4 &  1 &  3 &  4 &  1 & \\
 191102976 &  9 &  13 &  5 &  3 &  2 &  2 &  . &  3 &  1 &  1 &  3 &  2 &  1 &  3 &  2 &  1 & \\
 207491625 &  20 &  1 &  7 &  5 &  4 &  3 &  1 &  4 &  2 &  2 &  1 &  . &  1 &  2 &  1 &  . & \\
 215547904 &  15 &  8 &  6 &  5 &  4 &  3 &  1 &  2 &  . &  2 &  2 &  . &  2 &  1 &  . &  . & \\
 219648000 &  10 &  5 &  7 &  1 &  . &  2 &  . &  2 &  3 &  . &  . &  . &  1 &  3 &  . &  . & \\
 299710125 &  21 &  16 &  8 &  4 &  3 &  3 &  . &  5 &  2 &  1 &  4 &  3 &  1 &  5 &  3 &  1 & \\
 326956500 &  16 &  11 &  9 &  4 &  3 &  4 &  . &  4 &  3 &  1 &  3 &  2 &  2 &  4 &  2 &  . & 
\rule[- 7pt]{0pt}{  5pt} \\ \hline
\end{array}
$$
\end{table}
\begin{table}[tb]
\caption{\label{expr2} Relations in the principal block}
$$\begin{array}{r@{\hspace{ 6pt}}*{21}{r@{\hspace{ 3.5pt}}}} \hline
 \Lambda:  &  1 &  2 &  3 &  4 &  5 &  6 &  7 &  8 &  9 & 10 & 11
 & 12 & 13 & 14 & 15 & 16 & 17 & 18 & 19 & 20 & 21
 \rule[- 7pt]{0pt}{ 20pt} \\ \hline
\Lambda_{22} &  1 & -1 & -1 &  1 &  1 &  . &  . &  3 & -2 &  2 & -1
  &  3 &  1 &  1 &  . &  3 &  4 &  7 & -9 &  2 &-8 \rule[  0pt]{0pt}{ 13pt} \\
\Lambda_{23} &  1 & -1 &  1 &  . &  . &  2 &  6 &  5 &  1 &  . &  2
  &  . &  5 &  . &  . &  8 & -1 &  . &  2 & -4 & 4 \\
\Lambda_{24} &  . &  . &  . &  . &  . &  . &  . &  . &  . &  . &  .
  &  1 &  . &  1 &  . & -1 &  . &  . &  1 &  3 & . \\
\Lambda_{25} &  . &  . &  . &  . &  . &  . &  . &  . &  . &  . &  .
  &  1 &  . &  1 &  . &  . &  1 &  1 & -1 & -1 & . \\
\Lambda_{26} &  1 & -1 &  . &  . &  . &  . &  . &  . &  . &  . & -1
  &  1 & -1 &  . &  . &  1 &  1 &  1 & -1 &  1 &-2 \\
\Lambda_{27} &  . &  . &  . &  . &  . &  . &  . &  . &  . &  . &  .
  &  . &  . &  1 &  . &  . & -1 &  1 &  2 &  . & 1 \\
\Lambda_{28} &  . &  . &  . &  1 &  2 &  . &  . &  . & -2 &  . &  .
  &  . &  . &  . &  . & -1 &  . &  . &  . &  . & . \\
\Lambda_{29} &  . &  . &  . &  . &  . &  . &  . &  . &  . &  . &  .
  &  . &  . &  1 &  . &  . &  1 & -1 & -2 &  1 &-1 \\
\Lambda_{30} &  . &  . &  . &  . &  . &  . &  . &  . &  . &  . &  .
  &  . &  . &  . &  1 &  . &  . &  . &  . &  2 &-2 \\
\Lambda_{31} &  . &  . &  . &  . &  . &  . &  . &  . &  . &  . &  .
  &  1 &  . &  . &  . &  . &  1 &  . & -1 &  . &-1 \\
\Lambda_{32} &  . &  . &  . &  . &  . &  . &  . &  . &  . &  . &  1
  &  . &  1 &  . &  . &  2 &  1 & -1 & -1 & -1 & . \\
\Lambda_{33} &  . &  . &  . &  . &  . &  . &  . &  . &  . &  . &  2
  &  . &  2 &  . &  . &  1 &  . & -1 &  1 & -2 & 1 \\
\Lambda_{34} &  . &  . &  . &  . &  . &  . &  . &  . &  1 &  . &  .
  &  . &  . &  . &  . &  . &  . & -1 &  1 &  1 & . \\
\Lambda_{35} &  . &  . &  . &  . &  . &  . &  . &  . &  . &  . &  1
  &  . &  1 &  . &  1 &  1 &  1 &  . & -2 &  . &-1 \\
\Lambda_{36} &  . &  . &  . &  . &  . &  . &  . &  . &  . &  . &  2
  &  . &  2 &  . &  . &  . &  1 & -2 &  . & -2 & 1 \\
\Lambda_{37} &  . &  . &  . &  . &  . &  . &  1 &  . &  . &  . &  .
  &  . &  1 &  . &  . &  1 &  . &  . &  . & -1 & .
\rule[- 7pt]{0pt}{  5pt} \\ \hline
\end{array}
$$
\end{table}
\begin{table}[b]
\caption{\label{orig2} Origin of projective characters of the principal block}
$$\begin{array}[t]{l@{\hspace{10pt}}cc} \hline
\mbox{\rm Char.} & \multicolumn{1}{c}{\mbox{\rm Origin}}
\rule[- 7pt]{0pt}{ 20pt} \\ \hline
   \Lambda   _{ 1 }: &    \chi   _{102} \otimes      \Phi   _{1,7} & \rule[  0pt]{0pt}{ 13pt} \\
   \Lambda   _{ 2 }: & {\rm Ind}_{\rm 2xCo2} (   \theta_{3} + \theta_{16}) & \\
   \Lambda   _{ 3 }: &    \chi   _{102} \otimes      \Phi   _{1,8} & \\
   \Lambda   _{ 4 }: &    \chi   _{102} \otimes      \Phi   _{7,7} & \\
   \Lambda   _{ 5 }: &    \chi   _{2} \otimes        \chi   _{39} & \\
   \Lambda   _{ 6 }: &    \chi   _{102} \otimes      \Phi   _{4,7} & \\
   \Lambda   _{ 7 }: &    \chi   _{102} \otimes      \chi   _{119} & \\
   \Lambda   _{ 8 }: &    \chi   _{2} \otimes        \chi   _{22} & \\
   \Lambda   _{ 9 }: &    \chi   _{2} \otimes        \chi   _{32} & \\
   \Lambda   _{10 }: &    \chi   _{6} \otimes        \chi   _{17} & \\
   \Lambda   _{11 }: &    \chi   _{102} \otimes      \chi   _{141} & \\
   \Lambda   _{12 }: &    \chi   _{4} \otimes        \chi   _{33} & \\
   \Lambda   _{13 }: &    \chi   _{2} \otimes        \chi   _{44} & \\
   \Lambda   _{14 }: & {\rm Ind}_{\rm 2xCo2} (   \theta   _{43}) & \\
   \Lambda   _{15 }: &    \chi   _{4} \otimes        \chi   _{26} & \\
   \Lambda   _{16 }: &    \chi   _{102} \otimes      \chi   _{149} & \\
   \Lambda   _{17 }: &    \chi   _{3} \otimes        \chi   _{73} & \\
   \Lambda   _{18 }: &    \chi   _{4} \otimes        \chi   _{32} & \\
   \Lambda   _{19 }: &    (\chi   _{2} \otimes        \chi   _{70})/2 & 
\rule[- 7pt]{0pt}{  5pt} \\ \hline
\end{array}
\quad\quad\quad
\begin{array}[t]{l@{\hspace{10pt}}cc} \hline
\mbox{\rm Char.} & \multicolumn{1}{c}{\mbox{\rm Origin}}
\rule[- 7pt]{0pt}{ 20pt} \\ \hline
   \Lambda   _{20 }: &    \chi   _{102} \otimes      \chi   _{161} & \rule[  0pt]{0pt}{ 13pt} \\
   \Lambda   _{21 }: &    \chi   _{2} \otimes        \chi   _{46} & \\
   \Lambda   _{22 }: &    \chi   _{16}^{2+} &  \\
   \Lambda   _{23 }: &    \chi   _{112}^{2+} & \\
   \Lambda   _{24 }: &    \chi   _{7} \otimes        \chi   _{33} & \\
   \Lambda   _{25 }: & {\Ind}_{\rm 2xCo2} (   \theta_{45}) & \\
   \Lambda   _{26 }: & {\Ind}_{\rm 2xCo2} (   \theta_{ 1} + \theta_{44}) & \\
   \Lambda   _{27 }: & {\Ind}_{\rm 2xCo2} (   \theta_{27} + \theta_{34}) & \\
   \Lambda   _{28 }: & {\Ind}_{\rm 2xCo2} (   \theta_{ 6} + \theta_{17}) & \\
   \Lambda   _{29 }: &    \chi   _{102} \otimes      \Phi   _{21,7} & \\
   \Lambda   _{30 }: &    \chi   _{102} \otimes      \Phi   _{20,7} & \\
   \Lambda   _{31 }: &    \chi   _{102} \otimes      \Phi   _{19,7} & \\
   \Lambda   _{32 }: &    \chi   _{102} \otimes      \Phi   _{16,7} & \\
   \Lambda   _{33 }: &    \chi   _{102} \otimes      \Phi   _{15,7} & \\
   \Lambda   _{34 }: &    \chi   _{102} \otimes      \Phi   _{14,7} & \\
   \Lambda   _{35 }: &    \chi   _{102} \otimes      \Phi   _{13,7} & \\
   \Lambda   _{36 }: &    \chi   _{102} \otimes      \Phi   _{10,7} & \\
   \Lambda   _{37 }: &    \chi   _{102} \otimes      \Phi   _{3,8} &
\rule[- 7pt]{0pt}{  5pt} \\ \hline
\end{array}
$$
\end{table}
\begin{table}[tb]
\caption{\label{tab6} Refined basic set of projectives of the principal block}
$$\begin{array}{r@{\hspace{ 7pt}}*{10}{r@{\hspace{ 5pt}}}*{11}{r@{\hspace{ 3.5pt}}}} \hline
 \Phi:   &  1 &  2 &  3 &  4 &  5 &  6 &  7 &  8 &  9 & 10 &
 11 & 12 & 13 & 14 & 15 & 16 & 17 & 18 & 19 & 20 & 21 \rule[- 7pt]{0pt}{ 20pt} \\ \hline
 1 &  1 &  . &  . &  . &  . &  . &  . &  . &  . &  . &  . &  . &  . &  . &  . &  . &  . &  . &  . &  . & . \rule[  0pt]{0pt}{ 13pt} \\
 276 &  . &  1 &  . &  . &  . &  . &  . &  . &  . &  . &  . &  . &  . &  . &  . &  . &  . &  . &  . &  . & . \\
 299 &  . &  . &  1 &  . &  . &  . &  . &  . &  . &  . &  . &  . &  . &  . &  . &  . &  . &  . &  . &  . & . \\
 17250 &  . &  . &  1 &  1 &  . &  . &  . &  . &  . &  . &  . &  . &  . &  . &  . &  . &  . &  . &  . &  . & . \\
 80730 &  1 &  . &  . &  . &  1 &  . &  . &  . &  . &  . &  . &  . &  . &  . &  . &  . &  . &  . &  . &  . & . \\
 94875 &  . &  1 &  . &  . &  . &  1 &  . &  . &  . &  . &  . &  . &  . &  . &  . &  . &  . &  . &  . &  . & . \\
 822250 &  . &  . &  1 &  . &  . &  . &  1 &  . &  . &  . &  . &  . &  . &  . &  . &  . &  . &  . &  . &  . & . \\
 871884 &  1 &  . &  . &  . &  . &  . &  . &  1 &  . &  . &  . &  . &  . &  . &  . &  . &  . &  . &  . &  . & . \\
 1821600 &  . &  1 &  . &  . &  . &  . &  . &  . &  1 &  . &  . &  . &  . &  . &  . &  . &  . &  . &  . &  . & . \\
 2055625 &  . &  . &  . &  1 &  . &  . &  . &  . &  . &  1 &  . &  . &  . &  . &  . &  . &  . &  . &  . &  . & . \\
 9221850 &  . &  . &  . &  . &  . &  1 &  1 &  . &  . &  . &  1 &  . &  . &  . &  . &  . &  . &  . &  . &  . & . \\
 16347825 &  . &  . &  . &  . &  1 &  . &  . &  . &  1 &  . &  . &  1 &  . &  . &  . &  . &  . &  . &  . &  . & . \\
 21528000 &  . &  . &  . &  . &  . &  1 &  . &  1 &  . &  . &  . &  . &  1 &  . &  . &  . &  . &  . &  . &  . & . \\
 21579129 &  . &  . &  . &  . &  . &  . &  . &  . &  . &  . &  . &  . &  . &  1 &  1 &  . &  . &  . &  . &  . & . \\
 24667500 &  1 &  . &  . &  . &  1 &  . &  . &  . &  . &  . &  . &  1 &  . &  1 &  . &  . &  . &  . &  . &  . & . \\
 31574400 &  1 &  . &  . &  . &  . &  . &  . &  1 &  . &  . &  . &  . &  1 &  1 &  . &  . &  . &  . &  . &  . & . \\
 57544344 &  . &  . &  1 &  . &  . &  . &  1 &  . &  . &  . &  1 &  . &  . &  . &  . &  1 &  . &  . &  . &  . & . \\
 66602250 &  . &  . &  . &  1 &  . &  . &  . &  . &  1 &  1 &  . &  . &  . &  . &  . &  . &  1 &  . &  . &  . & . \\
 85250880 &  . &  . &  1 &  1 &  . &  . &  . &  . &  . &  1 &  . &  . &  . &  . &  . &  1 &  . &  1 &  . &  . & . \\
 150732800 &  . &  . &  . &  . &  . &  . &  . &  . &  . &  1 &  . &  . &  . &  . &  . &  1 &  . &  . &  1 &  . & . \\
 163478250 &  . &  1 &  . &  . &  . &  1 &  . &  . &  . &  . &  1 &  . &  1 &  1 &  . &  . &  . &  . &  . &  1 & . \\
 191102976 &  . &  . &  . &  . &  . &  . &  . &  . &  . &  . &  . &  1 &  1 &  2 &  1 &  . &  . &  . &  . &  1 & . \\
 207491625 &  . &  . &  . &  1 &  . &  . &  . &  . &  . &  1 &  . &  . &  . &  . &  . &  . &  1 &  1 &  . &  . & 1 \\
 215547904 &  . &  1 &  . &  . &  . &  . &  . &  . &  1 &  1 &  . &  1 &  . &  1 &  . &  . &  1 &  . &  . &  1 & . \\
 219648000 &  . &  . &  . &  . &  . &  . &  . &  . &  . &  . &  . &  . &  . &  . &  1 &  . &  . &  . &  1 &  . & 1 \\
 299710125 &  . &  . &  . &  . &  . &  . &  . &  . &  . &  . &  1 &  . &  . &  . &  . &  1 &  . &  1 &  1 &  . & 1 \\
 326956500 &  . &  . &  . &  . &  . &  . &  . &  . &  . &  1 &  1 &  . &  . &  1 &  1 &  . &  1 &  . &  . &  1 & 1 
\rule[- 7pt]{0pt}{  5pt} \\ \hline
\end{array}
$$
\end{table}

\clearpage
\thispagestyle{empty}

%
%
%

\appendix

\clearpage
\thispagestyle{empty}

\chapter{}
\renewcommand{\arraystretch}{0.90}
\vspace*{-1cm}
\section{The $5$-decomposition numbers of $\mbox{\em Co}_2$} 
{\small
\[
\begin{array}{rrrrrrrrrrrrrrrrrr}
 1 &  1 &  . &  . &  . &  . &  . &  . &  . &  . &  . &  . &  . &  . &  . &  . & . \\
 23 &  . &  1 &  . &  . &  . &  . &  . &  . &  . &  . &  . &  . &  . &  . &  . & . \\
 253 &  . &  . &  1 &  . &  . &  . &  . &  . &  . &  . &  . &  . &  . &  . &  . & . \\
 1771 &  . &  . &  . &  1 &  . &  . &  . &  . &  . &  . &  . &  . &  . &  . &  . & . \\
 2024 &  1 &  . &  . &  . &  1 &  . &  . &  . &  . &  . &  . &  . &  . &  . &  . & . \\
 2277 &  . &  1 &  . &  . &  . &  1 &  . &  . &  . &  . &  . &  . &  . &  . &  . & . \\
 7084 &  . &  . &  . &  1 &  . &  . &  1 &  . &  . &  . &  . &  . &  . &  . &  . & . \\
 10395 &  . &  . &  . &  . &  . &  . &  . &  1 &  . &  . &  . &  . &  . &  . &  . & . \\
 10395 &  . &  . &  . &  . &  . &  . &  . &  1 &  . &  . &  . &  . &  . &  . &  . & . \\
 31878 &  . &  . &  . &  . &  . &  1 &  . &  . &  1 &  . &  . &  . &  . &  . &  . & . \\
 37422 &  . &  . &  1 &  . &  . &  . &  . &  . &  . &  1 &  . &  . &  . &  . &  . & . \\
 129536 &  . &  . &  1 &  1 &  . &  . &  1 &  . &  . &  . &  1 &  . &  . &  . &  . & . \\
 184437 &  1 &  . &  . &  . &  1 &  . &  1 &  . &  . &  . &  . &  1 &  . &  . &  . & . \\
 212520 &  . &  . &  . &  . &  . &  2 &  . &  . &  1 &  . &  . &  . &  1 &  . &  . & . \\
 226688 &  1 &  . &  . &  . &  1 &  . &  . &  1 &  . &  1 &  . &  1 &  . &  . &  . & . \\
 239085 &  . &  . &  . &  . &  . &  . &  . &  . &  . &  . &  . &  . &  . &  1 &  . & . \\
 239085 &  . &  . &  . &  . &  . &  . &  . &  . &  . &  . &  . &  . &  . &  1 &  . & . \\
 245916 &  . &  . &  . &  . &  1 &  . &  . &  . &  1 &  1 &  . &  1 &  . &  . &  . & . \\
 312984 &  . &  . &  . &  1 &  . &  . &  2 &  . &  . &  . &  1 &  . &  1 &  . &  . & . \\
 368874 &  . &  1 &  1 &  . &  . &  . &  . &  1 &  . &  . &  1 &  . &  . &  . &  1 & . \\
 430353 &  . &  . &  . &  . &  . &  1 &  2 &  . &  . &  . &  . &  . &  1 &  1 &  . & . \\
 637560 &  1 &  . &  . &  . &  1 &  . &  1 &  1 &  1 &  . &  . &  2 &  . &  . &  1 & . \\
 1291059 &  . &  . &  . &  . &  . &  1 &  . &  . &  1 &  1 &  . &  1 &  . &  . &  . & 1 \\
 1835008 &  . &  . &  . &  . &  . &  . &  1 &  1 &  . &  . &  1 &  1 &  . &  1 &  1 & 1 \\
 1943040 &  1 &  . &  1 &  . &  . &  . &  1 &  2 &  . &  1 &  2 &  2 &  . &  . &  1 & 1 \\
 2040192 &  . &  . &  . &  . &  . &  1 &  2 &  . &  1 &  . &  1 &  1 &  1 &  1 &  1 & 1 \\
 2072576 &  . &  1 &  . &  . &  . &  1 &  2 &  . &  . &  . &  1 &  . &  1 &  2 &  1 & 1 
\end{array}
\]
}

\pagebreak

\renewcommand{\arraystretch}{1.00}

\section{The $7$-decomposition numbers of $2\mbox{\em Co}_1$} 
\label{Co1mod7}
In this section we give the $7$-modular decomposition
matrices for the two blocks of maximal defect of
$2\mbox{\em Co}_1$\@. 
The decomposition 
numbers are followed by tables giving the degrees of the irreducible
Brauer characters as well as the prime factorizations
of these degrees.
The decomposition matrices for the blocks of non-maximal defect
can be found in~\cite{BBB}\@.

\subsection{Decomposition matrix of Block~1}
{\small
\noindent
$$\begin{array}{r@{\hspace{ 7pt}}*{10}{r@{\hspace{ 5pt}}}*{11}{r@{\hspace{ 3.5pt}}}} \hline
 \Phi:   &  1 &  2 &  3 &  4 &  5 &  6 &  7 &  8 &  9 & 10 &
 11 & 12 & 13 & 14 & 15 & 16 & 17 & 18 & 19 & 20 & 21 \rule[- 7pt]{0pt}{ 20pt} \\ \hline
 1 &  1 &  . &  . &  . &  . &  . &  . &  . &  . &  . &  . &  . &  . &  . &  . &  . &  . &  . &  . &  . & . \rule[  0pt]{0pt}{ 13pt} \\
 276 &  . &  1 &  . &  . &  . &  . &  . &  . &  . &  . &  . &  . &  . &  . &  . &  . &  . &  . &  . &  . & . \\
 299 &  . &  . &  1 &  . &  . &  . &  . &  . &  . &  . &  . &  . &  . &  . &  . &  . &  . &  . &  . &  . & . \\
 17250 &  . &  . &  1 &  1 &  . &  . &  . &  . &  . &  . &  . &  . &  . &  . &  . &  . &  . &  . &  . &  . & . \\
 80730 &  1 &  . &  . &  . &  1 &  . &  . &  . &  . &  . &  . &  . &  . &  . &  . &  . &  . &  . &  . &  . & . \\
 94875 &  . &  1 &  . &  . &  . &  1 &  . &  . &  . &  . &  . &  . &  . &  . &  . &  . &  . &  . &  . &  . & . \\
 822250 &  . &  . &  1 &  . &  . &  . &  1 &  . &  . &  . &  . &  . &  . &  . &  . &  . &  . &  . &  . &  . & . \\
 871884 &  1 &  . &  . &  . &  . &  . &  . &  1 &  . &  . &  . &  . &  . &  . &  . &  . &  . &  . &  . &  . & . \\
 1821600 &  . &  1 &  . &  . &  . &  . &  . &  . &  1 &  . &  . &  . &  . &  . &  . &  . &  . &  . &  . &  . & . \\
 2055625 &  . &  . &  . &  1 &  . &  . &  . &  . &  . &  1 &  . &  . &  . &  . &  . &  . &  . &  . &  . &  . & . \\
 9221850 &  . &  . &  . &  . &  . &  1 &  1 &  . &  . &  . &  1 &  . &  . &  . &  . &  . &  . &  . &  . &  . & . \\
 16347825 &  . &  . &  . &  . &  1 &  . &  . &  . &  1 &  . &  . &  1 &  . &  . &  . &  . &  . &  . &  . &  . & . \\
 21528000 &  . &  . &  . &  . &  . &  1 &  . &  1 &  . &  . &  . &  . &  1 &  . &  . &  . &  . &  . &  . &  . & . \\
 21579129 &  . &  . &  . &  . &  . &  . &  . &  . &  . &  . &  . &  . &  . &  1 &  1 &  . &  . &  . &  . &  . & . \\
 24667500 &  1 &  . &  . &  . &  1 &  . &  . &  . &  . &  . &  . &  1 &  . &  1 &  . &  . &  . &  . &  . &  . & . \\
 31574400 &  1 &  . &  . &  . &  . &  . &  . &  1 &  . &  . &  . &  . &  1 &  1 &  . &  . &  . &  . &  . &  . & . \\
 57544344 &  . &  . &  1 &  . &  . &  . &  1 &  . &  . &  . &  1 &  . &  . &  . &  . &  1 &  . &  . &  . &  . & . \\
 66602250 &  . &  . &  . &  1 &  . &  . &  . &  . &  1 &  . &  . &  . &  . &  . &  . &  . &  1 &  . &  . &  . & . \\
 85250880 &  . &  . &  1 &  1 &  . &  . &  . &  . &  . &  1 &  . &  . &  . &  . &  . &  1 &  . &  1 &  . &  . & . \\
 150732800 &  . &  . &  . &  . &  . &  . &  . &  . &  . &  1 &  . &  . &  . &  . &  . &  1 &  . &  . &  1 &  . & . \\
 163478250 &  . &  1 &  . &  . &  . &  1 &  . &  . &  . &  . &  1 &  . &  1 &  . &  . &  . &  . &  . &  . &  1 & . \\
 191102976 &  . &  . &  . &  . &  . &  . &  . &  . &  . &  . &  . &  1 &  1 &  1 &  1 &  . &  . &  . &  . &  1 & . \\
 207491625 &  . &  . &  . &  1 &  . &  . &  . &  . &  . &  . &  . &  . &  . &  . &  . &  . &  1 &  1 &  . &  . & 1 \\
 215547904 &  . &  1 &  . &  . &  . &  . &  . &  . &  1 &  . &  . &  1 &  . &  . &  . &  . &  1 &  . &  . &  1 & . \\
 219648000 &  . &  . &  . &  . &  . &  . &  . &  . &  . &  . &  . &  . &  . &  . &  1 &  . &  . &  . &  1 &  . & 1 \\
 299710125 &  . &  . &  . &  . &  . &  . &  . &  . &  . &  . &  1 &  . &  . &  . &  . &  1 &  . &  1 &  1 &  . & 1 \\
 326956500 &  . &  . &  . &  . &  . &  . &  . &  . &  . &  . &  1 &  . &  . &  . &  1 &  . &  1 &  . &  . &  1 & 1 
\rule[- 7pt]{0pt}{  5pt} \\ \hline
\end{array}$$
}

\subsection{Decomposition matrix of Block~7}

$$\begin{array}{r@{\hspace{ 7pt}}*{10}{r@{\hspace{ 5pt}}}*{11}{r@{\hspace{ 3.5pt}}}} \hline
  \Phi:  &  1 &  2 &  3 &  4 &  5 &  6 &  7 &  8 &  9 & 10 &
 11 & 12 & 13 & 14 & 15 & 16 & 17 & 18 & 19 & 20 & 21 \rule[- 7pt]{0pt}{ 20pt} \\ \hline
 24 &  1 &  . &  . &  . &  . &  . &  . &  . &  . &  . &  . &  . &  . &  . &  . &  . &  . &  . &  . &  . & . \rule[  0pt]{0pt}{ 13pt} \\

 2024 &  . &  1 &  . &  . &  . &  . &  . &  . &  . &  . &  . &  . &  . &  . &  . &  . &  . &  . &  . &  . & . \\
 4576 &  . &  . &  1 &  . &  . &  . &  . &  . &  . &  . &  . &  . &  . &  . &  . &  . &  . &  . &  . &  . & . \\
 40480 &  . &  1 &  . &  1 &  . &  . &  . &  . &  . &  . &  . &  . &  . &  . &  . &  . &  . &  . &  . &  . & . \\
 95680 &  1 &  . &  . &  . &  1 &  . &  . &  . &  . &  . &  . &  . &  . &  . &  . &  . &  . &  . &  . &  . & . \\
 299000 &  . &  . &  1 &  . &  . &  1 &  . &  . &  . &  . &  . &  . &  . &  . &  . &  . &  . &  . &  . &  . & . \\
 315744 &  . &  . &  1 &  . &  . &  . &  1 &  . &  . &  . &  . &  . &  . &  . &  . &  . &  . &  . &  . &  . & . \\
 789360 &  . &  . &  . &  . &  . &  1 &  . &  1 &  . &  . &  . &  . &  . &  . &  . &  . &  . &  . &  . &  . & . \\
 1937520 &  1 &  . &  . &  . &  . &  . &  . &  . &  1 &  . &  . &  . &  . &  . &  . &  . &  . &  . &  . &  . & . \\
 5051904 &  . &  . &  . &  1 &  . &  . &  . &  . &  1 &  1 &  . &  . &  . &  . &  . &  . &  . &  . &  . &  . & . \\
 7104240 &  . &  . &  . &  . &  1 &  . &  . &  . &  . &  . &  1 &  . &  . &  . &  . &  . &  . &  . &  . &  . & . \\
 9152000 &  . &  . &  . &  . &  . &  . &  . &  . &  . &  . &  . &  1 &  . &  . &  . &  . &  . &  . &  . &  . & . \\
 9152000 &  . &  . &  . &  . &  . &  . &  . &  . &  . &  . &  . &  . &  1 &  . &  . &  . &  . &  . &  . &  . & . \\
 13156000 &  . &  1 &  . &  1 &  . &  . &  . &  . &  . &  1 &  . &  . &  . &  1 &  . &  . &  . &  . &  . &  . & . \\
 17050176 &  . &  1 &  . &  . &  . &  . &  . &  . &  . &  . &  1 &  . &  . &  1 &  . &  . &  . &  . &  . &  . & . \\
 19734000 &  1 &  . &  . &  . &  . &  . &  . &  . &  1 &  1 &  . &  . &  . &  . &  1 &  . &  . &  . &  . &  . & . \\
 34155000 &  . &  . &  1 &  . &  . &  1 &  1 &  . &  . &  . &  . &  . &  . &  . &  . &  1 &  . &  . &  . &  . & . \\
 49335000 &  . &  . &  . &  . &  1 &  . &  1 &  . &  . &  . &  . &  . &  . &  . &  . &  . &  1 &  . &  . &  . & . \\
 50519040 &  . &  . &  . &  . &  . &  . &  1 &  . &  . &  . &  . &  . &  . &  . &  . &  . &  . &  1 &  . &  . & . \\
 67358720 &  . &  . &  . &  . &  . &  1 &  . &  1 &  . &  . &  . &  1 &  1 &  . &  1 &  1 &  . &  . &  . &  . & . \\
 210496000 &  . &  . &  . &  . &  1 &  . &  . &  . &  . &  . &  1 &  . &  . &  1 &  . &  . &  . &  . &  1 &  . & . \\
 215547904 &  . &  . &  . &  . &  . &  . &  . &  . &  . &  . &  . &  . &  . &  1 &  . &  . &  . &  . &  . &  1 & . \\
 313524224 &  1 &  . &  . &  . &  . &  . &  . &  1 &  . &  1 &  . &  1 &  1 &  . &  1 &  1 &  . &  . &  . &  . & 1 \\
 394680000 &  . &  . &  . &  . &  . &  . &  1 &  . &  . &  . &  . &  1 &  1 &  . &  . &  1 &  1 &  1 &  . &  . & 1 \\
 485760000 &  1 &  . &  . &  . &  1 &  . &  . &  . &  . &  . &  . &  . &  . &  . &  . &  . &  1 &  . &  1 &  . & 1 \\
 517899096 &  . &  . &  . &  . &  . &  . &  . &  1 &  . &  . &  . &  1 &  1 &  . &  . &  . &  . &  1 &  . &  1 & 1 \\
 655360000 &  . &  . &  . &  . &  . &  . &  . &  . &  . &  1 &  . &  . &  . &  1 &  . &  . &  . &  . &  1 &  1 & 1 
\rule[- 7pt]{0pt}{  5pt} \\ \hline
\end{array}
$$

\clearpage

\subsection{Degrees of irreducible Brauer characters}

The table on the left hand side lists the degrees of the
irreducible Brauer characters of the principal block, the
other table those of Block~$7$\@.

\bigskip

\noindent $\begin{array}{lrc} \hline
\mbox{\rm Char.} & \multicolumn{1}{c}{\mbox{\rm Degree}}
                 & \multicolumn{1}{c}{\mbox{\rm Factors}}
\rule[- 7pt]{0pt}{ 20pt} \\ \hline
\phi_{ 1} & 1  & 1  \rule[  0pt]{0pt}{ 13pt}  \\
\phi_{ 2} & 276  & 2^2 \, 3 \,23 \\
\phi_{ 3} & 299  & 13\, 23 \\
\phi_{ 4} & 16951  & 11\, 23\, 67 \\
\phi_{ 5} & 80729  & 11\, 41\, 179 \\
\phi_{ 6} & 94599  & 3^2 \, 23\, 457 \\
\phi_{ 7} & 821951 & 13\, 23\, 2749 \\
\phi_{ 8}  & 871883  & 871883 \\
\phi_{ 9}  & 1821324  & 2^2\,  3\, 23\, 6599 \\
\phi_{10}  & 2038674  & 2\, 3\, 11\, 17\, 23\, 79 \\
\phi_{11}  & 8305300  & 2^2\,  5^2\,  23^2\,  157 \\
\phi_{12}  & 14445772  & 2^2\,  11\, 569\, 577 \\
\phi_{13}  & 20561518  & 2\, 607\, 16937 \\
\phi_{14}  & 10140998  & 2\, 7\, 227\, 3191 \\
\phi_{15}  & 11438131  & 11438131 \\
\phi_{16}  & 48416794  & 2\, 23\, 37\, 28447 \\
\phi_{17}  & 64763975  & 5^2\,  23\, 163\, 691 \\
\phi_{18}  & 34778162  & 2\, 23\, 83\, 9109 \\
\phi_{19}  & 100277332  & 2^2\,  23\, 733\, 1487 \\
\phi_{20}  & 134516557  & 7\, 19216651 \\
\phi_{21}  & 107932537  & 23\, 43\, 109133
\rule[- 7pt]{0pt}{  5pt} \\ \hline
\end{array}$
\hfill
$
\begin{array}{rc} \hline
 \multicolumn{1}{c}{\mbox{\rm Degree}} 
                 & \multicolumn{1}{c}{\mbox{\rm Factors}}
\rule[- 7pt]{0pt}{ 20pt} \\ \hline
 24 & 2^3\,  3  \rule[  0pt]{0pt}{ 13pt} \\
 2024 & 2^3\,  11\, 23 \\
 4576 & 2^5\,  11\, 13 \\
 38456 & 2^3  11\, 19\, 23 \\
 95656 & 2^3\,  11\, 1087 \\
 294424 & 2^3\,  13\, 19\, 149 \\
 311168 & 2^7\,  11\, 13\, 17 \\
 494936 & 2^3\,  13\, 4759 \\
 1937496 & 2^3\,  3\, 11\, 41\, 179 \\
 3075952 & 2^4\,  11\, 17477 \\
 7008584 & 2^3\,  11\, 73\, 1091 \\
 9152000 & 2^9\,  5^3\,  11\, 13 \\
 9152000 & 2^9\,  5^3\,  11\, 13 \\
 10039568 & 2^4\,  7\, 11\, 29\, 281 \\
 14720528 & 2^4\,  283\, 3251 \\
 33544832 & 2^7\,  262069 \\
 48928176 & 2^4\,  3^2\,  11\, 17\, 23\, 79 \\
 50207872 & 2^7\,  11\, 13^2\,  211 \\
 193352192 & 2^9\,  11^2\,  3121 \\
 205508336 & 2^4\,  11^2\,  101\, 1051 \\
 243383952 & 2^4\,  3\, 7\, 227\, 3191 
\rule[- 7pt]{0pt}{  5pt} \\ \hline
\end{array}$

\pagebreak
\clearpage
 
\chapter{}
\vspace*{-1.3cm}
\section{The $5$-modular character table of $\mbox{\em Co}_2$} 
\begin{verbatim}
   ;        @      @      @      @    @    @      @    @    @ 
           42    743     41      1   46   15      3   12    7 
 305421312000 178240 287680 474560 6560 5520 096576 2880 3728 
      p power      A      A      A    A    A      A    B    B 
      p' part      A      A      A    A    A      A    A    A 
 ind       1A     2A     2B     2C   3A   3B     4A   4B   4C 
   +        1      1      1      1    1    1      1    1    1 
   +       23     -9      7     -1   -4    5      7   -5    3 
   +      253     29     13    -11   10   10     29    9    1 
   +      275     51     35     11    5   14     19   15    7 
   +     1771    -21    -21     11  -11   16     91   -5   -5 
   +     2023    231    103     39   -2   25      7   23   23 
   +     2254   -210    126    -10   13   31     14  -30   10 
   +     4025   -231    105      1  -25   29    105  -35    5 
   +     5313    -63    -63     33   21    3     49    1    1 
   o     9625   -455    105    -15   40   -5     -7    5   29 
   o     9625   -455    105    -15   40   -5     -7    5   29 
   +    10395    315    -21    -45   27    0    -21   -9   15 
   +    12650    554    330     26  -40   59     -6   50    2 
   +    23000    600    280    120   50    5    184   40    8 
   +    29624   -168    392    -16   32   14    -56    4  -36 
   +    31625    265    -55    -55   35   35    377   25   -7 
   +    31625   1385    505    145   35   35     41   45   53 
   +    37169   1105    449    -55  -10   71     97   77    5 
   +    44275  -1869    595    -29   -5   94     35  -85   59 
   +    63250   -110    210   -110  -65  -20    322  -30    2 
\end{verbatim}

\pagebreak

\noindent\begin{verbatim}
   o    91125    405   -315     45    0    0    -27  -15    9 
   o    91125    405   -315     45    0    0    -27  -15    9 
   +   113575    903   -105    143   40  -23    119   -5    3 
   +   122199    567   -441    -33  -84   15    343   -5    3 
   +   173075   -973     35     43    5  -58    483   15    7 
   +   177100   2828    364    188  -20  -29    -84  -28   20 
   +   178388  -1932    420     76   53   44    -28  -48  -24 
   +   221375   4095    735     15 -160  -25    -49  -25   87 
   +   236004    -92   -188   -108   51  -30     20  -20   -4 
   +   239085  -2835   -147     45 -108    0   -147   45   45 
   +   253000   2120   -440    200 -125   10    104  -40   24 
   +   284625  -3855   1505    -55  -90   45    273 -115    5 
   +   354200  -4872    840    -32  -40  -58    -56   60   36 
   +   385825   3457   1105    233   -5  -23    -15    5  -19 
   +   442750   -770   1470   -130 -185   40   -210   30  -66 
   +   462000   5040    560   -400   30  120    112   80   16 
   +   664125  -1155    -35   -275  195  -30    637   85    5 
   +   664125  -1155    -35    365  195  -30    -35    5  -11 
   +   664125   2205   1645   -315  195  -30     77  -15  -39 
   +   853875   7155    435    -45  135    0   -237  -45   51 
   +  1044912   -720    496    -80  -15  -87   -160   16  -80 
   +  1288000  -2240  -2240    320  100  100    448    0    0 
   +  1771000 -12040   1400   -200  205  -20   -168   40   40 
   +  1771000   1400   -840    -40 -200  115   -168   40  -56 
   +  1992375   3255  -1225    215  180   45   -441  -25  -25 
   +  2004750   8910  -1170   -450    0    0   -162  -90   54 
   +  2095875  -3645  -2205    -45    0    0   -189   75   27 
\end{verbatim}
\newpage
\noindent\begin{verbatim}
   @    @    @    @    @    @    @    @    @    @    @    @    @ 
   4                                                             
9152 6144 6144 1280 5760 5184 4320 3456  576  288   56  768  768 
   A    B    B    C   AB   AA   BA   BA   BB   BC    A    A    C 
   A    A    A    A   AB   AA   BA   BA   BB   BC    A    A    A 
  4D   4E   4F   4G   6A   6B   6C   6D   6E   6F   7A   8A   8B 
   1    1    1    1    1    1    1    1    1    1    1    1    1 
  -1    3   -1   -1    4    0    3   -3    1   -1    2   -1   -3 
   5    1   -3    1   10    2    2    2   -2   -2    1   -3    3 
  -5    7    3   -1    5   -3    6    6    2    2    2    3    5 
  -5   -5    3   -1   21   -3   -6    0    0    2    0    3   -1 
   7    7    7   -1   -2    6    3    9    1    3    0   -1    3 
   6   10   -2    2   -3   -3    3   -9    3   -1    0   -2   -4 
   1    5    1    1   15    3    9   -3   -3    1    0    1   -5 
  17    1   -7   -3   -3    9   -3    3    3   -3    0    1    1 
   1   -3   -7    5    0    4   -5   -5    3    3    0    1    3 
   1   -3   -7    5    0    4   -5   -5    3    3    0    1    3 
  19   -1    3   -5    3   -9    0    0    0    0    0    3    1 
  10   18    2    6    0   -4   -1   11    3   -1    1    2    6 
  24    8    8    0   10    6   15   -3    1    3   -2    8    0 
 -16   12    0   -8    8   12  -18    6    2    2    0    0    2 
  -7   -7    1    5   35   -5   -5   -5   -1   -1   -1   -7    1 
   1    5   25    5   -5   -1    5   11    7    1   -1    1   -1 
 -23    5   -7   -7  -10   -2    7    7   -1   -1   -1    1    3 
 -13   11   -5   -1   -5    3    6  -18   -2   -2    0    3   -5 
   2    2   10  -10   15    7   10    4    0   -2   -2   -6    2 
\end{verbatim} 
 
\pagebreak 
 
\noindent\begin{verbatim} 
 -51    9   -3    5    0    0    0    0    0    0   -1   -3    3 
 -51    9   -3    5    0    0    0    0    0    0   -1   -3    3 
 -17    3   -1    7    0   12    3  -15   -3   -1    0   -1   -7 
 -17    3    7    3  -12    0   -9   -9    3    3    0   -1   -3 
 -21    7  -21    3    5  -19    2   14    2   -2    0    3    5 
  -4  -12   -4   12    4  -16   -1  -13   -5   -1    0   -4   -8 
  12  -24   12    4   -3  -15   12   12  -12    4    0    4    2 
  -1   -9    7   -5    0    0  -15   15    3   -3    0   -1   -9 
  36   12   12   -8   -5    7    4   10   -2    0   -1   -4    8 
 -19   -3   -3    5   12    0    0    0    0    0    0   -3   -3 
  40    8  -24    0   -5   23  -10    2   -2    2   -1    0    4 
  -7    5    1    5  -10    6   15   -3    5   -1   -2   -7    3 
  32  -12  -16   -8    0  -12  -12    6    6    4    0    0   10 
   9  -19   25   -7   -5   19    7  -23    1   -1   -1    1   -5 
  46   14    6   10   15   -5  -20   -8    0   -4    0    6    2 
  48  -16  -16    0  -10   18    0    0   -4   -4    0    0    0 
  13    5   -3    5   -5   15    0   -6   -2    4    0   -3   -7 
  45    5   -3    5   -5   15    0   -6   -2   -4    0   -3    5 
  37   -7  -11    5   -5   -9    0   18   10    0    0   -3   -5 
 -77  -13    3   -5   15   27    0    0    0    0    1    3   -1 
 -32    0   -8   -4    1    9    9    9    1    1    1    8    4 
 -64    0    0    0  -20   -8  -20    4    4   -4    0    0    0 
 -40   -8   24    0    5   -7  -10   -4   -4   -2    0    0    4 
  -8  -24    8    0    0   -4    5  -13    3    5    0   -8    0 
  39   -9  -17   -5   20  -12   15   -3    5   -1    0   -1    3 
 -18    6    6   10    0    0    0    0    0    0   -1    6    6 
  51   27    3   -5    0    0    0    0    0    0   -2    3   -9 
\end{verbatim}
\newpage
\noindent\begin{verbatim}
   @    @    @    @    @    @    @    @    @    @    @    @    @ 
                                                                 
 512  512  256   64   54   11  864  288  288  288   96   96   48 
   D    D    C    E    A    A   BA   AC   DA   EC   BD   EB   AE 
   A    A    A    A    A    A   AA   AC   BA   BC   AD   BB   AE 
  8C   8D   8E   8F   9A  11A  12A  12B  12C  12D  12E  12F  12G 
   1    1    1    1    1    1    1    1    1    1    1    1    1 
   3   -1    1    1    2    1   -2    0    1   -3    2    1    0 
   5    1   -1   -1    1    0    2   -2    2    4    2    0   -2 
   3   -1    1    1    2    0    1    1   -2    4    1    0    1 
   7   -1   -1   -1   -2    0    1    1    4   -2    1   -2    1 
  -1   -1    3   -1   -2   -1   -2    2    1   -1   -2   -1   -2 
  -2    2    0    0   -2   -1    5    1   -1    1   -3   -3    1 
   5    1   -1   -1    2   -1   -3   -1   -3   -7    1    1   -1 
  -3    5    1    1    0    0   -5    1   -5    1   -1    1    1 
   1   -3   -1   -1    1    0    2   -4   -1   -1   -2   -1    0 
   1   -3   -1   -1    1    0    2   -4   -1   -1   -2   -1    0 
  -1    3   -3    1    0    0   -3    3    0    0    1    0   -1 
  -2   -2   -2    2   -1    0   -6   -4    3   -1   -2   -1    0 
   0    0    0    0    2   -1    4    2    1    5    0    1    2 
  -4    0   -2    2    2    1   -2    0   -2    0    2    4    0 
   5   -3    1    1   -1    0   -1   -1   -1    5   -1    1   -1 
   1    5    3   -1   -1    0    5   -1   -1   -1    1    3   -1 
  -3    1   -1   -1   -1    0   -2    2   -5   -1   -2   -1    2 
  -1   -1    3   -1   -2    0   -1   -1    2    2   -1    2   -1 
  -2    6    2    2    1    0    7   -1    4   -4   -1    0   -1 
\end{verbatim} 
 
\pagebreak 
 
\noindent\begin{verbatim}
   1    5   -1   -1    0    1    0    0    0    0    0    0    0 
   1    5   -1   -1    0    1    0    0    0    0    0    0    0 
  -5    7   -3   -3   -2    0    2    0    5   -3   -2    1    0 
  -9   -5    1    1    0    0   10    0   -5   -3   -2    1    0 
  -9   -5    1    1   -1    1   -3    1    6    4   -3    0    1 
   4   -4    0    0    1    0    6   -4    3   -1    2   -1    0 
   8   -4   -2   -2    2    1   -1   -3    8    0    3    0   -3 
   3    3   -1   -1    2    0   -4    0   -1    3   -4   -1    0 
   0    0    0   -4   -3   -1   -7   -1    2    2   -3   -2    3 
  -3   -3   -3    1    0    0    6    0    0    0    2    0    0 
   0    0    4    0    1    0    5    3    2    0    1   -4   -1 
   1   -3   -1   -1    0    0   -6    2   -3   -1    2   -1    2 
   4    0   -2    2    2    0   -2    0   -2    0    2    0    0 
   5    1   -1   -1    1    0    3   -1   -3   -1    3   -1   -1 
  -6    2    2   -2    1    0   -3    3    0    0    1    0   -1 
   0    0    0    0    0    0    4   -2    4   -2    0    2    2 
  -3    5    1    1    0    0  -11   -1   -2    2    1   -2   -1 
   5   -3   -3    1    0    0    1    7   -2   -2   -3    2   -1 
   1   -3   -1   -1    0    0    5    3    2    0    1    0   -1 
  -5   -5   -1   -1    0    0   -3    3    0    0    1    0   -1 
  12    4    4    0    0    0   -7    1    5    1    1    1   -3 
   0    0    0    0    1   -1   -2    0    4    0    2    0    0 
   0    0    4    0    1    0    3    1    0   -2   -1   -2    1 
   0    0    0    0    1    0   -6    4    3    1   -2    1    0 
  -5    3    3   -1    0    0    0   -4   -3   -1    0   -1    0 
   2    2   -2    2    0    0    0    0    0    0    0    0    0 
   3   -5   -1   -1    0    1    0    0    0    0    0    0    0 
\end{verbatim}
\newpage
\noindent\begin{verbatim}
   @    @    @    @    @    @    @    @    @    @    @    @ 
                                                            
  48   56   28   28   32   32   18   23   23   24   24   28 
  EF   AA   AB   AB    D    C   AB    A    A   CA   BB   AA 
  BF   AA   AB   AB    A    A   AA    A    A   BA   AB   AA 
 12H  14A  14B  C**  16A  16B  18A  23A  B**  24A  24B  28A 
   1    1    1    1    1    1    1    1    1    1    1    1 
  -1   -2    0    0    1   -1    0    0    0   -1    0    0 
   0    1   -1   -1    1   -1   -1    0    0    0    0    1 
   0    2    0    0   -1    1    0   -1   -1    0   -1   -2 
   0    0    0    0    1    1    0    0    0    0   -1    0 
   1    0   -2   -2   -1   -1    0   -1   -1   -1    0    0 
   1    0    0    0   -2    0    0    0    0    1   -1    0 
   1    0    0    0   -1    1    0    0    0    1    1    0 
  -1    0    0    0   -1   -1    0    0    0    1    1    0 
  -1    0    0    0    1   -1    1  b23   **    1    0    0 
  -1    0    0    0    1   -1    1   **  b23    1    0    0 
   0    0    0    0    1   -1    0   -1   -1    0    1    0 
  -1    1    1    1    0    0   -1    0    0   -1    0    1 
  -1   -2    0    0    0    0    0    0    0   -1    0    2 
   0    0    0    0    2    0    0    0    0    0    2    0 
   1   -1    1    1   -1   -1    1    0    0   -1    1   -1 
   1   -1    1    1    1   -1   -1    0    0    1   -1   -1 
  -1   -1    1    1   -1    1    1    1    1    1    0   -1 
  -2    0    0    0    1    1    0    0    0    0    1    0 
  -2    2    0    0    0    0    1    0    0    0   -1    0 
\end{verbatim} 
 
\pagebreak 
 
\noindent\begin{verbatim}
   0   -1   i7  -i7   -1    1    0   -1   -1    0    0    1 
   0   -1  -i7   i7   -1    1    0   -1   -1    0    0    1 
  -1    0    0    0    1   -1    0    1    1   -1    2    0 
   1    0    0    0   -1    1    0    0    0   -1    0    0 
   0    0    0    0    1   -1   -1    0    0    0   -1    0 
  -1    0    0    0    0    0   -1    0    0   -1   -2    0 
   0    0    0    0    2    0    0    0    0   -2   -1    0 
   1    0    0    0    1    1    0    0    0   -1    0    0 
   0   -1    1    1   -2    2    1    1    1    2   -1   -1 
   0    0    0    0    1    1    0    0    0    0    0    0 
   0   -1    1    1    0    0   -1    0    0    0    1   -1 
   1    2    0    0    1   -1    0    0    0   -1    0    0 
   2    0    0    0   -2    0    0    0    0    0   -2    0 
   1   -1   -1   -1   -1    1    1    0    0    1    1   -1 
   0    0    0    0    0    0    1    0    0    0   -1    0 
   2    0    0    0    0    0    0   -1   -1    0    0    0 
   0    0    0    0    1    1    0    0    0    0   -1    0 
   0    0    0    0    1    1    0    0    0    0   -1    0 
  -2    0    0    0   -1    1    0    0    0    0    1    0 
   0    1    1    1   -1   -1    0    0    0    0   -1    1 
   1    1   -1   -1    2   -2    0   -1   -1   -1    1    1 
   0    0    0    0    0    0    1    0    0    0    0    0 
   0    0    0    0    0    0   -1    0    0    0    1    0 
  -1    0    0    0    0    0   -1    0    0    1    0    0 
   1    0    0    0    1    1    0    0    0   -1    0    0 
   0   -1   -1   -1    0    0    0    1    1    0    0   -1 
   0    2    0    0   -1   -1    0    0    0    0    0    0 
\end{verbatim}

\clearpage
\thispagestyle{empty}

\addcontentsline{toc}{chapter}{

\clearpage
\thispagestyle{empty}

\cleardoublepage

\addcontentsline{toc}{chapter}{\textbf{Index}}

\begin{theindex}
\item {atom} {27}
      \subitem {Brauer--- } {27}
      \subitem {projective--- } {27}
      \subitem {system of atoms} {27}
\indexspace
\item {basic set}
      \subitem {of Brauer characters} {24}
      \subitem {of projective characters} {24}
      \subitem {special--- } {36}
\item {block} {14}
      \subitem {strongly real--- } {90}
\item {Brauer graph} {14}
\indexspace
\item {character}
      \subitem {bit of a--- } {65}
      \subitem {Brauer--- } {10}
      \subitem {generalized--- } {11}
      \subitem {part of a--- } {56}
      \subitem {projective--- } {12}
      \subitem {projective indecomposable--- } {12}
      \subitem {subsum of a--- } {57}
      \subitem {symmetrized--- } {20}
      \subitem {virtual--- } {11}
\item {class function} {11}
\indexspace
\item {decomposition matrix} {11}
      \subitem {wedge shape of--- } {39}
\indexspace
\item {essential subset} {34}
\indexspace
\item {lattice} {91}
\indexspace
\item {maximal multiplicity} {65}
\item {multiplicity free} {65}
\indexspace
\item {orthogonality relations} {13}
\indexspace
\item {$p'$-element} {10}
\item {$p$-modular system} {10}
\item {$p$-regular element} {10}
\item {$p$-singular element} {10}
\indexspace
\item {relation} {35}
      \subitem {with respect to a basis} {35}
\indexspace
\item {Weyl module} {20}
\end{theindex}

\clearpage
\thispagestyle{empty}
\cleardoublepage

\thispagestyle{empty}

\address{Gerhard Hiss
        IWR
        Universit\"at Heidelberg
        Im Neuenheimer Feld 368
        69120 Heidelberg, Germany}
\address{Christoph Jansen
        Lehrstuhl D f\"ur Mathematik
        RWTH Aachen
        Templergraben 64
        52062 Aachen, Germany}
\address{Klaus Lux
        Lehrstuhl D f\"ur Mathematik
        RWTH Aachen
        Templergraben 64
        52062 Aachen, Germany}
\address{Richard Parker
        Perihelion Software Ltd.
        Shepton Mallet BA4 4PP
        England}

\end{document}